\documentclass{article}
\usepackage[english]{babel}
\usepackage{csquotes}

\usepackage[letterpaper,top=2cm,bottom=2cm,left=2.5cm,right=2.5cm,marginparwidth=1.75cm]{geometry}
\usepackage[backend=biber, style=numeric,  maxnames=99, sorting=nyt]{biblatex}
\AtEveryBibitem{
    \clearfield{url}
    \clearfield{doi}
    \clearfield{issn}
    \clearfield{eprint}
}

\usepackage{calrsfs}
\usepackage{amssymb}
\usepackage{amsmath}
\usepackage{amsthm}
\usepackage{subcaption}
\usepackage{multirow}
\usepackage{graphicx}
\graphicspath{{images/}}
\usepackage{amsfonts}
\usepackage{mathtools}
\usepackage{mathrsfs}
\usepackage{placeins}

\newtheorem{remark}{Remark}
\newtheorem{definition}{Definition}
\newtheorem{assumption}{Assumption}
\newtheorem{theorem}{Theorem}

\usepackage[mathscr]{euscript}
\usepackage[dvipsnames]{xcolor}
\usepackage[colorlinks=true, allcolors=OliveGreen]{hyperref}
\usepackage{tikz}
\usetikzlibrary{external}
\usepackage{authblk}
\usepackage{todonotes}
\usepackage{array}
\newcolumntype{C}{>{\centering\arraybackslash}p{2.5cm}}

\DeclareMathAlphabet{\mathcalligra}{T1}{calligra}{m}{n}

\newcommand{\avgb}[1]{\{\!\!\{#1\}\!\!\}}
\newcommand{\jumpb}[1]{[\![#1]\!]}
\newcommand{\vertii}[2]{{\left\lVert #1 \right\rVert_{#2}}}
\newcommand{\vertiii}[2]{{\left\vert\kern-0.25ex\left\vert\kern-0.25ex\left\vert #1 
        \right\vert\kern-0.25ex\right\vert\kern-0.25ex\right\vert_{#2}}}
\newcommand{\vertiiDG}[1]{{\lVert #1 \rVert_{\text{DG}}}}
\newcommand{\vertiiiDG}[1]{{\left\vert\kern-0.25ex\left\vert\kern-0.25ex\left\vert #1 
        \right\vert\kern-0.25ex\right\vert\kern-0.25ex\right\vert_{\text{DG}}}}

\usepackage{float}

\title{A discontinuous Galerkin method for the three-dimensional heterodimer model with application to prion-like proteins' dynamics}

\author[1]{Paola F. Antonietti\footnote{paola.antonietti@mail.polimi.it}}

\author[1]{Mattia Corti\footnote{mattia.corti@mail.polimi.it}}

\author[1]{Giacomo Lorenzon\footnote{giacomo.lorenzon@mail.polimi.it}}

\affil[1]{MOX-Dipartimento di Matematica, Politecnico di Milano, Piazza Leonardo da Vinci 32, Milan, 20133, Italy}

\begin{document}
    \maketitle
    
    \begin{abstract}
        Neurocognitive disorders, such as Alzheimer's and Parkinson's, have a wide social impact. These proteinopathies involve misfolded proteins accumulating into neurotoxic aggregates. Mathematical and computational models describing the prion-like dynamics offer an analytical basis to study the diseases' evolution and a computational framework for exploring potential therapies.
        This work focuses on the heterodimer model in a three-dimensional setting, a reactive-diffusive system of nonlinear partial differential equations describing the evolution of both healthy and misfolded proteins. We investigate traveling wave solutions and diffusion-driven instabilities as a mechanism of neurotoxic pattern formation. For the considered mathematical model, we propose a space discretization, relying on the Discontinuous Galerkin method on polytopal/polyhedral grids, allowing high-order accuracy and flexible handling of the complicated brain's geometry. Further, we present \textit{a priori} error estimates for the semi-discrete formulation and we perform convergence tests to verify the theoretical results. Finally, we conduct simulations using realistic data on a three-dimensional brain mesh reconstructed from medical images.
    \end{abstract} 
    
    \section{Introduction}
    \label{sec:introduction}
    Neurodegenerative disorders involve inevitable cognitive decline, often diagnosed late after substantial brain deterioration \cite{VILLEMAGNE2013357}. Many of these conditions, classified as proteinopathies, share a common neurodegenerative mechanism \cite{jucker2018}. Despite varied symptoms, they exhibit similar molecular behaviors: misfolded proteins act as templates, causing further protein misfolding, neurotoxic aggregation, and cell loss \cite{STOPSCHINSKI2017323, Prusiner1998}. Diseases such as Alzheimer's, Parkinson's, and Amyotrophic Lateral Sclerosis follow this pattern, with specific proteins showing stereotypical progression \cite{JuckerWalker2013}. In Parkinson's disease, $\alpha$-synuclein lesions first appear in the dorsal motor and anterior olfactory nuclei before spreading \cite{braak2003}. This resembles prion diseases, suggesting a shared molecular mechanism \cite{JuckerWalker2011}, with axonal transport potentially aiding this spread \cite{JuckerWalker2011}.
    \par
    In Parkinson's, $\alpha$-synuclein monomers form insoluble Lewy bodies \cite{spillantini1997alpha}. Targeting misfolding and aggregation could offer therapeutic avenues \cite{JuckerWalker2013, hasegawa2017prion}. However, detecting $\alpha$-synuclein remains a challenge as current methods are post-mortem and PET imaging ligands are lacking \cite{korat2021alpha, smith2023alpha}. Advances in detecting $\alpha$-synucleinopathies, including Parkinson's disease, are limited \cite{smith2023alpha}. Consequently, mathematical and numerical modeling is essential for understanding prion-like disease mechanisms \cite{hist}.
    \par
    Physics-based models have recently begun shedding light on neurodegenerative processes \cite{WEICKENMEIER2019264, hist}, correlating molecular scales with symptoms \cite{WeickenmeierEtAl2018}. The heterogeneity of spatial and temporal scales involved poses a challenge for clinical observations, but mathematical and computational models offer analytical tools for tailoring personalized therapies. Aggregate models often lack accuracy in depicting prion concentration evolution, as spatial domain description is critical for estimating disease progression \cite{mat}. The prion paradigm is essential for designing physics-based models that consistently describe neurodegenerative disorder traits \cite{WEICKENMEIER2019264}. Several models have been proposed to understand prion kinetics, including monomeric \cite{cohen1994structural} and polymeric \cite{jarrett1993seeding} seeding models. This work focuses on a simpler kinetic model describing the conversion from native to misfolded proteins at a fixed rate \cite{WEICKENMEIER2019264}. Despite the difficulty in quantifying diffusion and transport mechanisms \textit{in vivo} \cite{SCHAFER2019369}, the heterodimer model offers a quantitative framework for prion-like propagation \cite{WEICKENMEIER2019264, corti}. Frequently used models in the literature include the Smoluchowski and Fisher-Kolmogorov models \cite{fornari2019prion}. The Smoluchowski model, while well-suited to describe polymeric fragmentation and aggregation, is complex, prompting the use of simpler models with high spatial detail for intricate domains like the brain. The Fisher-Kolmogorov model, a well-known non-linear reaction-diffusion equation, effectively describes population dynamics as an extension of the logistic equation with a diffusive term. Numerous analytical \cite{Murray2002, Murray2003} and numerical \cite{WEICKENMEIER2019264, cortiAntonietti2023, Corti2023} studies are available. A historical overview of mathematical models for prion-like dynamics is provided in \cite{hist}, with a comprehensive survey in \cite{fornari2019prion}.
    \par
    This work focuses on the heterodimer model \cite{Prusiner1998}, aiming to model molecular interactions and describe macroscopic prion-like propagation across complex spatial domains such as the human brain. In particular, we extend and validate the work of \cite{corti} in a three-dimensional, patient-specific setting, demonstrating that our model not only replicates but also refines these findings. Indeed, our approach considers both native and misfolded protein fields, simplifying the Smoluchowski approach while enhancing the infection mechanism description compared to the Fisher-Kolmogorov model. The heterodimer model explicitly includes physical coefficients representing molecular kinetic processes (production, destruction, and conversion) \cite{fornari2019prion}, crucial for developing effective therapeutic interventions. The bistability of the dynamical system justifies seeking traveling wave-like solutions, for which a closed expression in linear approximation can be obtained, and a velocity estimate derived, generalizing considerations in \cite{Murray2003}. Moreover, we validate our results thanks to Braak staging theory and we delve further into the sensitivity analysis results found in the literature.
    \par
    For the numerical approximation of the heterodimer model, finite element methods are commonly used \cite{fornari2019prion, mat}, alongside reduced-order network diffusion models \cite{thompson2020protein, mat}. Recently, discontinuous Galerkin (DG) formulations of the heterodimer model have been proposed on polytopal/polyhedral grids \cite{corti}. These partial differential equations typically admit wavefront propagating solutions \cite{Murray2002}, necessitating high-order approximating schemes capable of handling complex geometries with locally varying discretization parameters. This justifies the use of a polytopal discontinuous (PolyDG) approach \cite{Antonietti2021}, employing arbitrarily shaped elements and enabling mesh construction from clinical images through agglomeration techniques \cite{antonietti2022agglomeration}. For time integration, the $\theta$-method scheme is adopted, with the nonlinear term treated semi-implicitly. Stability results and {a priori} convergence error estimates are extensively studied by \cite{corti}.
    \par
    \bigskip
    This work is organized as follows. In Section~\ref{sec:model}, we present the heterodimer model and we discuss some theoretical results. In Section~\ref{sec:discretization}, we derive the polytopal discontinuous Galerkin formulation and we present theoretical results concerning stability and {a priori} error estimates. Also, we discuss the fully discrete formulation. In Section~\ref{sec:numerics}, we present numerical results verifying the theoretical convergence rates, and in Section~\ref{sec:results} we deal with a realistic simulation of $\alpha$-synuclein spreading in a three-dimensional domain obtained from magnetic resonance imaging and tractography data. In Section~\ref{sec:conclusions}, we conclude and discuss future developments.
    
    \section{Mathematical modelling of protein spreading}
    \label{sec:model}
    The mathematical modeling of the prion-like pathogenesis kinetics establishes a coherent connection between different spatial and temporal scales, including microscopic biochemical reactions and the accumulation of misfolded proteins across extensive brain regions, causing neuronal demise followed by the disease's symptomatic expressions. Whilst the prion conversion mechanism remains poorly understood, a straightforward approach to unravel the interactions between native and prionic proteins is the heterodimer model, initially proposed by Prusiner \cite{prusiner1991molecular, prusiner1982novel}. According to this model, a prion molecule interacts with a healthy prion protein, forming a heterodimer. Subsequently, the prion protein undergoes conversion into its misfolded variant, and the resultant polymer fragments into two singular misfolded prions that perpetuate as infectious seeds, propagating the infection \cite{prionmech}. The heterodimer model offers one of the simplest descriptions of the proteins' kinetic in their native state $c$ and misfolded one $q$. As suggested in \cite{fornari2019prion}, the conversion rate between them is globally described by a coefficient $k_{12}$, representing the single-step reaction
    \begin{equation*}
        c + q \xrightarrow[]{k_{12}} q + q,
    \end{equation*}
    in which the rate $k_{12}$ is proportional both to the concentrations of native and misfolded proteins \cite{WEICKENMEIER2019264}.
    
    \subsection{Strong formulation}
    The dynamics of the aggregate concentration of healthy and misfolded proteins within a bounded, open domain $\Omega \subset \mathbb{R}^N, \ N=2,3$, with Lipschitz boundary, over a time span $(0, T], \ T \in \mathbb{R}^+$ is described by a system of non-linear partial differential equations (\ref{eq: heterodimer}) denoted by heterodimer model \cite{fornari2019prion}. It reads:
    \begin{equation}
        \label{eq: heterodimer}
        \left\{
        \begin{aligned}
            \frac{\partial c}{\partial t} &= \nabla \cdot (\mathbf{D}\nabla c) - k_1 c - k_{12} cq + k_0, \qquad &\text{in } \Omega \times (0,T],\\
            \frac{\partial q}{\partial t} &= \nabla \cdot (\mathbf{D}\nabla q) - \widetilde{k}_1 q + k_{12} cq, \qquad &\text{in } \Omega \times (0,T].
        \end{aligned}
        \right.
    \end{equation}
    The aggregate concentration of native proteins and those that are misfolded are represented here by the functions $c,q: \Omega \times (0, T] \xrightarrow{} \mathbb{R}^+$, respectively. A comprehensive mathematical model is obtained by imposing the set of natural border and initial conditions. Coefficients $k_0, k_1, \widetilde{k}_1, k_{12}$ and data are all non-negative and sufficiently regular. In particular, the diffusion tensor $\mathbf{D}$ is chosen uniformly elliptic.
    
    In specific brain regions, the modeled protein is highly expressed at a rate of $k_0 \geq 0$. Clearance can occur through various biological pathways. This elimination mechanism is described by the clearance rate coefficient $k_1 \geq 0$. The overall conversion from healthy to misfolded proteins happens at a rate of $k_{12} \geq 0$, as previously described. A similar clearance process also interests the misfolded proteins. Notice that, being stereoisomers of the healthy ones, their three-dimensional conformation is different, and so is their chemical functionality. In particular, their solubility is strongly affected. The clearance coefficient $\widetilde{k}_1 \geq 0$ is thus in general different from $k_1$. The diffusion term in (\ref{eq: heterodimer}) describes the nature of the spreading for both healthy and misfolded proteins from the macroscopic perspective. We take into account two contributions, unlike in \cite{mat}: extracellular and axonal diffusion. The first one is concentration gradient-driven and it is supposed isotropic, typically, it is associated with neuronal proximity \cite{frost2009propagation}. The second mechanism of spread is anisotropic, occurring along anatomically defined paths called axons \cite{seeley2009neurodegenerative}. Let $\boldsymbol{a}: \Omega \subset \mathbb{R}^N \xrightarrow{} \mathbb{R}^N : |\boldsymbol{a}(\boldsymbol{x})|=1 \quad \forall \boldsymbol{x} \in \Omega$ be the local axonal direction at $\boldsymbol{x}$ in $\Omega$, with $|\cdot|$ being the Euclidean norm in $\mathbb{R}^{N}$. The diffusion tensor is then defined as
    \begin{equation*}
        \label{eq:diffanis}
        \mathbf{D}(\boldsymbol{x}) = d_{\text{ext}}\mathbf{I} + d_{\text{axn}}\boldsymbol{a}(\boldsymbol{x}) \otimes \boldsymbol{a}(\boldsymbol{x}), \quad \boldsymbol{x} \in \Omega,
    \end{equation*}
    with $d_{\text{ext}} > 0$, and with $\mathbf{I}$ being the identity. Notice that being the axonal diffusion faster than the extracellular one, typically $d_{\text{axn}} > d_{\text{ext}} \geq 0$ is assumed \cite{WEICKENMEIER2019264}.
    \par
    The bistability of the system suggests the following assumption.
    \begin{assumption}[Data assumptions\label{ass1}]
        To have physically relevant solutions with equilibria lying in the first quadrant of the phase space $(c,q)$, the inequality $k_0 k_{12} > k_1 \widetilde{k}_1 > 0$ must be satisfied.
    \end{assumption}
    \begin{remark}
        Under Assumption \ref{ass1}, no bifurcations can be observed in system (\ref{eq: heterodimer}). Namely, no qualitative changes in the system's dynamics emerge as its parameters are varied.
    \end{remark}
    We observe that, under Assumption \ref{ass1}, the equilibrium point $(\widetilde{k}_1/k_{12}, k_0/\widetilde{k}_1 - k_1/k_{12})$ is always stable. Under the same condition, $(k_0/k_1, 0)$ is instead always unstable and, in particular, it is a saddle point. This is coherent with the disease's phenomenological point of view, since from the time of the onset, the disease always evolves inevitably toward a pathological condition. This motivates the search for a solution in the form of a traveling wave. Also, we obtain an explicit expression for the minimum velocity of the traveling wave solution that reads:
    \begin{equation}
        \label{eq: minveltrue}
        \nu_{\text{min}} = 2\sqrt{d_\text{ext} \widetilde{k}_1\left(\frac{k_0 k_{12}}{k_1 \widetilde{k}_1}-1\right)}.
    \end{equation}
    According to \cite{Murray2002}, it can be derived the traveling wavefront's velocity as a function of the initial solution. In particular, if the wave profile has compact support, the solution will travel with minimum speed $\nu_{\text{min}}$ defined as in (\ref{eq: minveltrue}). Further, it can be shown that the microstructure of the domain, and so the diffusion tensor, together with the front shape influence the propagation velocity.
    
    \section{Polytopal discontinuous Galerkin semi-discrete formulation}
    \label{sec:discretization}
    In this section, we present the polytopal discontinuous Galerkin semi-discrete formulation of model (\ref{eq: heterodimer}), along with the stability and a priori error analysis \cite{corti}. We then derive the fully discrete formulation based on the Discontinuous Galerkin method.
    
    \subsection{PolyDG semi-discrete formulation}
    Discontinuous Galerkin methods on polytopal grids are a family of finite element methods for the approximation of partial differential equations that can support arbitrarily shaped mesh elements \cite{Riviere2008, HesthavenWarburton2008, Antonietti2016} and extend the classical paradigm of Discontinuous Galerkin methods \cite{AntoniettiHoustonGiani2010, HoustonCangianiDongGeorgoulis2016, HoustonCangianiDongGeorgoulis2017}. The idea is to look for the numerical solution in a discrete space of polynomials, possibly discontinuous across mesh elements. We introduce some preliminary results and definitions to derive the PolyDG formulation. Let $\mathcal{T}_h$ be a partition of $\Omega$. Its elements $K \in \mathcal{T}_h$ are a finite number of non-empty and disjoint polyhedral. In what follows, for each of them, $|K|$ is the element's measure, namely the $N$-dimensional Lebesgue measure of $K$, $h_K = \operatorname{diam}(K)$ is the element's diameter, $h = \operatorname{max}\{h_K : K \in \mathcal{T}_h\}$ is the global mesh size parameter.
    
    We will refer to mesh interfaces as the intersection of the $(N-1)$-dimensional facets of two neighboring elements $K,K' \in \mathcal{T}_h$. In dimension $N=3$, each interface consists of a general polygon. By assumption, we suppose that each polygon may be decomposed into a set of co-planar triangles. We will refer to $\mathcal{F}_h$ as the collection of all the triangles of all the sub-triangulations of each interface. We will use the term face to refer to one of the elements in $\mathcal{F}_h$ \cite{Antonietti2021}. Internal faces are denoted as $\mathcal{F}^I_h=\left\{ F \in \mathcal{F}_h : F \subset \Omega \right\}$ and the boundary ones as $\mathcal{F}_h^B = \mathcal{F}_h \setminus \mathcal{F}_h^I $. The latter can be further subdivided into Dirichlet boundary faces $\mathcal{F}_h^D=\left\{ F \in \mathcal{F}_h : F \subset \partial \Omega_D \right\}$ and Neumann boundary faces $\mathcal{F}_h^N=\left\{ F \in \mathcal{F}_h : F \subset \partial \Omega_N \right\}$. Here it is assumed that each face $F$ belongs entirely either to $\mathcal{F}_h^D$ or $\mathcal{F}_h^N$. We also require $\mathcal{T}_h$ to satisfy some regularity conditions, following the work of \cite{hpgalerkin, corti}. This allows for exploiting the trace, the inverse, and the discrete Gagliardo-Nirenberg inequalities, and eventually deriving the upcoming estimates.
    
    \begin{assumption}
        \label{gridass1}
        We assume that the family of meshes $\{\mathcal{T}_h\}_h$ is uniformly polytopic-regular.
    \end{assumption}
    
    The following assumption will allow us to state the inverse trace estimate \cite{hpgalerkin}. For more details, we refer to \cite{Antonietti2021}. Also, we require the next assumption that will be instead useful \cite{hpgalerkin} for the approximation results presented in the next sections.
    \begin{definition}
        A covering $\mathcal{T}_\# = \{T_K\}$ related to the polytopic mesh $\mathcal{T}_h$ is a set of shape-regular $N$-dimensional simplices $T_K$, such that for each $K\in \mathcal{T}_h$, there exists a $T_K \in \mathcal{T}_\#$ such that $K\subsetneq T_K$.
    \end{definition}
    From now on, the binary operator $\lesssim$ from now on denotes the following relation: $\alpha \lesssim \beta$ if and only in there exists $ C=C(m,p)>0: \alpha \leq C \beta$, where $p$ is the degree of the polynomial approximation and $m$ are the model parameters. Notice that $C$ is independent of the mesh size $h$.
    \begin{assumption}
        \label{gridass2}
        There exists a covering $\mathcal{T}_\#$ of $\mathcal{T}_h$, such that
        \begin{equation*}
            \max_{K \in \mathcal{T}_h} \operatorname{card} \{ K' \in \mathcal{T}_h : K' \cap T_K \neq \varnothing, T_K \in \mathcal{T}_\# \ s.t. \ K \subset T_K \} \lesssim 1,
        \end{equation*}
        and $h_{T_K} \lesssim h_K$ for each pair $K \in \mathcal{T}_h$ and $T_K \in \mathcal{T}_\#$ with $K \subset T_K$.
    \end{assumption}
    
    For the subsequent analysis, we require some more definitions and results. We define the broken Sobolev space as follows:
    \begin{equation*}
        W = H^1(\Omega, \mathcal{T}_h) = \{ w \in L^2(\Omega): w|_K \in H^1(\Omega) \quad \forall K \in \mathcal{T}_h \},
    \end{equation*}
    and the DG-broken version:
    \begin{equation*}
        W_h^\text{DG} = \{ w \in L^2(\Omega): w|_K \in \mathbb{P}_p(K) \quad \forall K \in \mathcal{T}_h \},
    \end{equation*}
    where $\mathbb{P}_p(K)$ be the space of polynomials of total degree $p \geq 1$ on $K$. We will also make use of two standard trace-related operators: the jump and average operators, applicable to both scalar and vector-valued functions. These will be referred to as $\jumpb{\cdot}$ and $\avgb{\cdot}$, respectively. For further details, refer to \cite{arnold}.
    
    \subsubsection{PolyDG semi-discrete approximation of the heterodimer model}
    Here we present the semi-discrete formulation of model (\ref{eq: heterodimer}). Let $c_h^0, q_h^0 \in W_h^\text{DG}$ be the $L^2(\Omega)$-projection of the initial conditions $c^0$ and $q^0$ onto $W_h^\text{DG}$, respectively. The semi-discrete problem reads as follows.
    \par
    Given the initial conditions $(c_h^0,q_h^0)$, for each $t \in (0,T]$ find $(c_h(t), q_h(t)) \in W_h^\text{DG} \times W_h^\text{DG}$ so that:
    \begin{equation}
        \label{eq:semidiscr}
        \left\{
        \begin{aligned}
            &\displaystyle\int_\Omega \displaystyle\frac{\partial c_h}{\partial t} w_h + \mathcal{A}(c_h, w_h) + r_L(c_h, w_h) + r_N(c_h, q_h, w_h) = F_c(w_h), \ &\forall w_h \in W_h^\text{DG},\\
            &\displaystyle\int_\Omega \displaystyle\frac{\partial q_h}{\partial t} v_h + \mathcal{A}(q_h, v_h) + \widetilde{r}_L(q_h, v_h) - r_N(q_h, c_h, v_h) = F_q(v_h), \ &\forall v_h \in W_h^\text{DG}.
        \end{aligned}
        \right.
    \end{equation}
    The bilinear and trilinear forms are defined, for a fixed triangulation $\mathcal{T}_h$, for any $u_h, \varphi_h, \phi_h \in W_h^\text{DG}$ in the following way:
    \begin{equation*}
        \label{eq:bilin}
        \begin{aligned}
            \mathcal{A}(u_h, \varphi_h) =& \int_\Omega (\mathbf{D} \nabla_h u_h)\cdot \nabla_h \varphi_h
            - \displaystyle\sum_{F \in \mathcal{F}_h^I \cup \mathcal{F}_h^D} \int_F \jumpb{u_h}\cdot\avgb{\mathbf{D} \nabla_h \varphi_h}\\
            &+ \displaystyle\sum_{F \in \mathcal{F}_h^I \cup \mathcal{F}_h^D} \int_F \gamma_F \jumpb{u_h}\cdot\jumpb{\varphi_h}
            - \jumpb{\varphi_h} \cdot \avgb{\mathbf{D} \nabla_h u_h},
        \end{aligned}
    \end{equation*}
    \begin{equation*}
        r_L(u_h,\varphi_h) = \int_\Omega k_1 u_h \varphi_h,
        \qquad
        \widetilde{r}_L(u_h,\varphi_h) = \int_\Omega \widetilde{k}_1 u_h \varphi_h,
        \qquad
        r_N(u_h, \varphi_h, \phi_h) = \int_\Omega k_{12} u_h\varphi_h\phi_h,
    \end{equation*}
    \begin{equation*}
        \begin{aligned}
            F_c(\varphi_h) &= \int_\Omega f_c \varphi_h,
            \qquad
            F_q(\varphi_h) &= \int_\Omega f_q \varphi_h,
        \end{aligned}
    \end{equation*}
    where $\nabla_h$ is the element-wise gradient, and $\gamma_F$ is the discontinuity penalization function defined as
    \begin{equation*}
        \label{eq:penltycoeff}
        \gamma_F = \gamma_0 
        \begin{dcases}
            \max\left\{\{d^K\}_\mathrm{H},\{k^K\}_\mathrm{H}\right\}\dfrac{p^2}{\{h\}_\mathrm{H}}, & \mathrm{on} \; F \in \mathcal{F}_h^I, \\
            \max\left\{d^K, k^K\right\}\dfrac{p^2}{h},                               & \mathrm{on} \; F \in \mathcal{F}_h^D.
        \end{dcases}
    \end{equation*}
    The coefficient $\gamma_0$ is a constant parameter that should be chosen sufficiently large to ensure the stability of the discrete formulation, $d^K = \|\sqrt{\mathbf{D}|_K}\|^2$, and $k^K = \|\,(1 + k_{12}|_K)(k_1|_K + \widetilde{k}_1|_K)\|$. Finally, $\{\cdot\}_H$ is the harmonic average operator defined as $\{v\}_\mathrm{H} = (2 v^+ v^-)/(v^+ + v^-)$.
    
    \subsection{Semi-discrete PolyDG formulation's analysis}
    \label{subsection:semi_analysis}
    We recall some useful definitions to prove the problem's well-posedness, the error estimates, and the scheme stability. For sufficiently regular functions, we can define:
    \begin{equation}
        \label{eq:dg2}
        \vertiiDG{w} = \vertii{\sqrt{\mathbf{D}}\nabla_h w}{} + \vertii{\gamma_F^{1/2} \jumpb{w}}{\mathcal{F}^I_h\cup\mathcal{F}^D_h},
    \end{equation}
    \begin{equation*}
        \label{eq:dg3}
        \vertiiiDG{w} = \vertiiDG{w} + \vertii{\gamma_F^{-1/2}\avgb{\mathbf{D}\nabla_h w}}{\mathcal{F}^I_h\cup\mathcal{F}^D_h},
    \end{equation*}
    where $\nabla_hv$ is the element-wise gradient, namely $\nabla_h v |_K = \nabla v|_K$. These two norms are well defined and also equivalent in the space of discontinuous functions $W_h^{\text{DG}}$, provided the regularity Assumptions \ref{gridass1} and \ref{gridass2} on the mesh  \cite{corti}. Two more norms will be useful in our analysis:
    \begin{equation*}
        \label{eq:eps2}
        \vertii{w(t)}{\epsilon}^2 = \vertii{w(t)}{}^2 + \int_0^t \vertiiDG{w(s)}^2 \mathrm{d}s,
    \end{equation*}
    \begin{equation*}
        \label{eq:eps3}
        \vertiii{w(t)}{\epsilon}^2 = \vertii{w(t)}{}^2 + \int_0^t \vertiiiDG{w(s)}^2 \mathrm{d}s,
    \end{equation*}
    for $t>0$. Under the same regularity Assumptions \ref{gridass1} and \ref{gridass2}, the bilinear form $\mathcal{A}$ is coercive \cite{hpgalerkin} with respect to DG-norm (\ref{eq:dg2}) and continuous with respect to (\ref{eq:dg3}), provided that the penalty parameter is large enough.
    
    The stability estimate for the DG-discrete solution is proven in \cite{corti}. Here we omit the proof and just recall that the solution in the energy norm $\vertii{\cdot}{\epsilon}$ is bounded by data, provided that $\gamma_0$ is large enough, and under the regularity Assumptions \ref{gridass1} and \ref{gridass2} on the mesh. The detailed analysis of the error estimates can be found in \cite{corti}. For our purpose, the following result is sufficient.\vspace{0.5em}
    \begin{theorem}[Error bounds \cite{corti}]
        Let $\mathcal{T}_h$ be the partition of $\Omega$, so that it induces a quasi-uniform grid. If the exact solution $c(t), q(t) \in H^{p+1}(\Omega) \ \forall t \in (0, T]$, with $p$ polynomial degree of approximation, then
        \begin{equation*}
            \label{eq:conv_h}
            \vertiii{c(t) - c_h(t)}{\epsilon}^2 + \vertiii{q(t) - q_h(t)}{\epsilon}^2 \lesssim h^p,
        \end{equation*}
        provided that the penalty parameter $\gamma_0$ is large enough.
    \end{theorem}
    
    In this section, we present the algebraic formulation of model (\ref{eq: heterodimer}), and the fully discrete formulation.
    \subsubsection{Algebraic form}
    Let $N_h = \operatorname{dim} W_h^\text{DG}$ be the dimension of the DG space. Fix a basis $\{ \varphi_j \}_{j=1}^{N_h}$. Let $\boldsymbol{C},\boldsymbol{Q}: \mathcal{T}_h \times (0,T] \xrightarrow[]{} \mathbb{R}^{N_h}$ be the vectors of all the coefficients of the solution expansion in $W_h^{DG}$. More precisely, the solution expansion with respect to the chosen basis reads as
    \begin{equation*}
        c_h(\boldsymbol{x},t)=\sum^{N_h}_{j=1}C_j(t)\varphi_j(\boldsymbol{x}),
        \qquad q_h(\boldsymbol{x},t)=\sum^{N_h}_{j=1}Q_j(t)\varphi_j(\boldsymbol{x}).
    \end{equation*}
    The algebraic form of problem (\ref{eq:semidiscr}) follows.
    \begin{equation}
        \label{eq:algebraic}
        \left\{
        \begin{aligned}
            &\mathbf{M} \frac{\partial \boldsymbol{C}(t)}{\partial t} + \mathbf{A} \boldsymbol{C}(t) + \mathbf{R}_L \boldsymbol{C}(t) + \mathbf{R}_N\big(\boldsymbol{C}(t)\big) \boldsymbol{Q}(t)&= \boldsymbol{F}_c(t),
            \qquad t \in (0,T],\\
            &\mathbf{M} \frac{\partial \boldsymbol{Q}(t)}{\partial t} + \mathbf{A} \boldsymbol{Q}(t) + \widetilde{\mathbf{R}}_L \boldsymbol{Q}(t) - \mathbf{R}_N\big(\boldsymbol{Q}(t)\big) \boldsymbol{C}(t)&= \boldsymbol{F}_q(t),
            \qquad t \in (0,T],\\
            &\boldsymbol{C}(0) = \boldsymbol{C}^0, \quad \boldsymbol{Q}(0)=\boldsymbol{Q}^0,
        \end{aligned}
        \right.
    \end{equation}
    where $\boldsymbol{C}^0$ and $\boldsymbol{Q}^0$ are the expansion coefficients in the chosen basis of the projected initial data. In system (\ref{eq:algebraic}), $\mathbf{M}, \mathbf{A}, \mathbf{R}_L, \widetilde{\mathbf{R}}_L, \mathbf{R}_N  \in \mathbb{R}^{{N_h}\times{N_h}}$, whereas $\boldsymbol{F}_c, \boldsymbol{F}_q \in \mathbb{R}^{N_h}$. The component-wise definitions are given by:
    \begin{equation*}
        \mathrm{M}_{ij} = (\varphi_i, \varphi_j),
        \qquad
        \mathrm{A}_{ij} = \mathcal{A}(\varphi_i, \varphi_j),
        \qquad
        \mathrm{R}_{L,ij} = r_L(\varphi_i, \varphi_j),
        \qquad
        \widetilde{\mathrm{R}}_{L,ij} = \widetilde{r}_L(\varphi_i, \varphi_j),
    \end{equation*}
    \begin{equation*}
        \mathrm{R_N}(\psi(t))_{ij} = r_N(\psi(t), \varphi_i, \varphi_j),
        \qquad
        \mathrm{F}_{c,i} = F_c(\varphi_i),
        \qquad
        \mathrm{F}_{q,i} = F_q(\varphi_i),
    \end{equation*}
    with $i,j$ ranging in $\{1, \dots, N_h\}$.
    
    \subsubsection{Time discretization}
    Let $\{ t_i \}_{i=0}^{N_T}$ be the partition of the time interval $[0,T]$, so that $0=t_0<t_1<\dots<t_{N_T}=T$ and $\Delta t = t_{k+1} - t_k, \quad \forall k \in \{0,\dots, N_T-1\}$. Denote with $(\boldsymbol{C}^0, \boldsymbol{Q}^0) = (\boldsymbol{C}(0), \boldsymbol{Q}(0))$ the initial datum and define $(\boldsymbol{C}^k, \boldsymbol{Q}^k) \approx (\boldsymbol{C}(t_k), \boldsymbol{Q}(t_k)), \quad \forall k \in \{1,\dots, N_T\}$. Now we discretize the time derivative with a classic forward scheme and introduce the theta-method scheme \cite{dgaflio} to integrate in $t$.
    
    Given the initial datum $(\boldsymbol{C}^0, \boldsymbol{Q}^0)$, find $(\boldsymbol{C}^k, \boldsymbol{Q}^k), \ \forall k=1,2,\dots, N_T$ solving system (\ref{eq:matrix_theta}):
    \begin{equation}
        \label{eq:matrix_theta}
        \left\{
        \begin{aligned}
            &\left(\frac{\mathbf{M}}{\Delta t} + \theta \mathbf{A} + \theta \mathbf{R}_L \right)\boldsymbol{C}^k + \mathbf{R}_{c,N}^\theta\boldsymbol{Q}^{k} = \boldsymbol{F}_c^\theta,\\
            &\left(\frac{\mathbf{M}}{\Delta t} + \theta \mathbf{A} + \theta \widetilde{\mathbf{R}}_L \right)\boldsymbol{Q}^k - \mathbf{R}_{q,N}^\theta\boldsymbol{C}^{k} = \boldsymbol{F}_q^\theta,
        \end{aligned}
        \right.
    \end{equation}
    where the forcing term is defined and the non-linear term is linearized according to the chosen method.
    \begin{itemize}
        \item Implicit Euler method $\theta = 1$:
        \begin{equation*}
            \mathbf{R}_{c,N}^{\theta =1} = \mathbf{R}_N(\boldsymbol{C}^{k-1}) = \boldsymbol{C}^{k-1},
            \qquad
            \mathbf{R}_{q,N}^{\theta =1} = \mathbf{R}_N(\boldsymbol{Q}^{k-1}) = \boldsymbol{Q}^{k-1},
        \end{equation*}
        \begin{equation*}
                \boldsymbol{F}_c^{\theta =1} = \boldsymbol{F}_c(t_k) + \frac{\mathbf{M}}{\Delta t}{\boldsymbol{C}^{k-1}},
                \qquad
                \boldsymbol{F}_q^{\theta =1} = \boldsymbol{F}_q(t_k) + \frac{\mathbf{M}}{\Delta t}{\boldsymbol{Q}^{k-1}}.
        \end{equation*}
        \item Crank-Nicolson method $\theta = 1/2$:
        \begin{equation*}
            \begin{aligned}
                &\mathbf{R}_{c,N}^{\theta=1/2} = \frac{3}{4} \mathbf{R}_N(\boldsymbol{C}^{k-1}) - \frac{1}{4} \mathbf{R}_N(\boldsymbol{C}^{k-2}) = \frac{3 \boldsymbol{C}^{k-1} - \boldsymbol{C}^{k-2}}{4},\\
                &\mathbf{R}_{q,N}^{\theta=1/2} = \frac{3}{4} \mathbf{R}_N(\boldsymbol{Q}^{k-1}) - \frac{1}{4} \mathbf{R}_N(\boldsymbol{Q}^{k-2}) = \frac{3\boldsymbol{Q}^{k-1} - \boldsymbol{Q}^{k-2}}{4} ,
            \end{aligned}
        \end{equation*}
        \begin{equation*}
            \label{eq:force_theta}
            \begin{aligned}
                &\boldsymbol{F}_c^{\theta=1/2} =
                 \frac{\boldsymbol{F}_c(t_k) + \boldsymbol{F}_c(t_{k-1})}{2}
                -
                \frac{3\boldsymbol{C}^{k-1}-\boldsymbol{C}^{k-2}}{4}\boldsymbol{Q}^{k-1}
                +
                 \left( \frac{\mathbf{M}}{\Delta t} -\frac{\mathbf{A} + {\mathbf{R}_L}}{2} \right){\boldsymbol{C}^{k-1}},\\
                &\boldsymbol{F}_q^{\theta=1/2} =
                \frac{\boldsymbol{F}_q(t_k) + \boldsymbol{F}_q(t_{k-1})}{2}
                -
                \frac{3\boldsymbol{Q}^{k-1}-\boldsymbol{Q}^{k-2}}{4}\boldsymbol{C}^{k-1}
                +
                \left(\frac{\mathbf{M}}{\Delta t}
                -\frac{\mathbf{A} + \widetilde{\mathbf{R}}_L}{2}
                \right)
                {\boldsymbol{Q}^{k-1}}.
            \end{aligned}
        \end{equation*}
    \end{itemize}
    Notice that the Crank-Nicolson scheme is defined for \(k \geq 2\). Thus, the first iteration should be performed with the Euler Implicit method nonetheless.
  
    \section{Three-dimensional numerical results: verification}
    \label{sec:numerics}
    This section focuses on validating the numerical results against theoretical predictions. Firstly, we have reproduced the results of \cite{corti} in two dimensions. Here we just discuss the three-dimensional generalization.
    
    \subsection{Convergence rates verification}
    \label{subsec:convergence_rates_verification}
    For the first test case, we set $\Omega = (0,1)^N, \ N=2,3$. Diffusion is set isotropic and the final time is set sufficiently small. Since the nuanced interplay between the diverse contributions of the reaction term and the diffusion term may introduce heightened numerical challenges, here we discuss a realistic test case. The detailed setting of the simulations follows in Table \ref{tab:phys_param_re}.
    \begin{table}[ht]
        \centering
        \begin{tabular}{ccc|ccc}
            \textbf{Coefficient} & \textbf{Value} & \textbf{Dimension} & \textbf{Coefficient} & \textbf{Value} & \textbf{Dimension} \\ \hline
            $d_{\text{axn}}$ & 0.0 & $[\mathrm{mm}^2 \mathrm{years}^{-1}]$ & $k_0$ & 0.6 & $ [\mathrm{\mu g}^{1}\mathrm{years}^{-1}\mathrm{mm}^{-3}]$ \\ \hline
            $d_{\text{ext}}$ & 8.0 & $[\mathrm{mm}^2 \mathrm{years}^{-1}]$ & $k_1$ & 0.5 & $ [\mathrm{years}^{-1}]$ \\ \hline
            $k_{12}$ & 1.0 & $[\mathrm{\mu g}^{-1}\mathrm{mm}^3\mathrm{years}^{-2}]$ & $\widetilde{k}_1$ & 0.3 & $ [\mathrm{years}^{-1}]$
        \end{tabular}
        \caption{Test case of Section \ref{subsubsection:realistic_coefficients_trivial_domain}. Physical parameters, taken from \cite{corti}.}
        \label{tab:phys_param_re}
    \end{table}
    From the numerical viewpoint refer to Table \ref{tab:num_param_re_dom_re}.
    \begin{table}[ht]
        \centering
        \begin{tabular}{lll} \hline
            \textbf{Temporal discretization}& Time step& $\Delta t =5 \cdot 10^{-6}$\\ 
            &Final time&  $T = 5 \cdot 10^{-5}$\\ 
            &Theta method&  $\theta = 1/2$\\ \hline 
            \textbf{DG parameters}&Interior Penalty Parameter&  $\eta = 1$\\ 
            &Penalty Parameter&  $\gamma_0 = 10$\\ 
        \end{tabular}
        \caption{Test case of Section \ref{subsubsection:realistic_coefficients_trivial_domain}. Numerical parameters.}
        \label{tab:num_param_re_dom_re}
    \end{table}
    For these test cases, we used uniform cubical cell grids. We chose the exact solution as
    \begin{equation*}
        \label{eq:exact}
        \begin{aligned}
            c(\boldsymbol{x},t) &= \cos(t) \sum_{i=1}^{N} \cos(2 \pi x_i),
            \qquad
            q(\boldsymbol{x},t) &= \left(\prod_{i=1}^{N}\cos(6\pi x_i)+2\right)e^{-t}, \qquad &\text{in } \Omega \subset \mathbb{R}^N \times (0,T],
        \end{aligned}
    \end{equation*}
    and set the right-hand side and the boundary data accordingly.
    
    \subsubsection{Realistic coefficients, idealized domain}
    \label{subsubsection:realistic_coefficients_trivial_domain}
    
\begin{figure}[ht]
        \centering
        \begin{subfigure}{0.32\textwidth}
            \includegraphics[width=\linewidth]{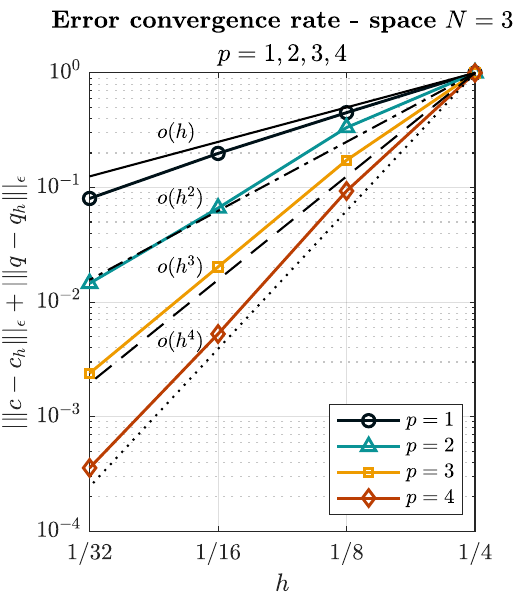}
            \caption{}
            \label{fig:3energy_p_re}
        \end{subfigure}
        \begin{subfigure}{0.32\textwidth}
            \includegraphics[width=\linewidth]{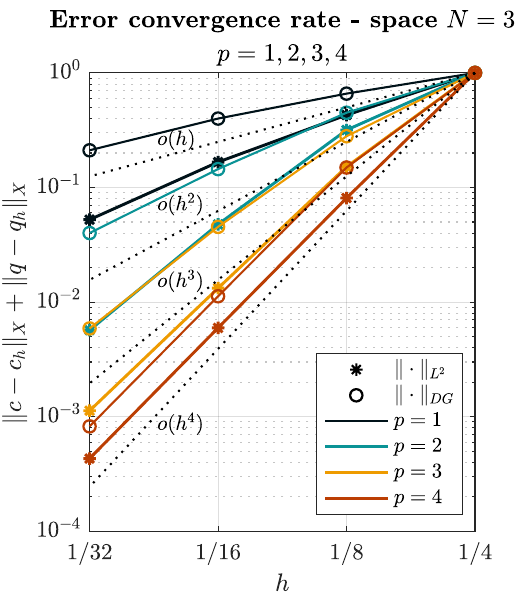}
            \caption{}
            \label{fig:3L2DG_p_re}
        \end{subfigure}
        \begin{subfigure}{0.32\textwidth}
            \includegraphics[width=\linewidth]{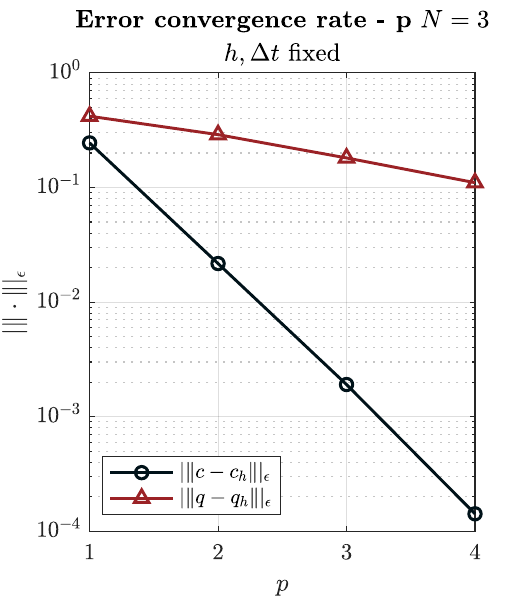}
            \caption{}
            \label{fig:p_3d_re}
        \end{subfigure}
        \caption{Test case of Section \ref{subsubsection:realistic_coefficients_trivial_domain}. (a) Computed errors in the $\vertiii{\cdot}{\epsilon}$ norm versus the mesh size $h$ (loglog scale) for different polynomial approximation degrees $p=1,2,3,4$. (b) Computed errors in the $\|\cdot \|_{L^2(\Omega)}$ and $\|\cdot \|_{DG}$ norms, versus the mesh size $h$ (loglog scale) for different polynomial approximation degrees $p=1,2,3,4$. (c) Computed errors in the $\vertiii{\cdot}{\epsilon}$ norm versus the polynomial degree of approximation (semilog scale) for $h=1/4, \Delta t = 5\cdot 10^{-6}$.}
        \label{fig:err_qua}
    \end{figure}

    Figures \ref{fig:3energy_p_re} and \ref{fig:3L2DG_p_re} show the computed errors in $\vertiii{\cdot}{\epsilon}$, $\| \cdot \|_{L^2(\Omega)}$ and $\| \cdot \|_{DG}$ norms as a function of $h$, for $p=1,2,3,4$, whereas Figure \ref{fig:p_3d_re} shows the computed errors in $\vertiii{\cdot}{\epsilon}$ as a function of $p$, for $h=1/4, \ \Delta t=5\cdot 10^{-6}$. Notably, Figure \ref{fig:3energy_p_re} shows that the theoretical convergence rate is consistently attained, in some cases even outperformed. Both Figures \ref{fig:3energy_p_re}, and \ref{fig:3L2DG_p_re} depict errors that have been normalized to highlight the convergence rate. Also, from Figure \ref{fig:p_3d_re} we observe an exponential convergence for the error. Notice that this case is not covered by the theoretical bounds of Section \ref{subsection:semi_analysis}. The error for the component $q(\boldsymbol{x},t)$ is always greater than that of $c(\boldsymbol{x},t)$ because its exact representation has a much higher frequency of oscillation, and thus it is more challenging to be accurately reproduced numerically.
    
    \subsubsection{Realistic coefficients, realistic domain}
    \label{subsubsection:realistic_coefficients_realistic_domain}
    We perform a convergence test on a realistic domain, based on data taken from Magnetic Resonance Imaging (MRI) within the OASIS-3 database \cite{LaMontagne2019}. The mesh of the brain has been reconstructed as explained in \cite{Corti2023, corti}. Since the mesh is fixed, we only check that the error gets lower for a higher degree of polynomial approximation. Figure \ref{fig:p_3d_brain} shows the computed errors in $\vertiii{\cdot}{\epsilon}$ as a function of $p$, for $\Delta t=5\cdot 10^{-6}$ fixed. An exponential convergence rate is clearly observable. Moreover, one can notice that the value of the error norm $\vertiii{\cdot}{\epsilon}$ is similar for the two solution components $c(\boldsymbol{x},t)$ and $q(\boldsymbol{x},t)$ since they both have the same frequency of oscillation.
    
        \begin{minipage}{0.47\linewidth}
            \centering
            \begin{figure}[H]
                \centering
                \includegraphics[width=0.71111\linewidth]{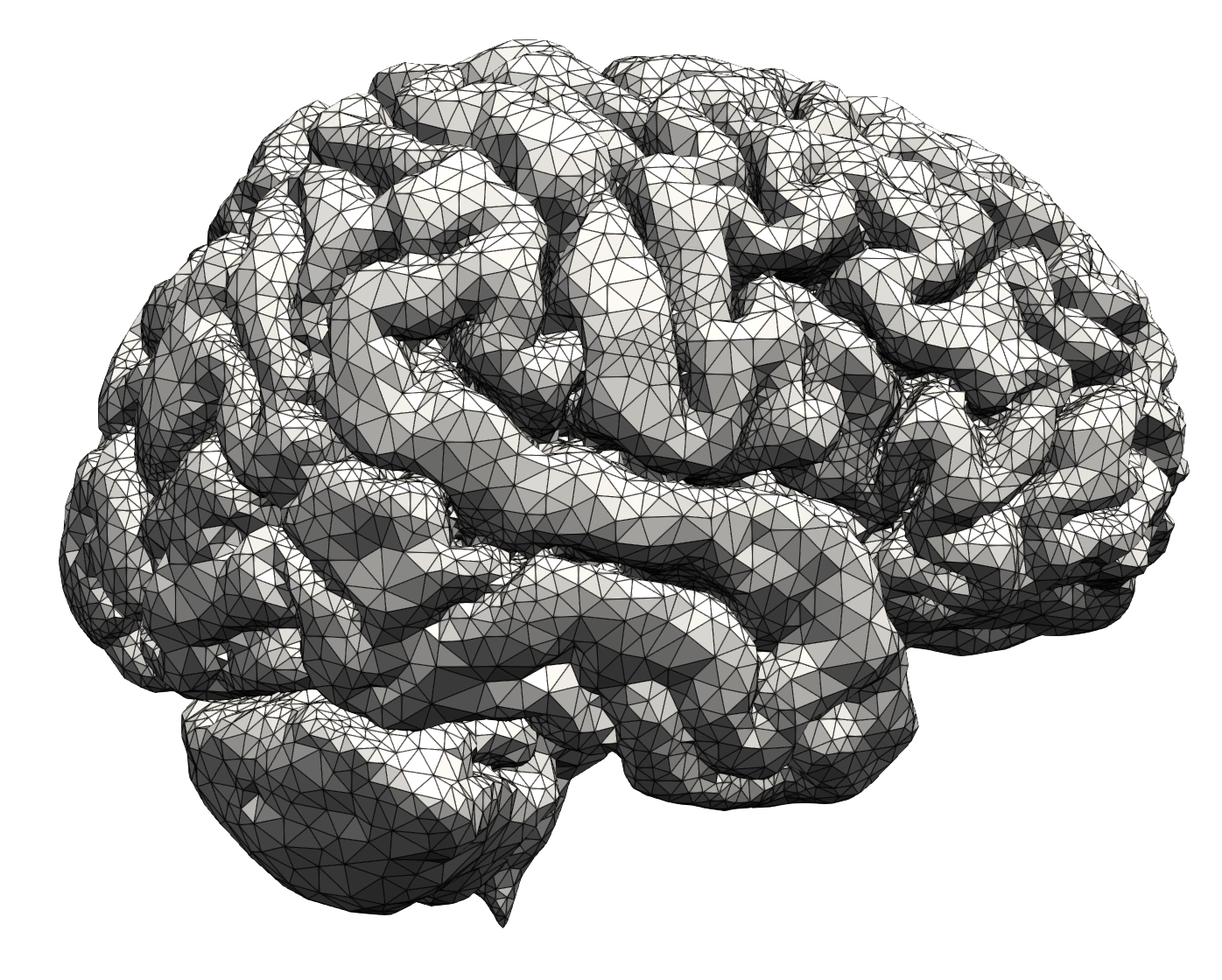}
                \caption{3D domain partition.}
                \label{fig:mesh}
            \end{figure}
        \end{minipage}
        \begin{minipage}{0.47\linewidth}
            \centering
            \begin{figure}[H]
                \centering
                \includegraphics[width=0.71111\linewidth]{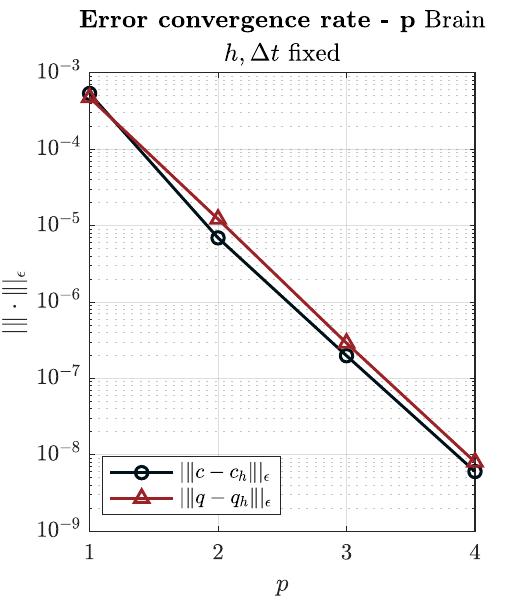}
                \caption{Test case of Section \ref{subsubsection:realistic_coefficients_realistic_domain}. Computed errors in the $\vertiii{\cdot}{\epsilon}$ norm versus the polynomial degree of approximation (semilog scale) for $\Delta t = 5\cdot 10^{-6}$.}
                \label{fig:p_3d_brain}
            \end{figure}
        \end{minipage}
    
    \section{Numerical simulation of \texorpdfstring{$\alpha$}{alfa}-synuclein spreading}
    \label{sec:results}
    In this section, we analyze the numerical results obtained from the simulation of the diffusion of the $\alpha$-synuclein protein \cite{spillantini1997alpha} within the above-mentioned three-dimensional brain mesh. For the details refer to \cite{cortiAntonietti2023, mardal2022mathematical}. Thanks to the high-order DG method used, the triangulation obtained is coarsened, following the work of \cite{antonietti2022agglomeration}. This way, the computational load is kept low. Diffusion-weighted imaging is used to obtain the connectome wiring, as explained in \cite{cortiAntonietti2023}, following the work of \cite{mardal2022mathematical}.
    
    As for the initial condition of the simulation, the seeding region is inferred from post-mortem histological analysis of infected brains \cite{dickson2018neuropathology}, and it is located in the dorsal motor nucleus. Since the patient under analysis is healthy, the initial concentration $c_0(\boldsymbol{x})$ is set equal to the healthy component of the unstable equilibrium $\boldsymbol{e}_2$. Thus, $c_0(\boldsymbol{x}) = k_0 / k_1 = 1.2$. As for the boundary conditions, it is assumed that no flux is allowed across the brain boundary. This translates into null Neumann boundary conditions for both healthy and misfolded proteins, namely $\partial \Omega$ coincides with $\Gamma_N$.  The specific choice of the parameters of model (\ref{eq: heterodimer}) is reported in Table \ref{tab:phys_param_re}. Some studies have tried to measure each parameter of Prusiner's model \cite{
        iljina2016kinetic,
        li2005quantifying,
        masel1999quantifying,
        perrino2019quantitative}, and others have investigated the contribution of each component through sensitivity analyses \cite{SCHAFER2019369, WEICKENMEIER2019264}. All these considerations lead to the choices shown in Table \ref{tab:phys_param_re}.
    In particular, the ratio $d_\text{axn}/d_\text{ext}=10$ is assumed as a modelling choice, following \cite{SCHAFER2019369}. As \cite{fornari2019prion} shows in Section 2.2, the heterodimer model can be reduced to the Fisher-Kolmogorov model, with a single partial differential equation. This can be done assuming the healthy protein concentration constant, and thus focusing solely on the description of the misfolded protein field evolution. The resulting model comprises a diffusive and a reaction term, multiplied by a constant, say $\alpha$, which can be expressed in terms of the coefficients of the heterodimer model:
    \begin{equation*}
        \alpha = k_{12} \frac{k_0}{k_1} - \widetilde{k}_1.
    \end{equation*}
    Through a sensitivity analysis, \cite{SCHAFER2019369, WEICKENMEIER2019264} argue that a high value of $\alpha$ effectively describes the underlying phenomenon, and thus set $\alpha=0.9/\mathrm{years}$. The reactive coefficients are thus not independent. The work of \cite{corti, cortiAntonietti2023} independently confirms the appropriateness of this choice.\vspace{0.5em}
    
    From the numerical perspective, we employ the MUltifrontal Massively Parallel sparse direct Solver \cite{AMESTOY2000501} (MUMPS). This solver is specifically tailored for addressing large sparse systems to solve algebraic equations on distributed memory parallel computers. The time step is fixed at $\Delta t = 1/120 \ \mathrm{years}$, while we set the final simulation time at $T=30 \ \mathrm{years}$. We employ a Crank-Nicolson time integration scheme and we fix the polynomial order of approximation to $p=2$.
    
\begin{figure}[ht]
        \centering
        \includegraphics[width=0.3\textwidth ]{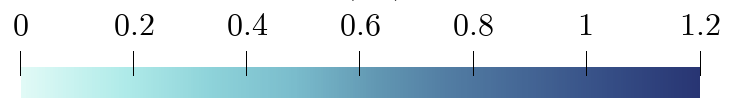}\\
        
        \includegraphics[width=0.14\textwidth ]{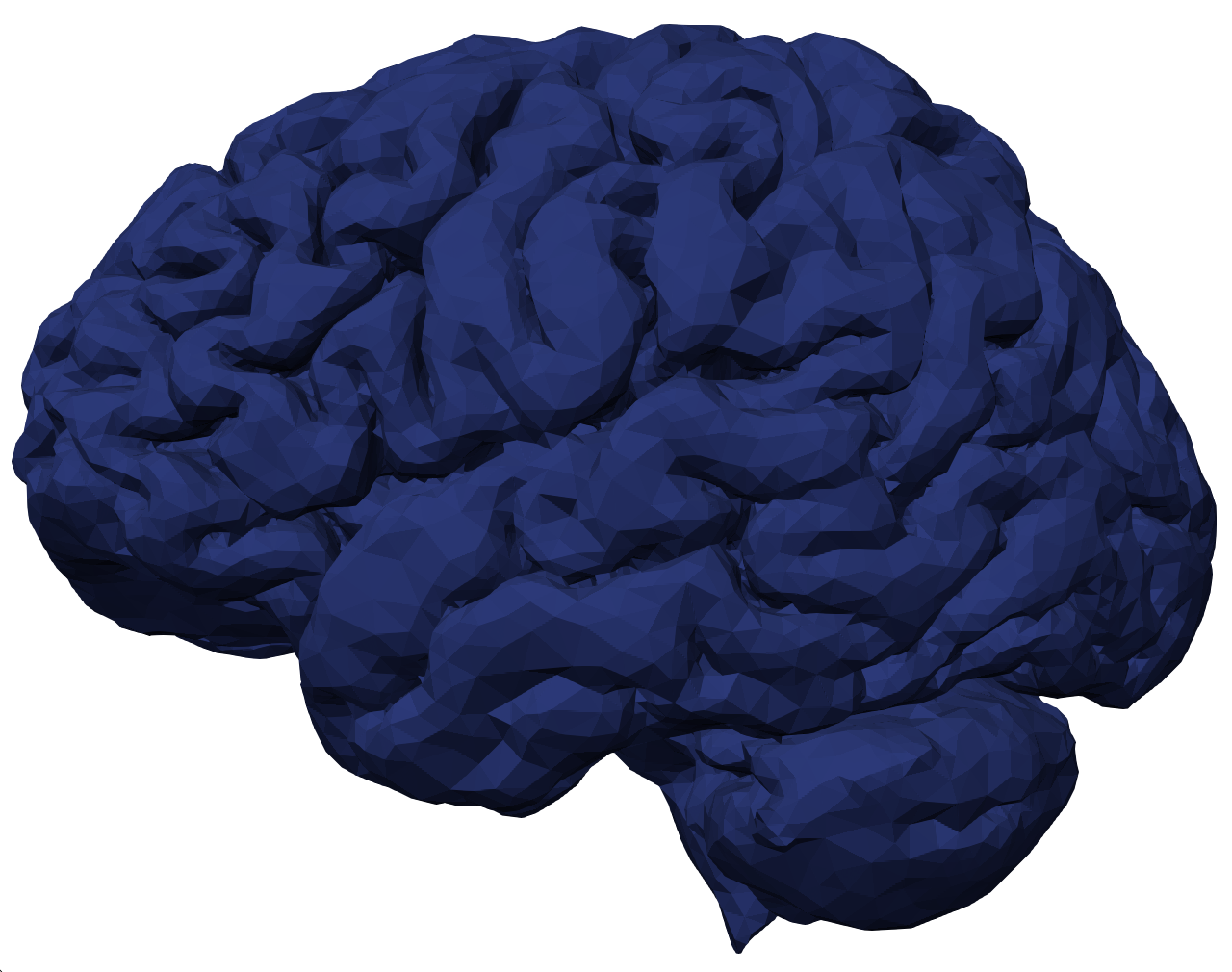}
        \includegraphics[width=0.14\textwidth ]{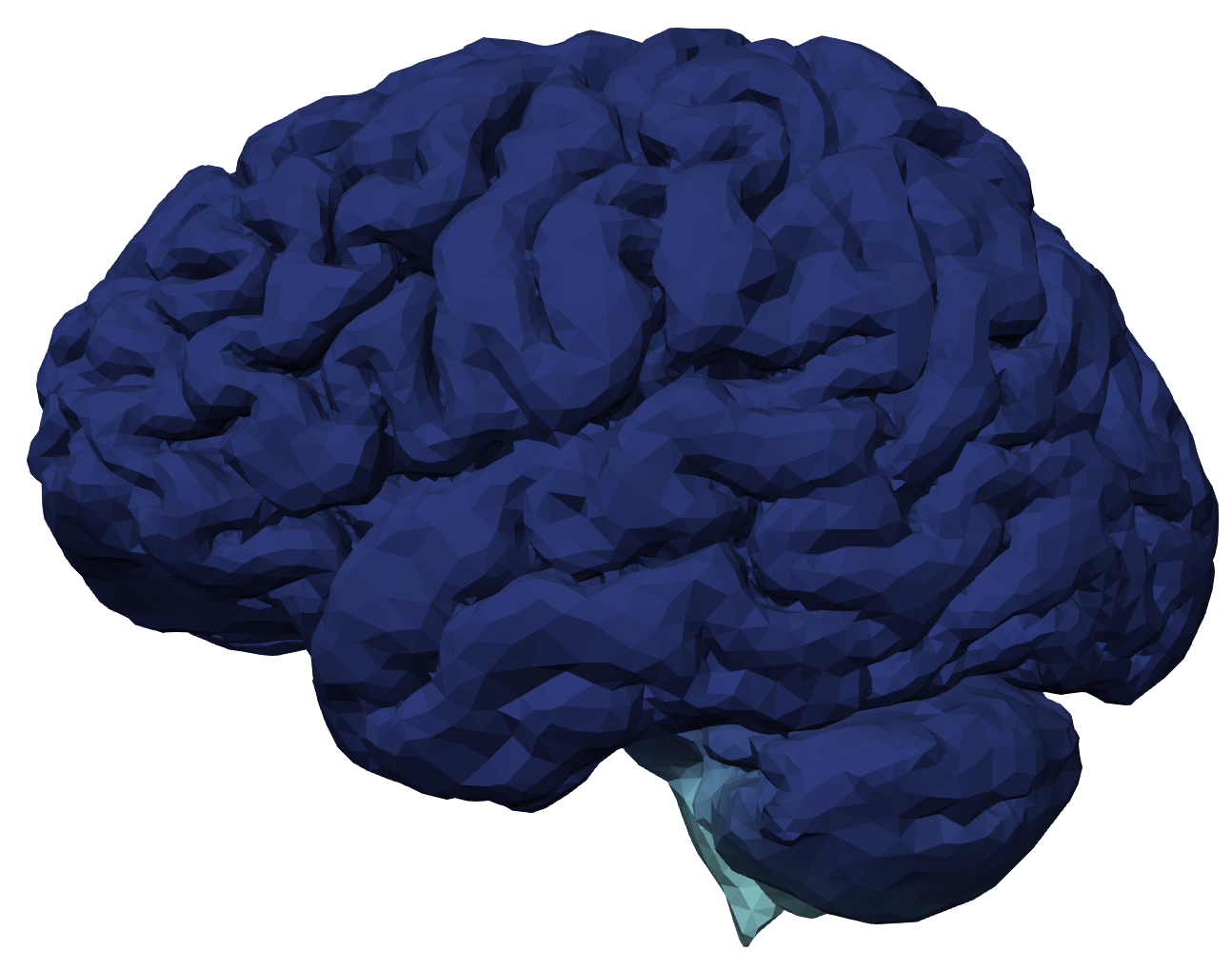}
        \includegraphics[width=0.14\textwidth ]{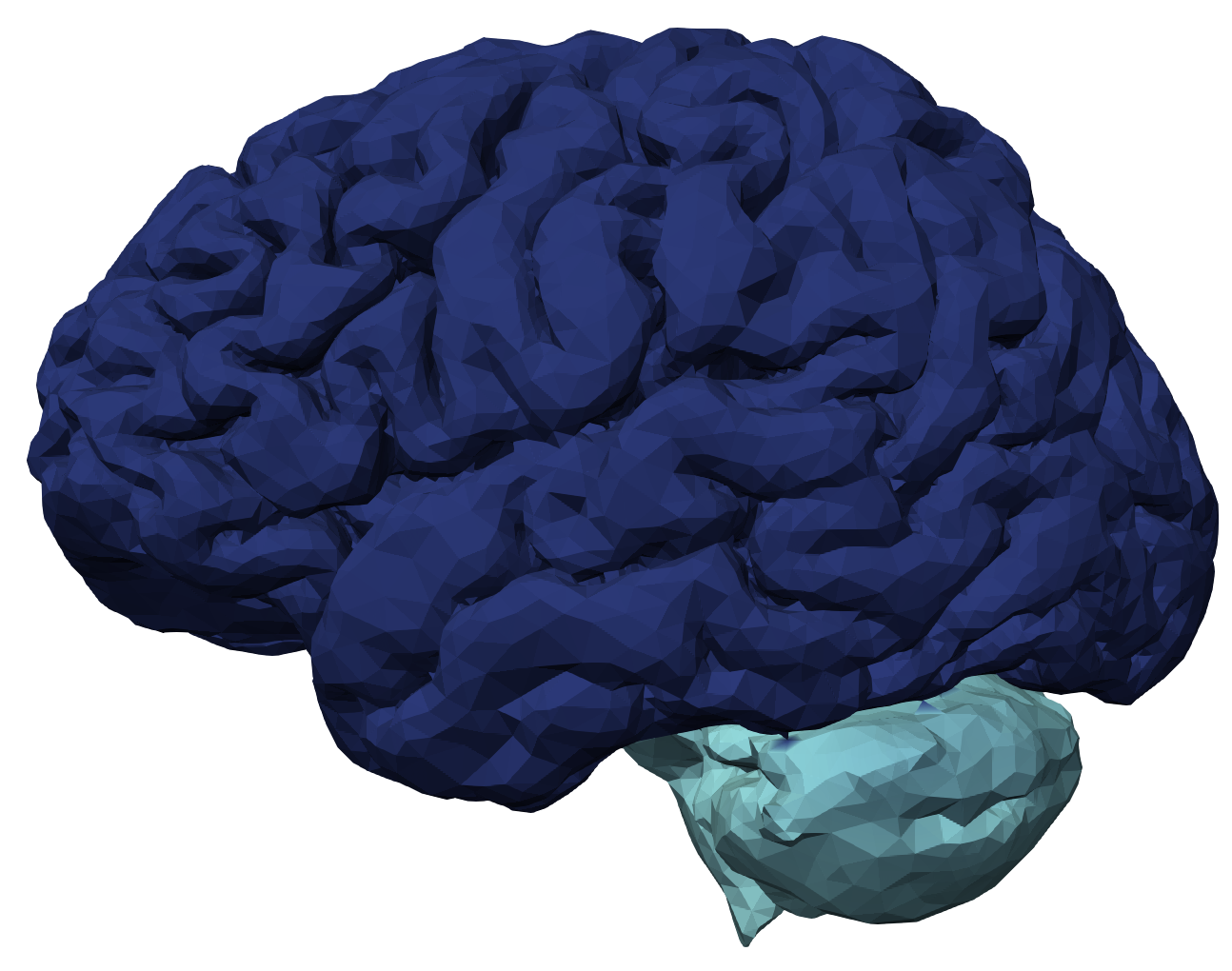}
        \includegraphics[width=0.14\textwidth ]{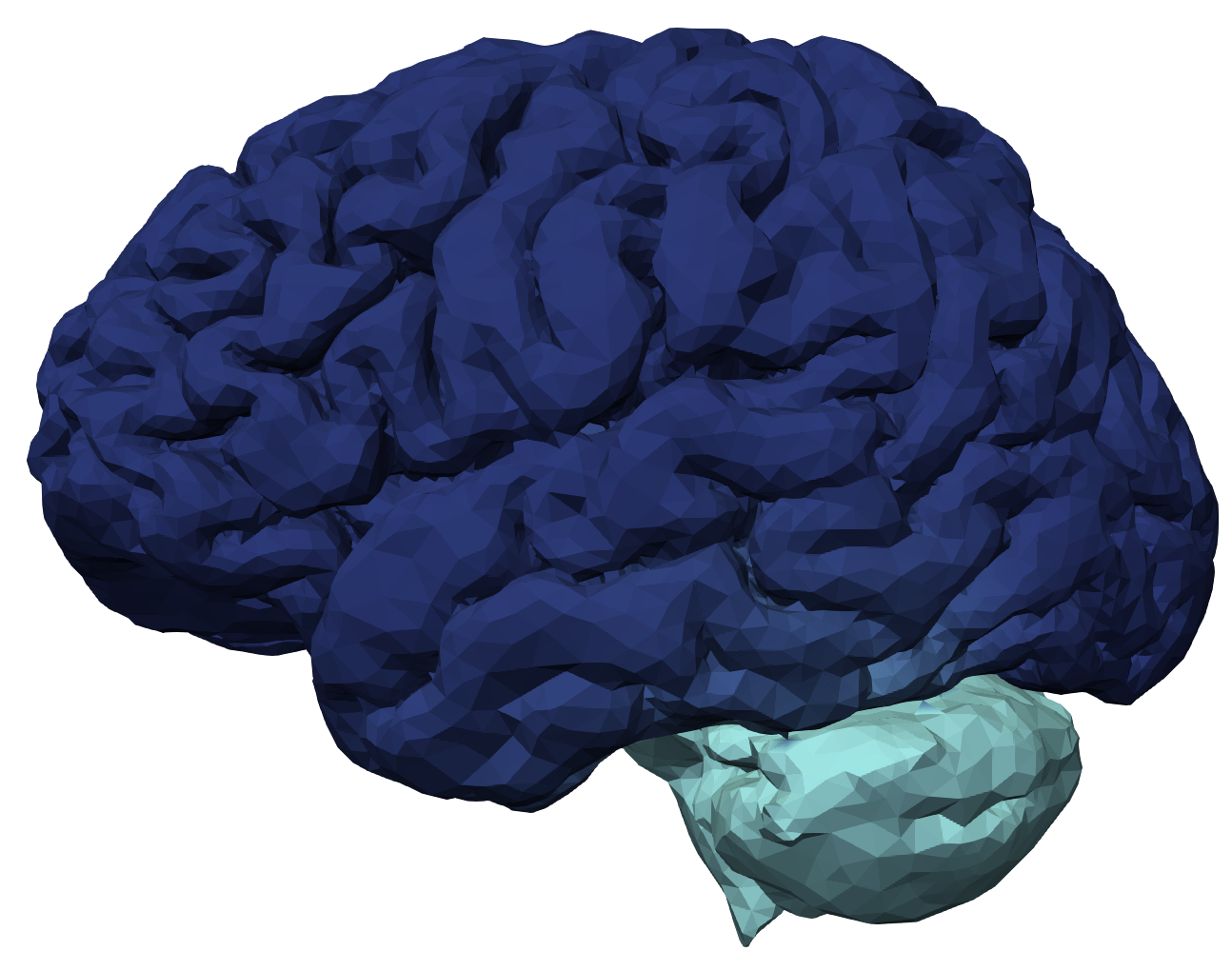}
        \includegraphics[width=0.14\textwidth ]{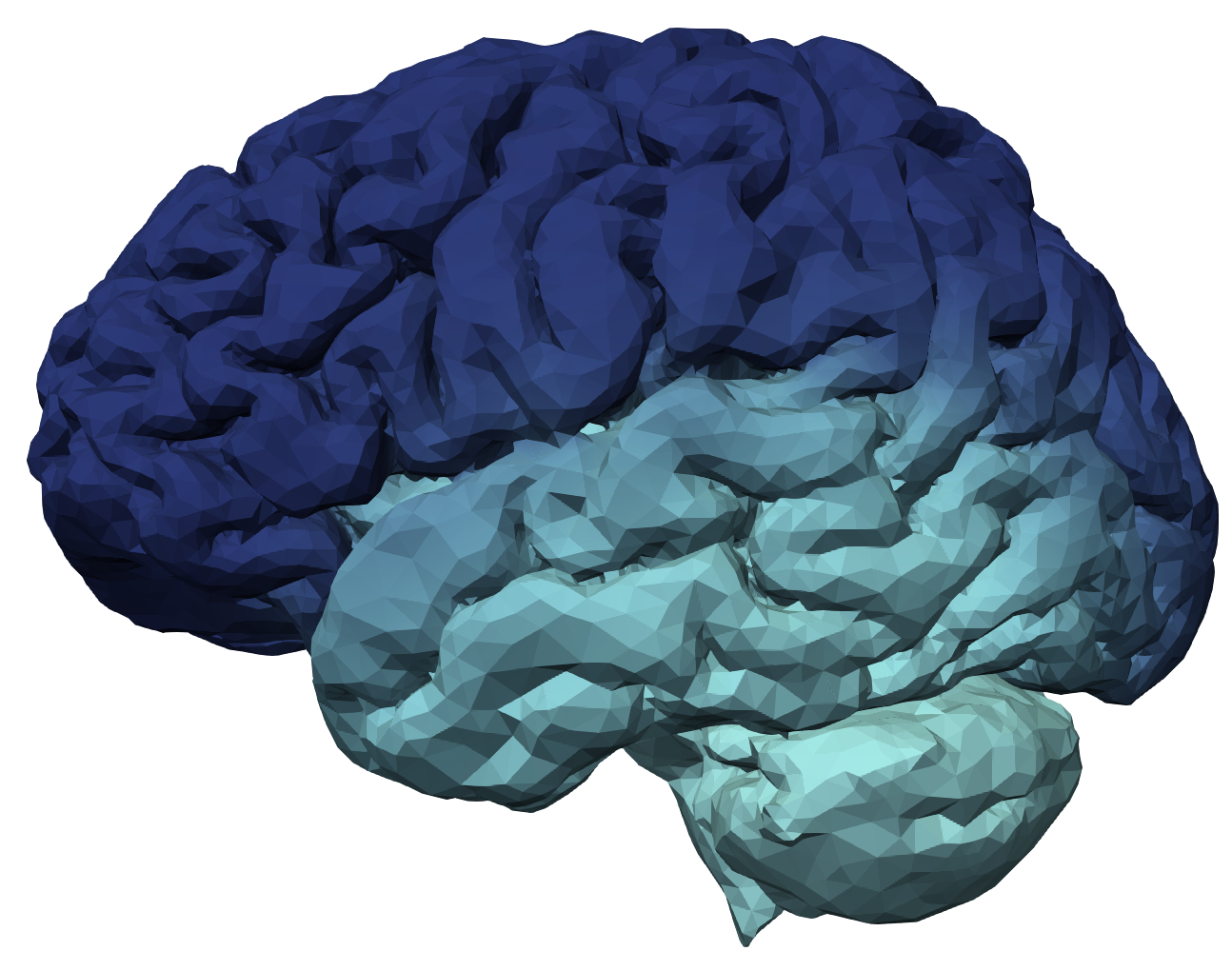}
        \includegraphics[width=0.14\textwidth ]{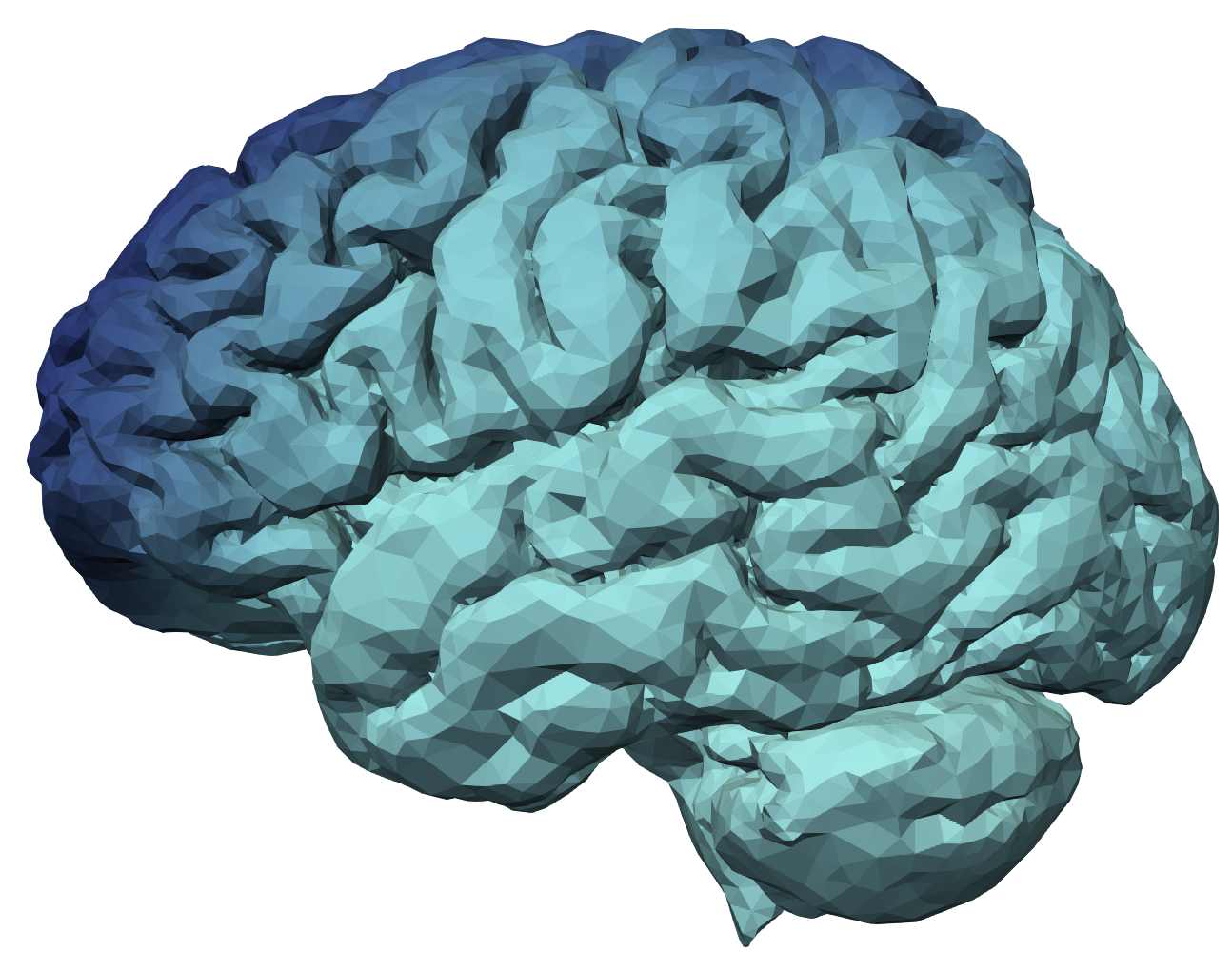}\\
        
        \includegraphics[width=0.14\textwidth ]{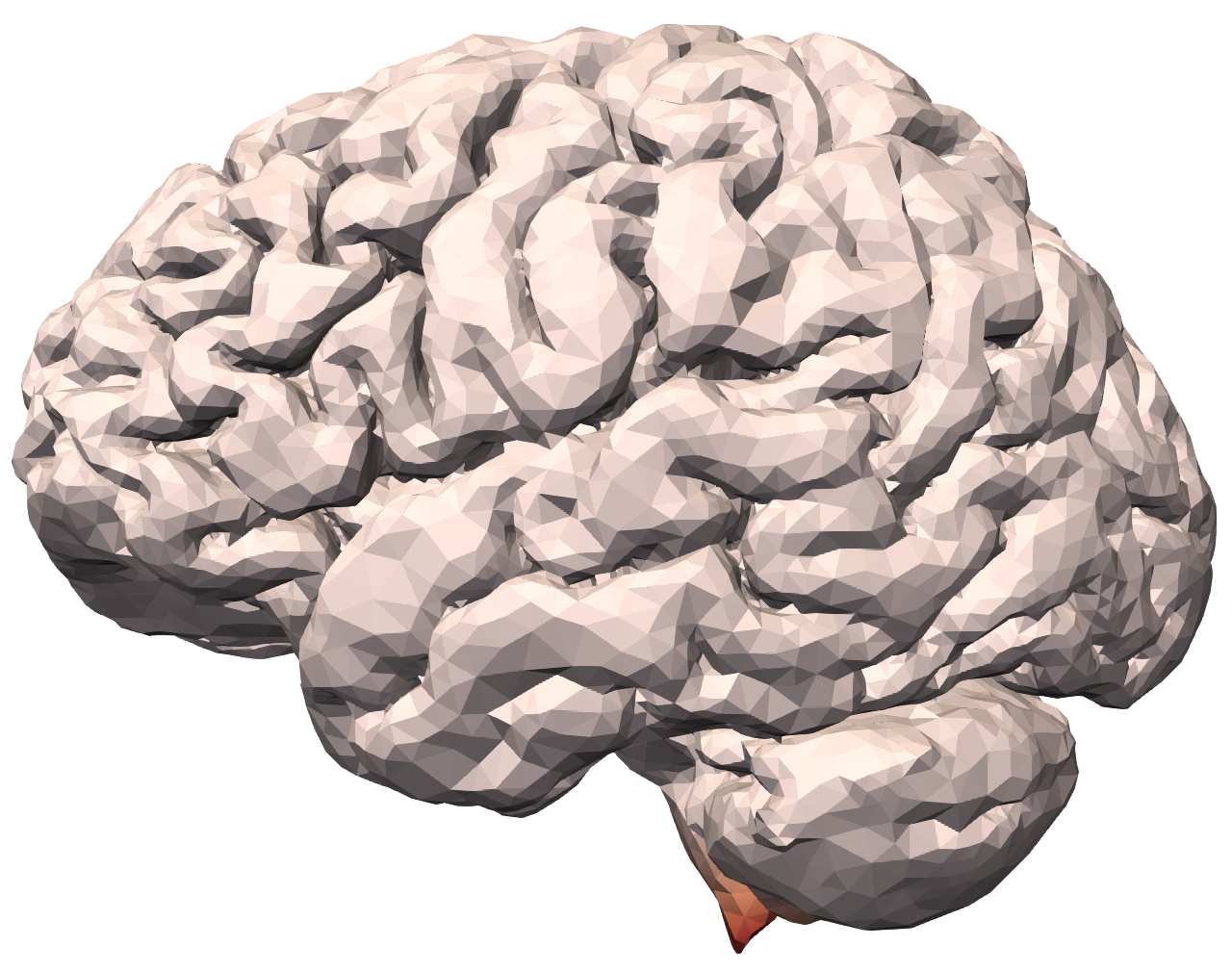}
        \includegraphics[width=0.14\textwidth ]{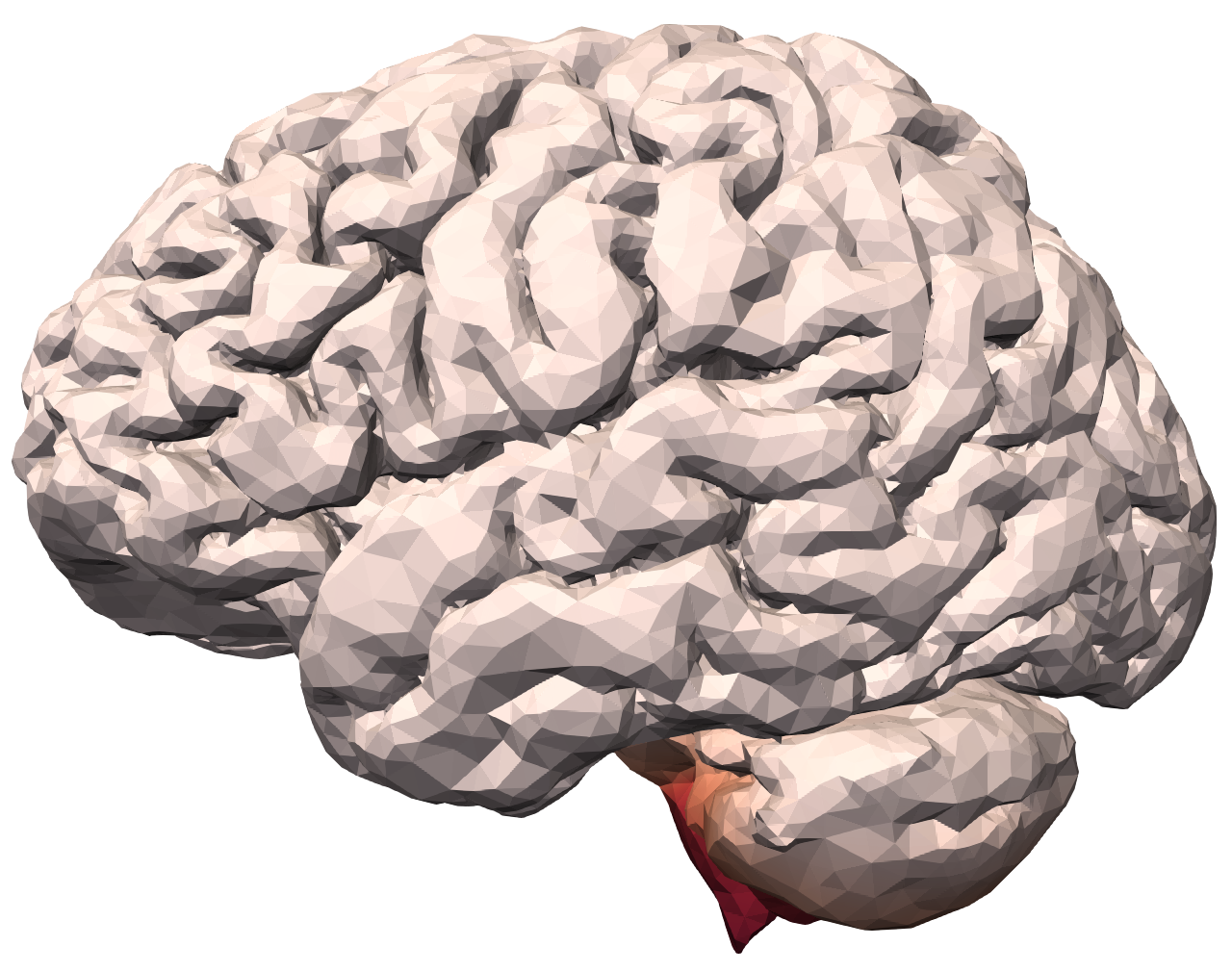}
        \includegraphics[width=0.14\textwidth ]{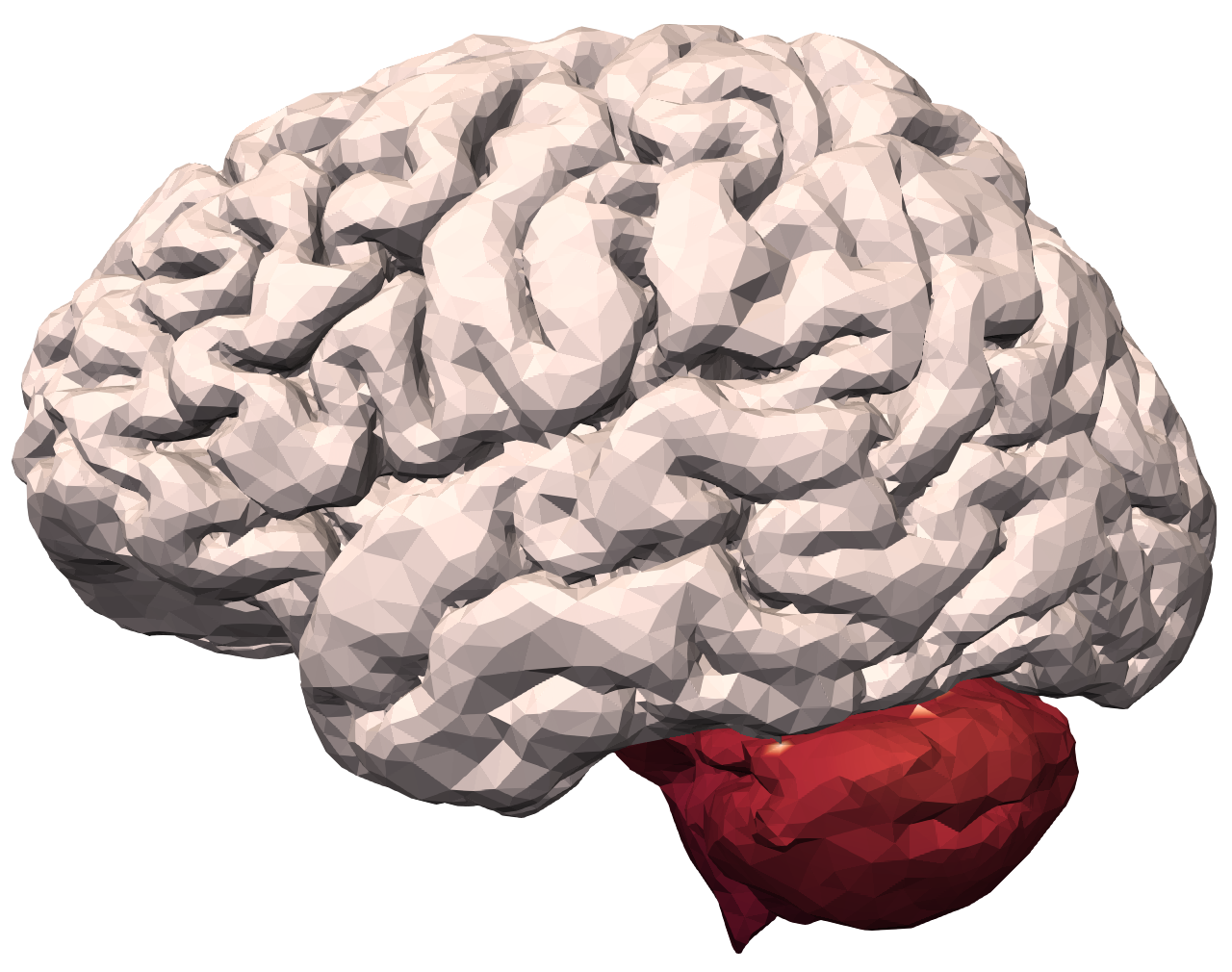}
        \includegraphics[width=0.14\textwidth ]{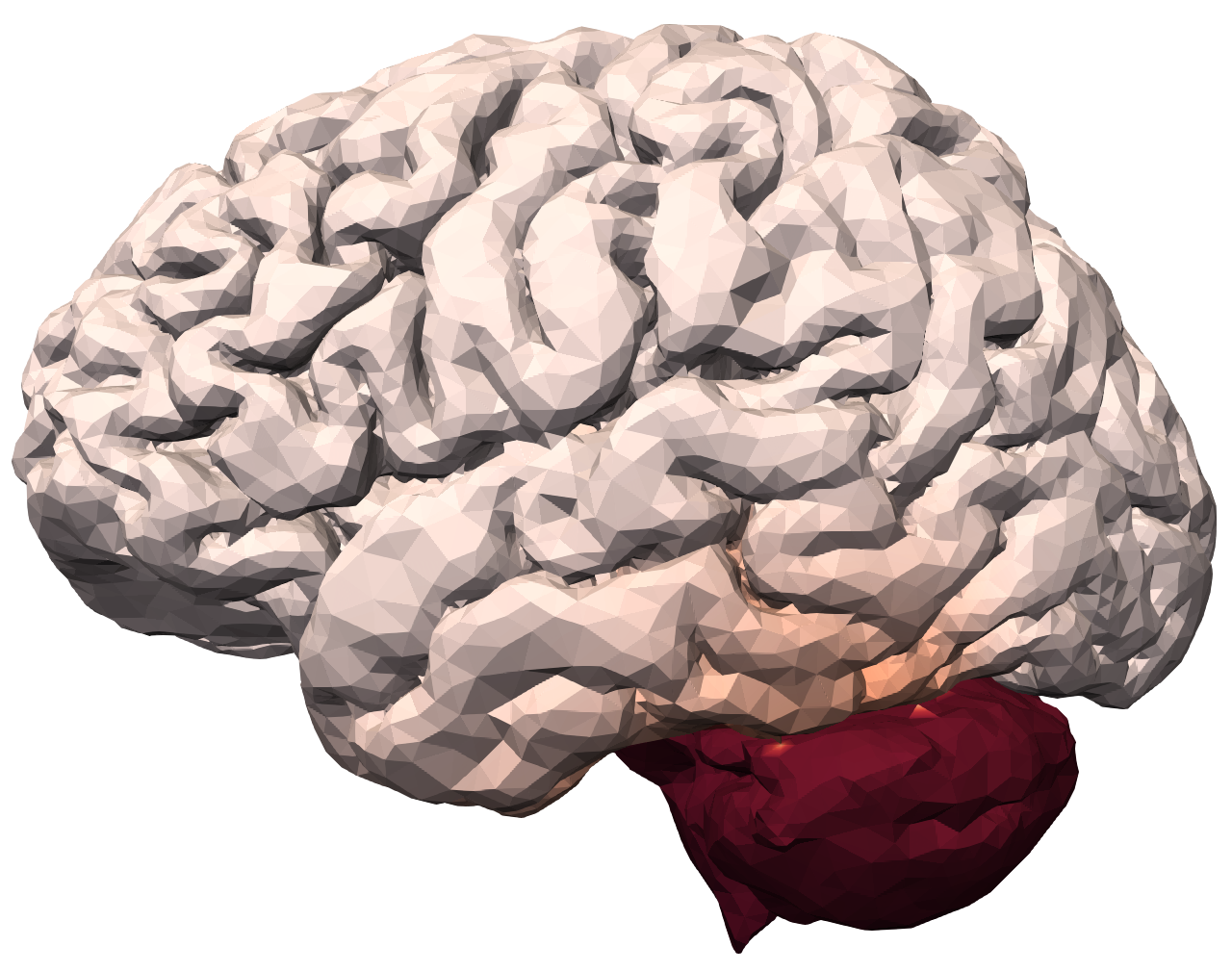}
        \includegraphics[width=0.14\textwidth ]{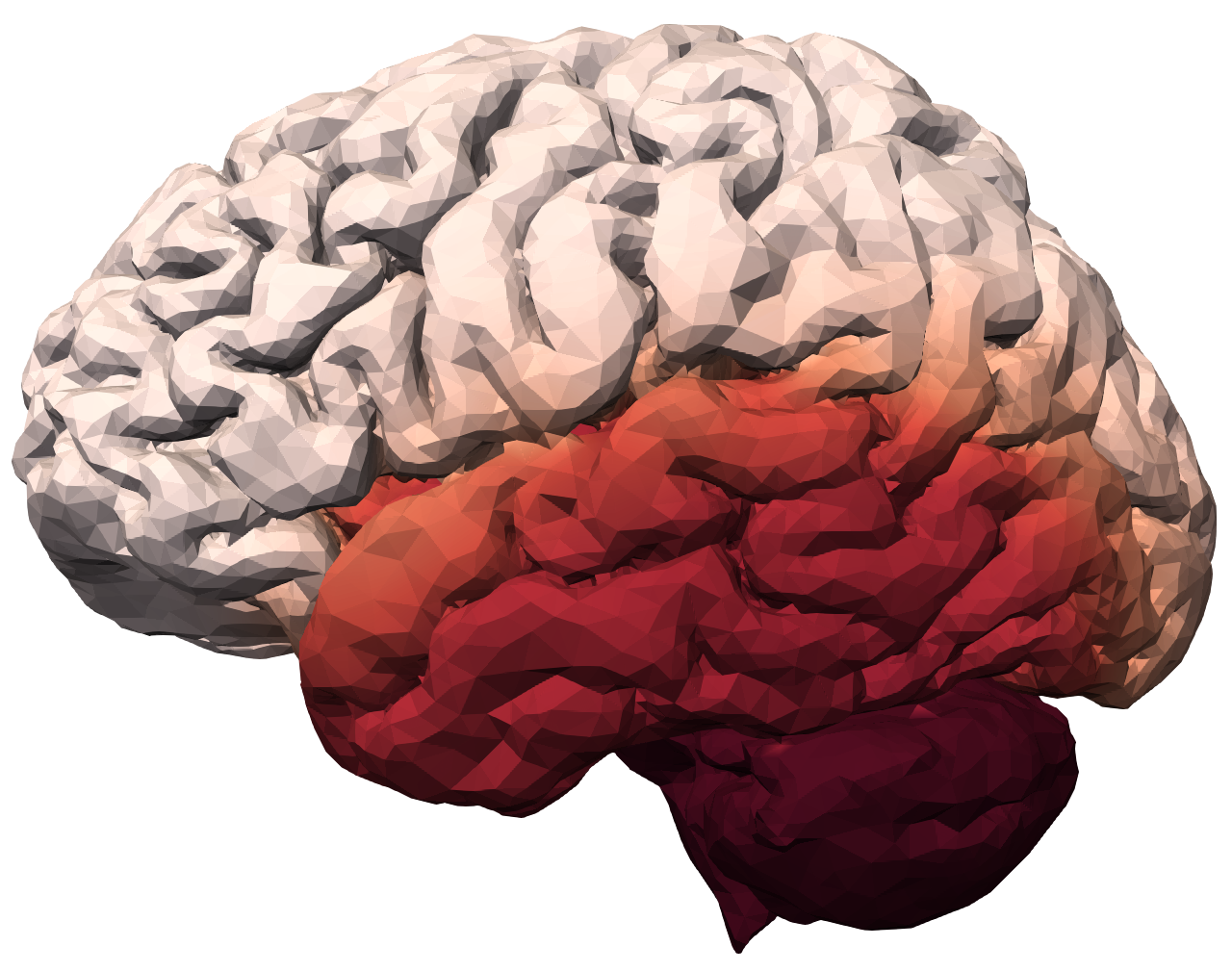}
        \includegraphics[width=0.14\textwidth ]{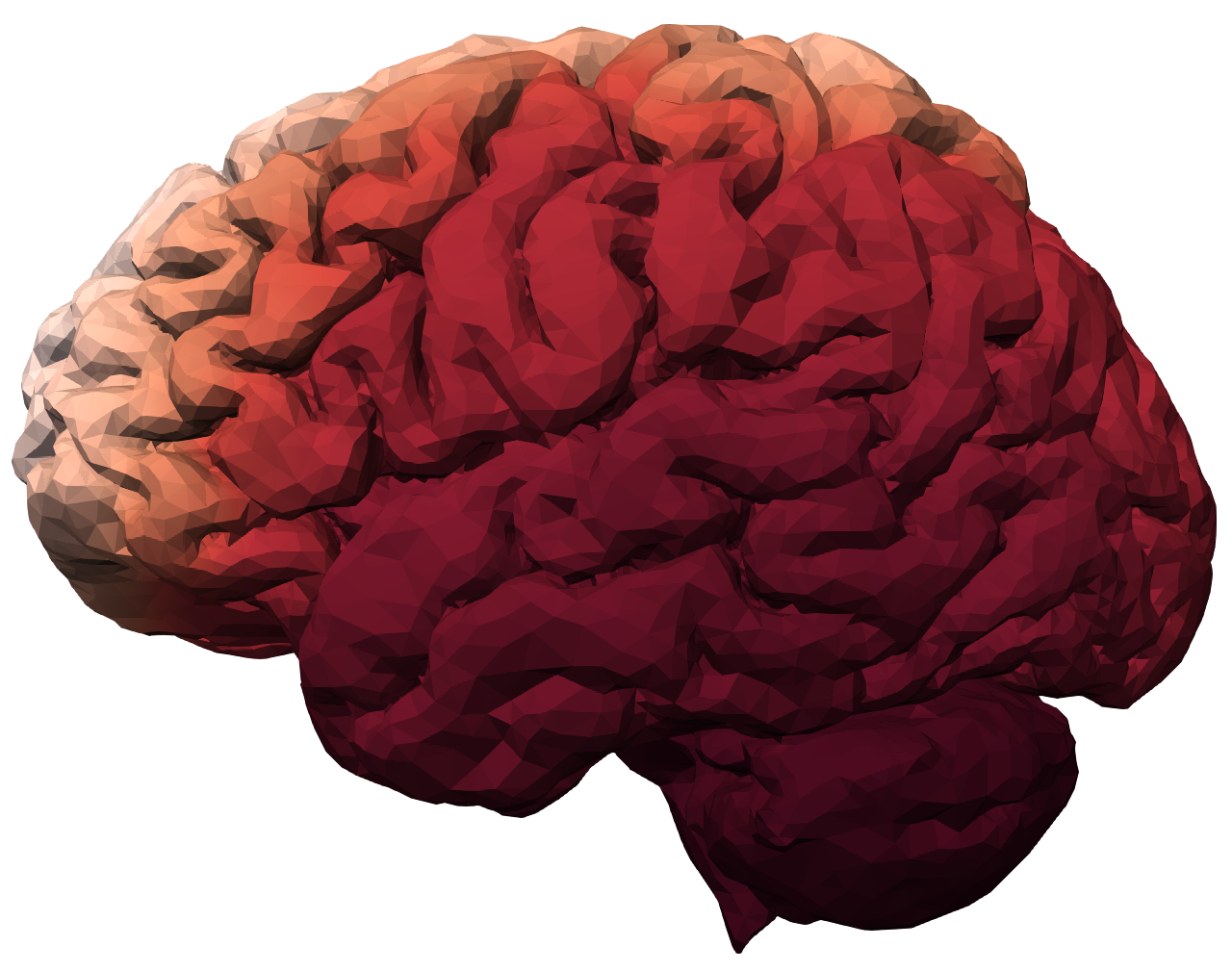}\\
        
        \includegraphics[width=0.3\textwidth ]{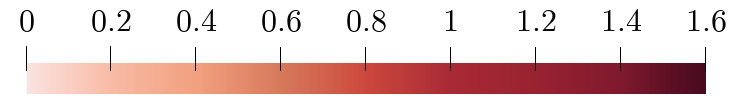}
        \caption{Test case of Section \ref{sec:results}. Snapshot of the computed solution $c_h(\boldsymbol{x},t)$ (blue) and $q_h(\boldsymbol{x},t)$ (red), for $t \in \{ 0, 5, 10, 15, 20, 25 \}$ years.}
        \label{fig:res5}
    \end{figure}

\begin{figure}[ht]
        \centering
        \includegraphics[width=0.3\textwidth ]{final_simulation/scale_c_h.png}\\
        
        \includegraphics[width=0.14\textwidth ]{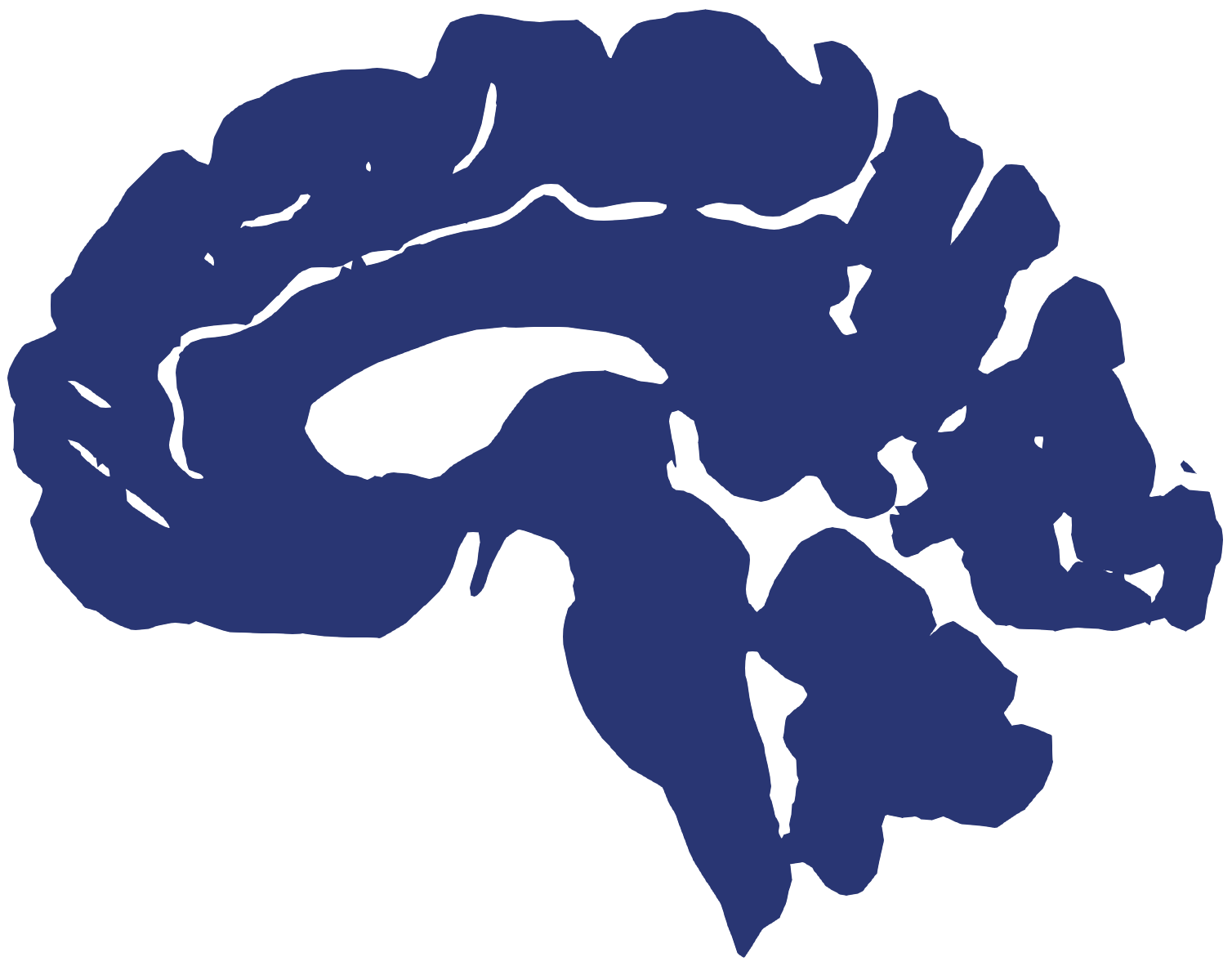}
        \includegraphics[width=0.14\textwidth ]{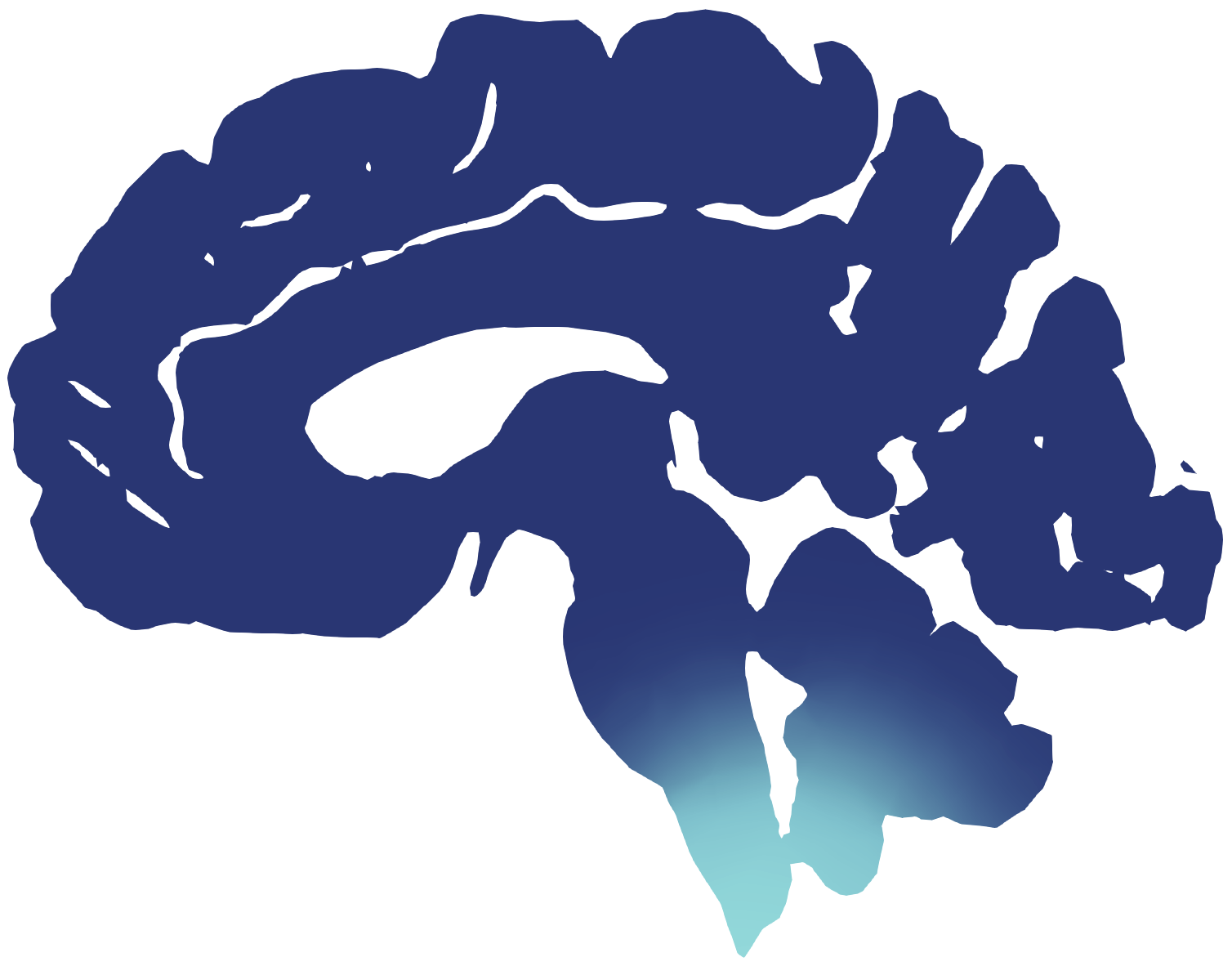}
        \includegraphics[width=0.14\textwidth ]{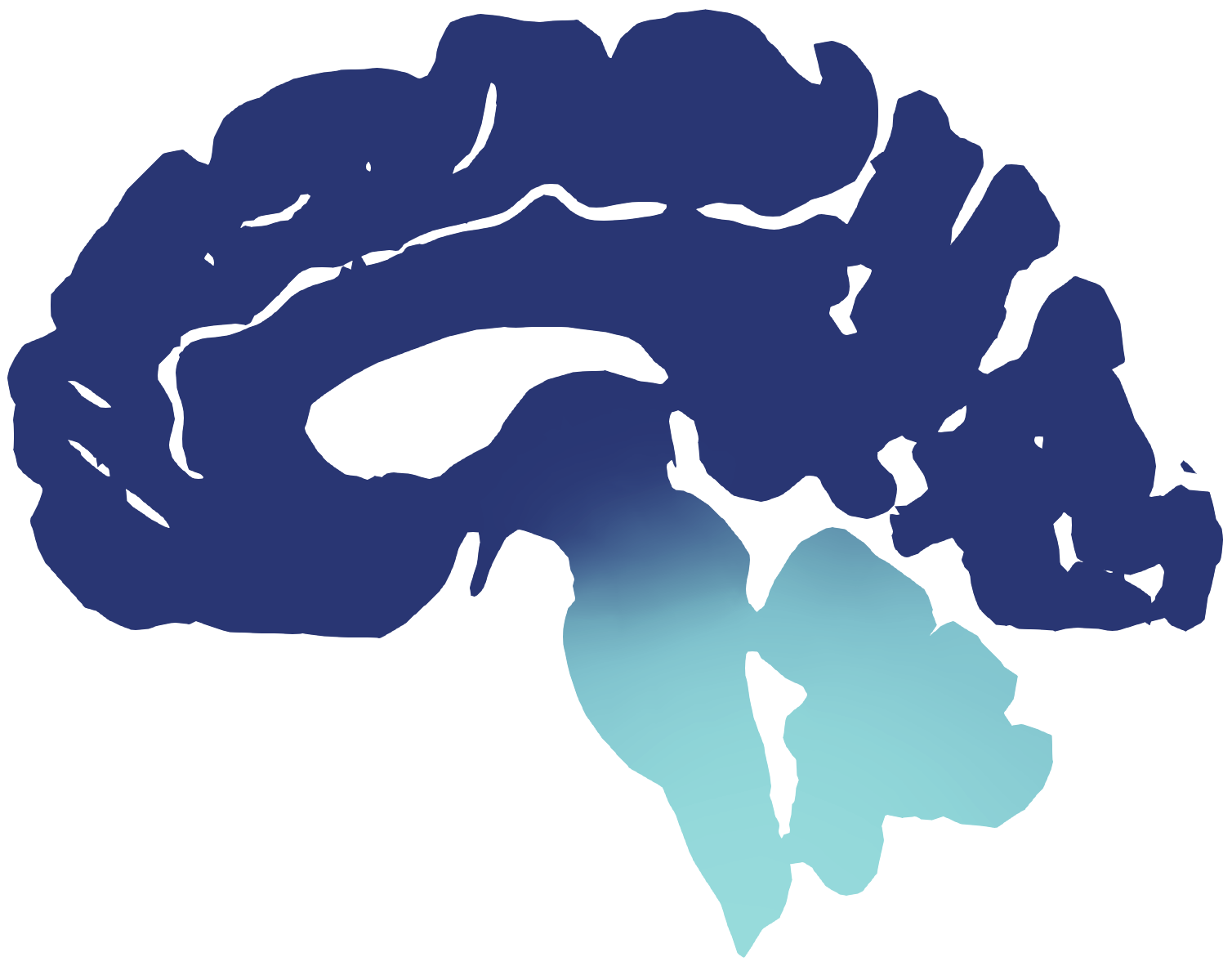}
        \includegraphics[width=0.14\textwidth ]{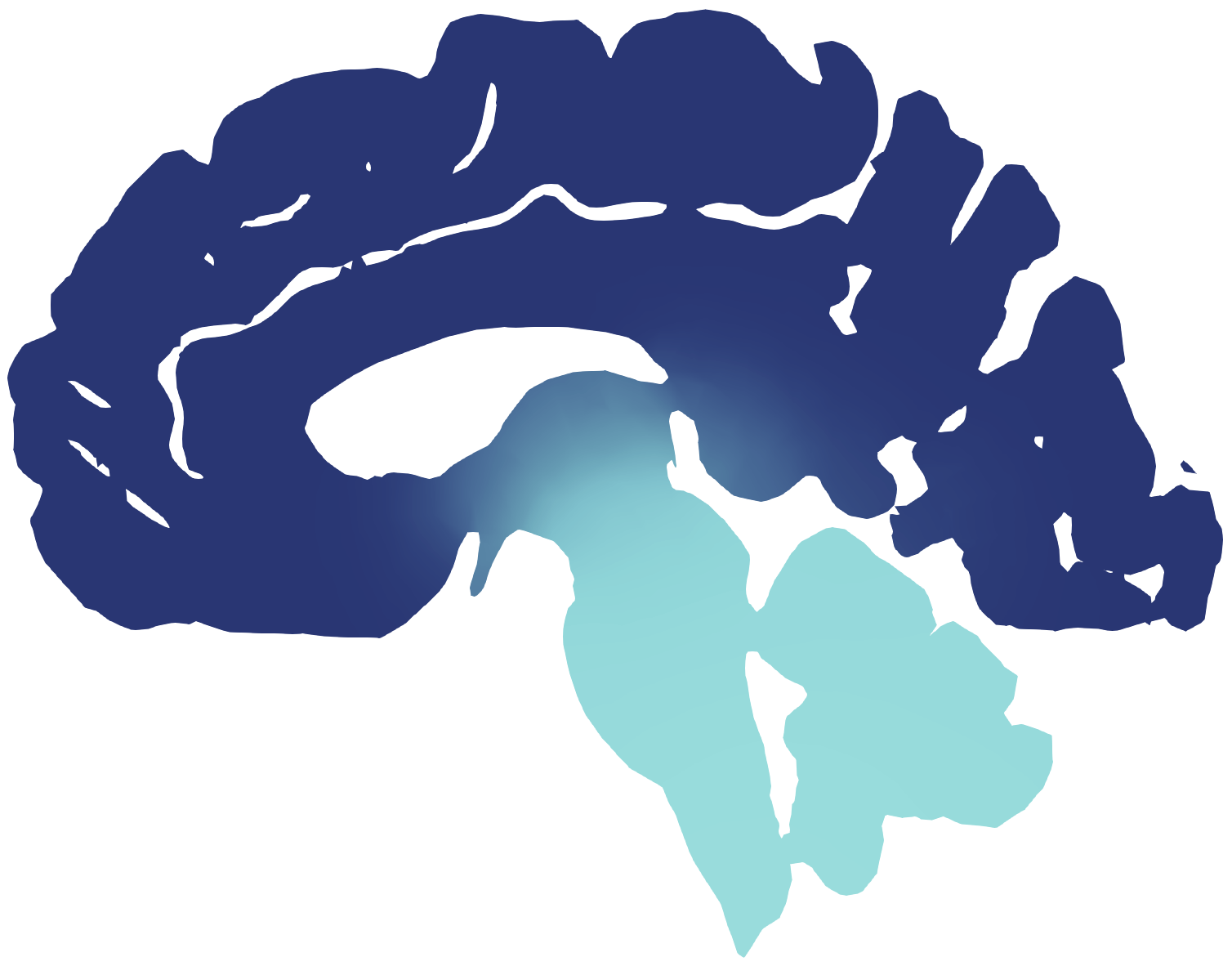}
        \includegraphics[width=0.14\textwidth ]{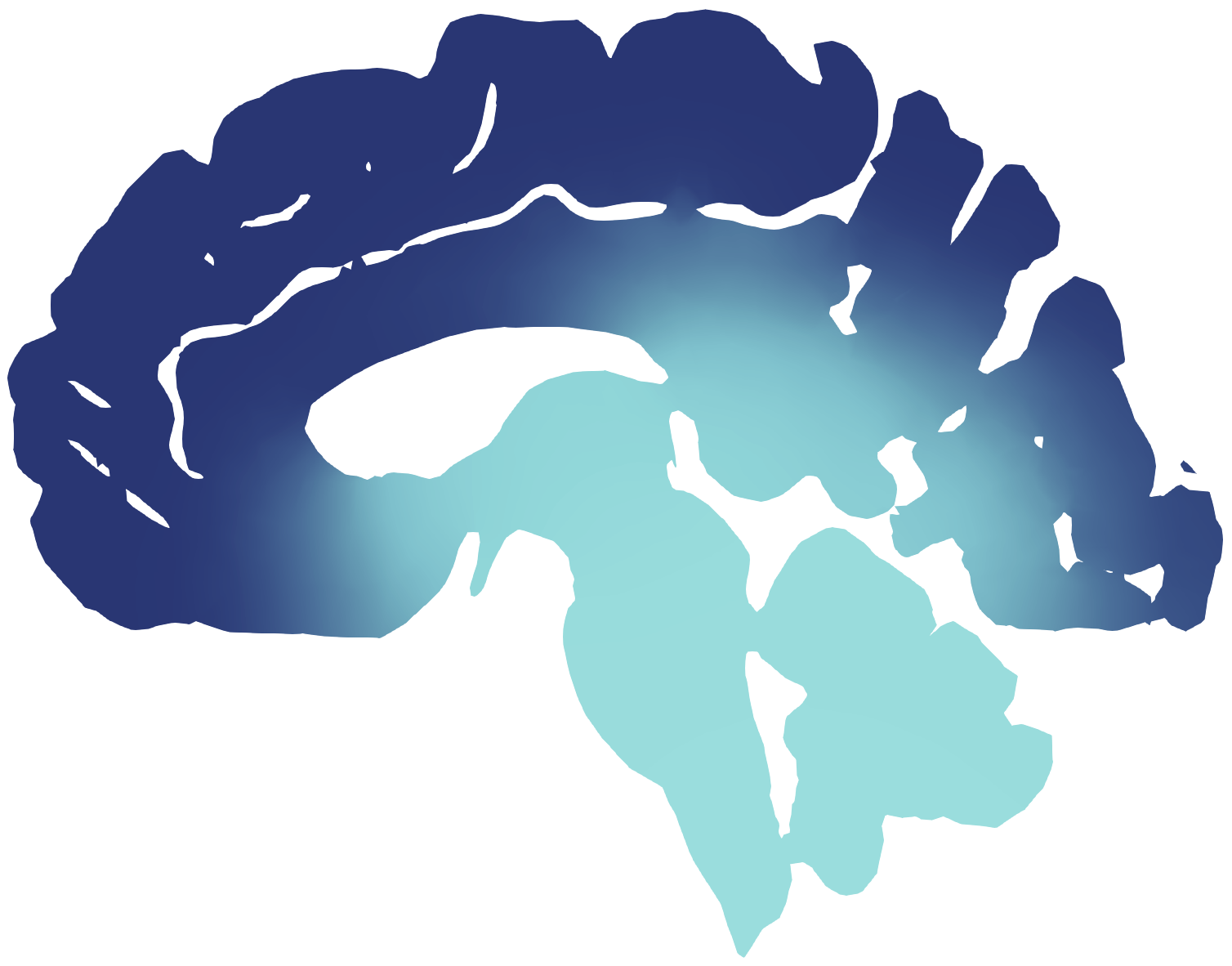}
        \includegraphics[width=0.14\textwidth ]{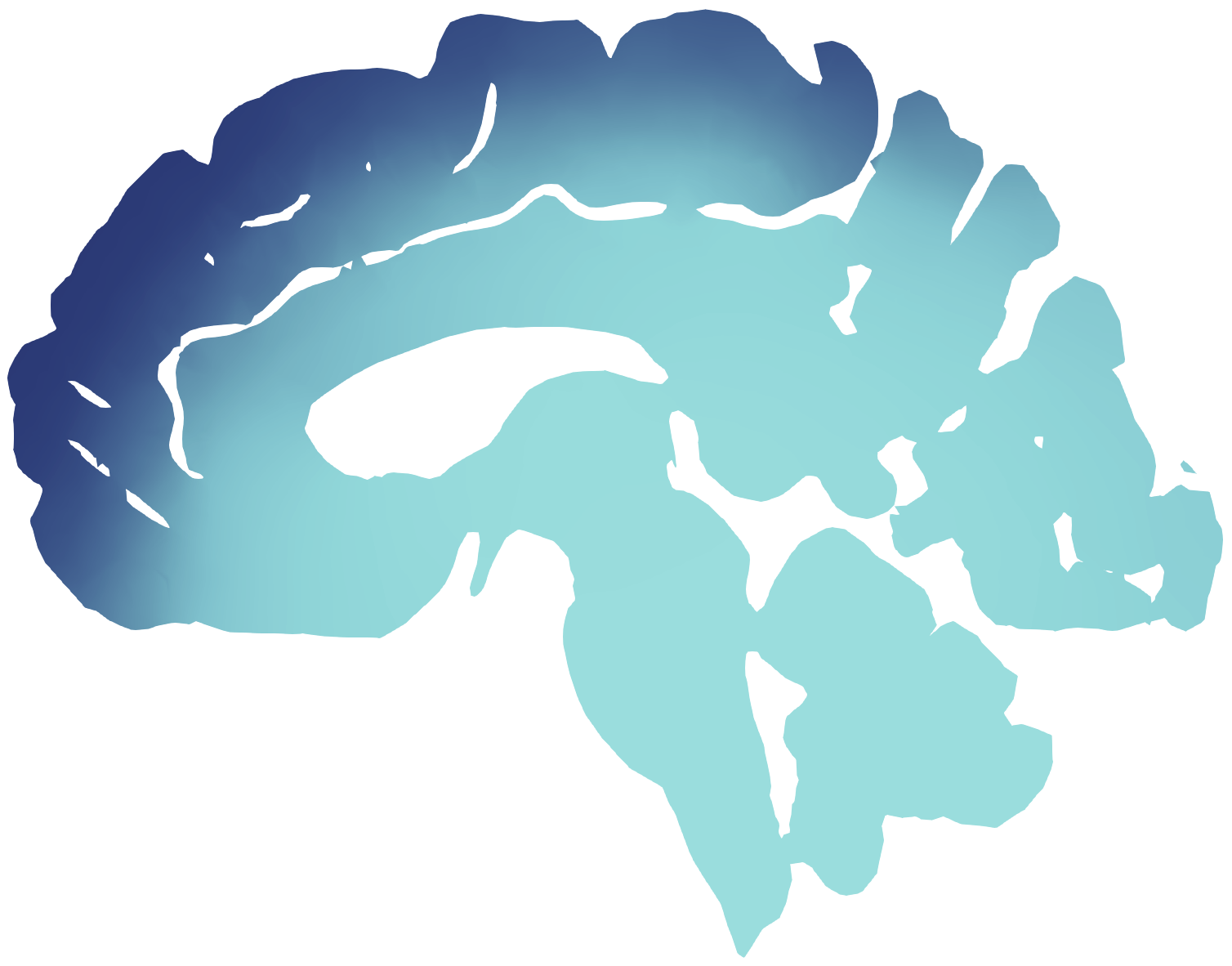}\\
        
        \includegraphics[width=0.14\textwidth ]{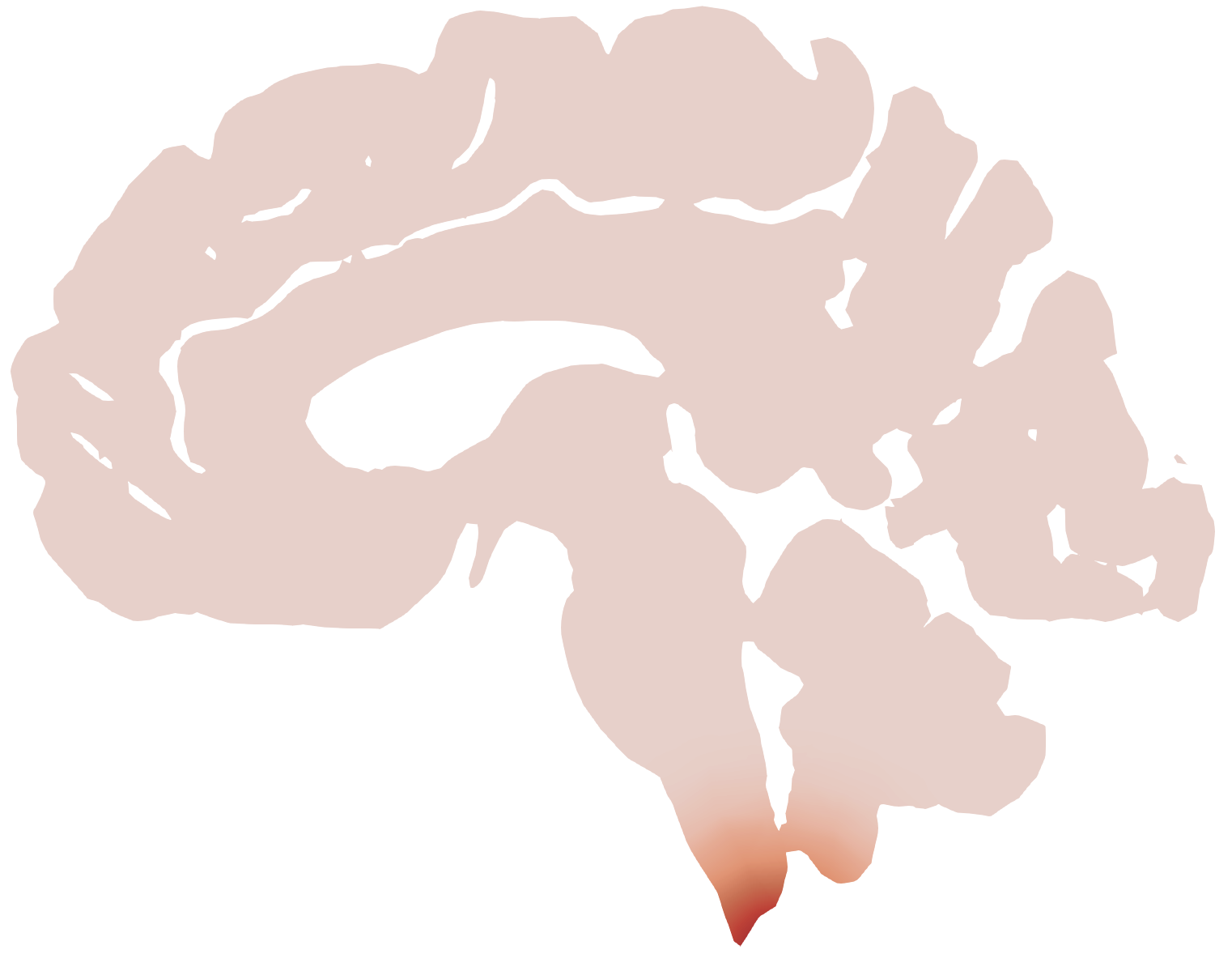}
        \includegraphics[width=0.14\textwidth ]{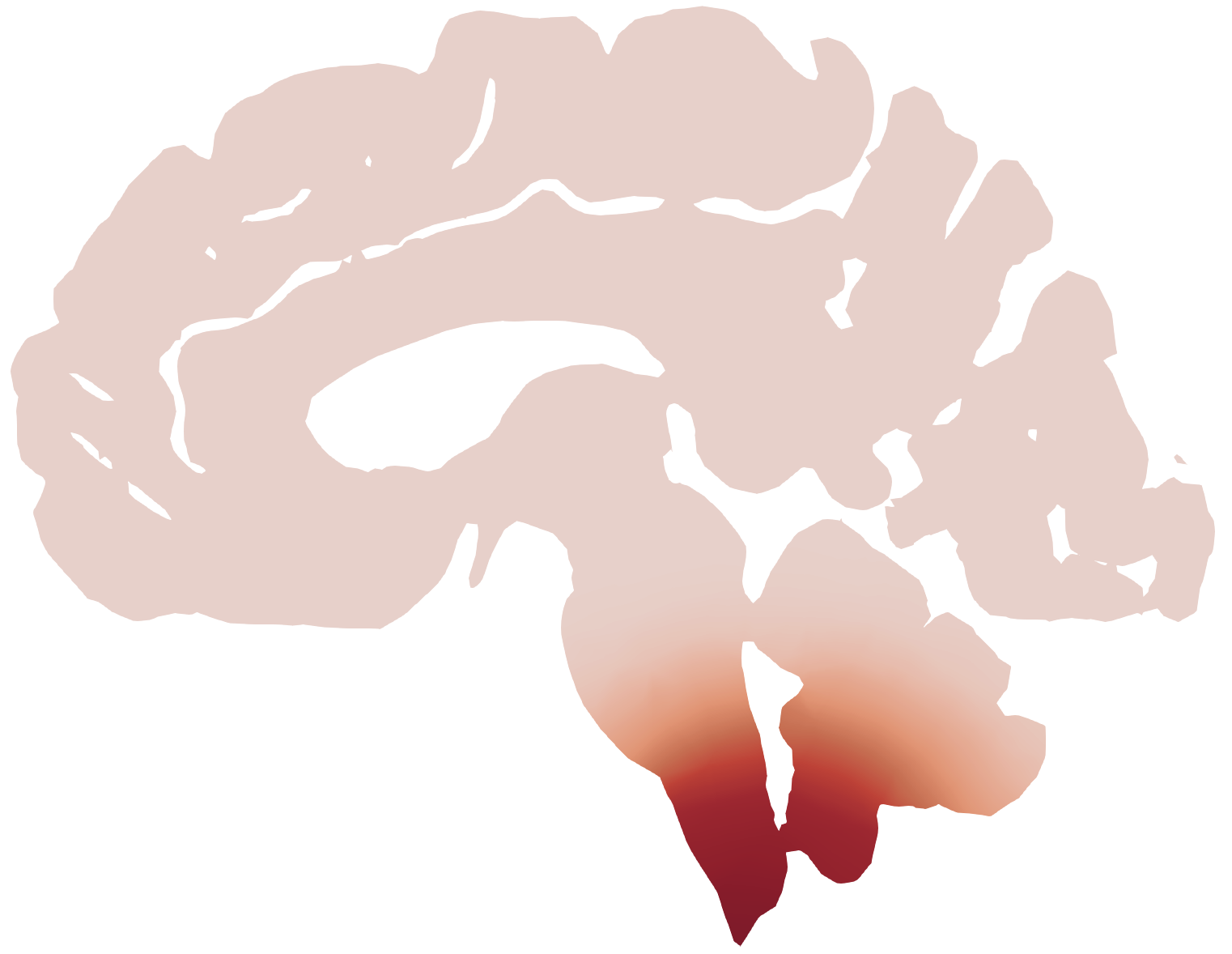}
        \includegraphics[width=0.14\textwidth ]{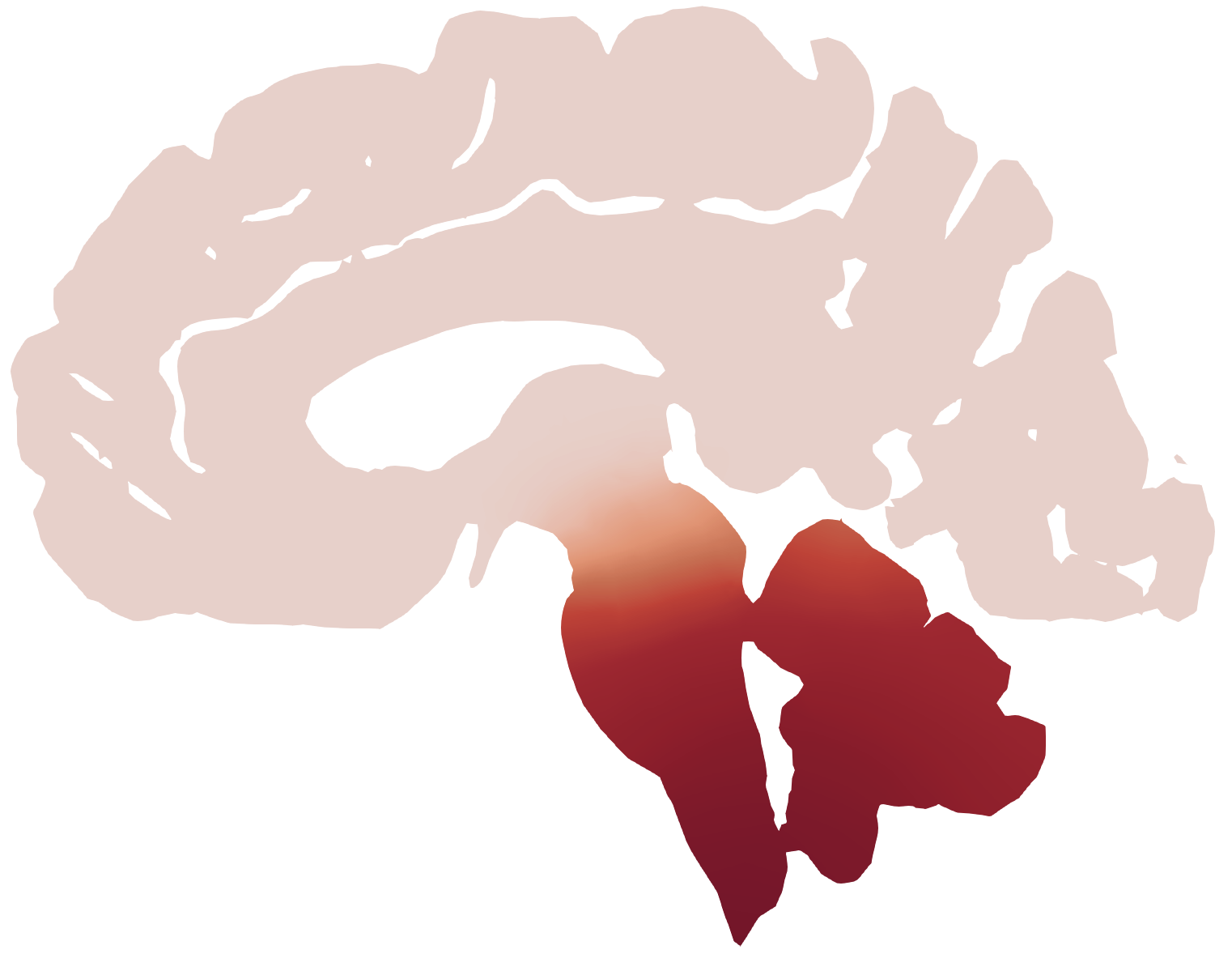}
        \includegraphics[width=0.14\textwidth ]{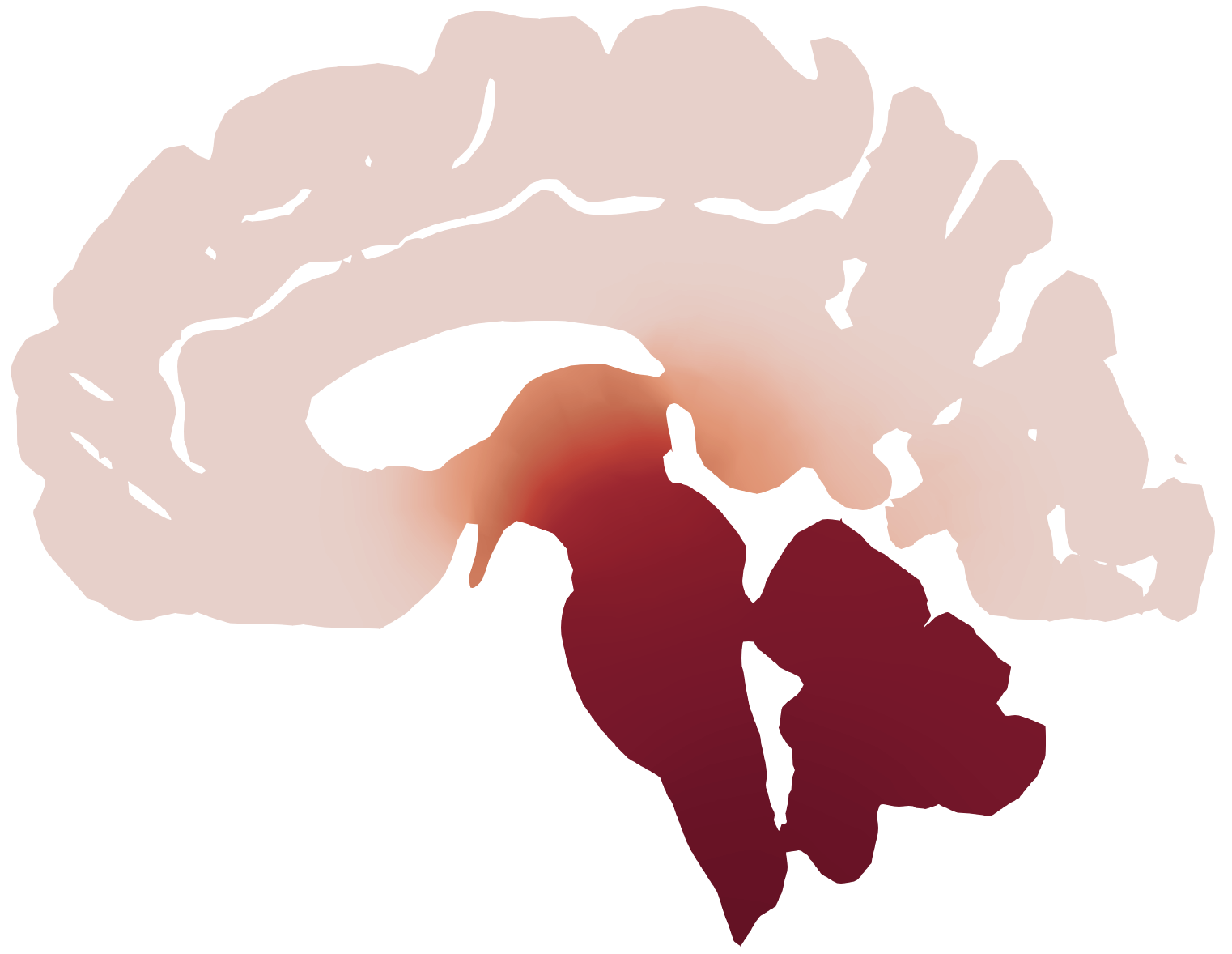}
        \includegraphics[width=0.14\textwidth ]{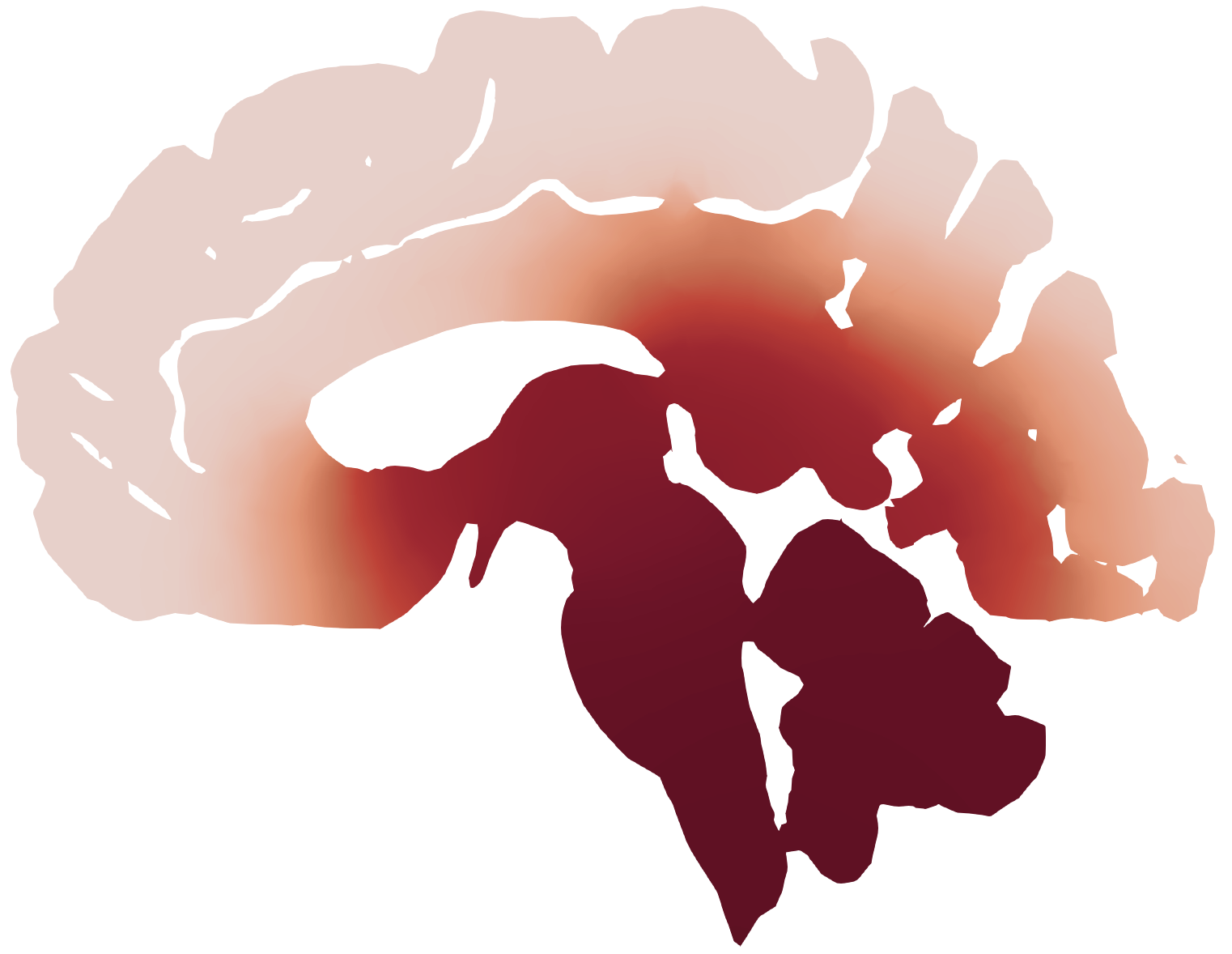}
        \includegraphics[width=0.14\textwidth ]{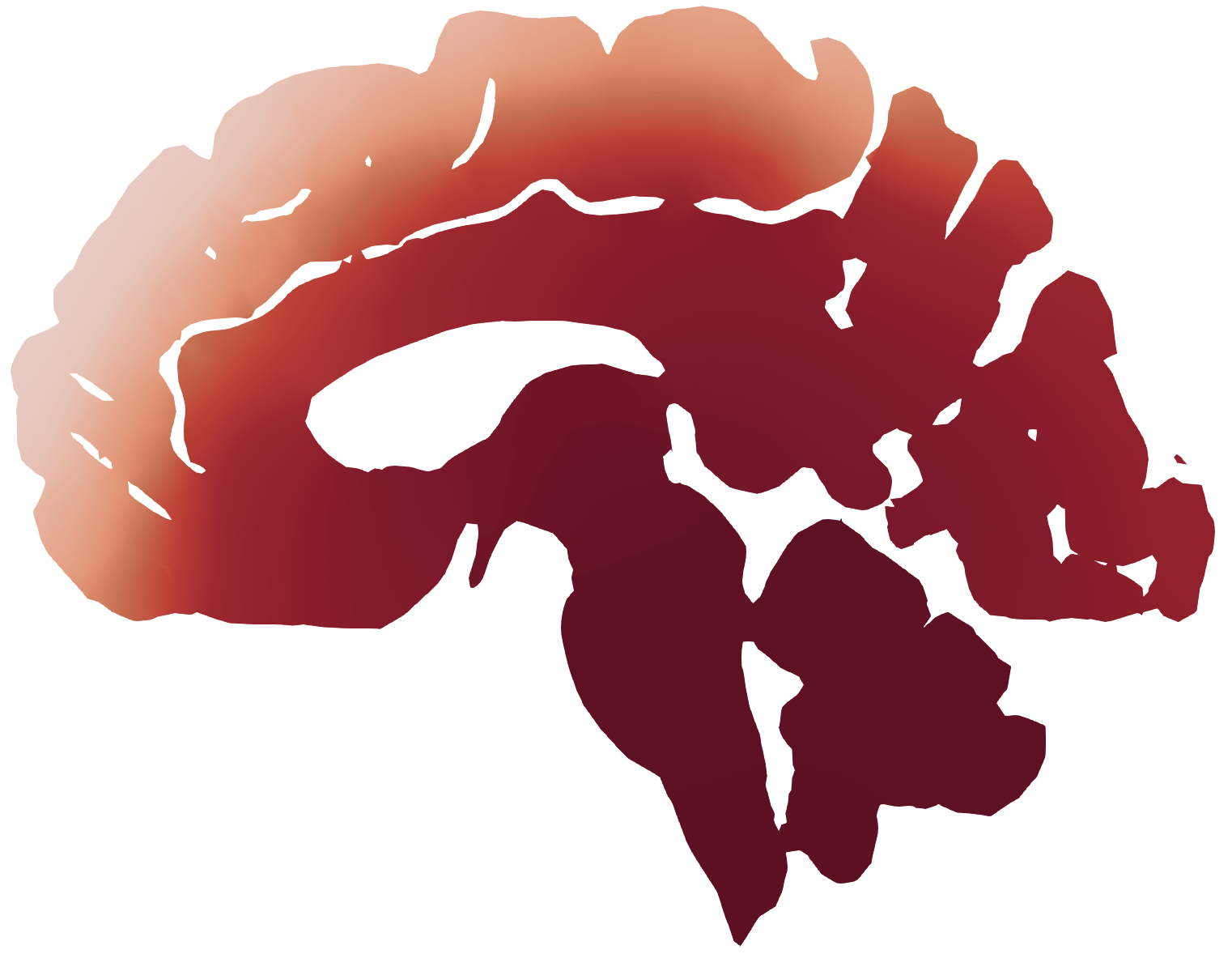}\\
        
        \includegraphics[width=0.3\textwidth ]{final_simulation/scale_q_h.png}
        
        \caption{Test case of Section \ref{sec:results}. Snapshot of the computed solution $c_h(\boldsymbol{x},t)$ (blue) and $q_h(\boldsymbol{x},t)$ (red), for $t \in \{ 0, 5, 10, 15, 20, 25 \}$ years in a sagittal medial section.}
        \label{fig:res6}
    \end{figure}

    As discussed in Section \ref{sec:model}, the system is expected to evolve from the unstable initial equilibrium point to the stable pathological equilibrium $\boldsymbol{e}_1=(0.3,1.5)$. Figure \ref{fig:res5} shows the temporal evolution of healthy and misfolded proteins across the domain, with a temporal resolution of $5$ years, accurately reaching the expected equilibrium. Figure \ref{fig:res6} depicts instead the numerical solution obtained over the medial sagittal section of the brain mesh.
    The results are in alignment with the ones presented in \cite{cortiAntonietti2023}, which were obtained starting from the same data but exploiting the Fisher-Kolmogorov model. Instead, a major discrepancy can be noticed with the spreading pattern of two-dimensional models. For example, consider \cite{corti}: not simply connected brain sections strongly affect protein diffusion (see Figure \ref{fig:res6}). Lateral ventricles would represent an obstacle to protein propagation in two dimensions. However, not simply connected points in 2D can be reached with continuity in the 3D brain mesh. Obviously, this affects the overall spreading of the wave solution.
    Also, the evolution time is of the order of 25--30 years, which is comparable to the medical references \cite{braak2003}, and similar to other simulation results \cite{cortiAntonietti2023, corti, fornari2019prion}. Furthermore, diffusion directions are globally coherent with the medical literature since they reproduce the stereotypical pattern described in \cite{braak2003}. For instance, there is no cerebral cortex involvement in the absence of lesions in the brain stem, and early disease stages affect only the \textit{medulla oblongata} and the \textit{pontine tegmentum}.\vspace{0.5em}
    
    From the phenomenological point of view, \textit{nigral} and {extranigral} induction sites vary in their susceptibility to disease-related changes and follow a well-described temporal sequence of degeneration \cite{braak2003}. Starting from the motor nucleus, the misfolded protein concentration locally increases in the seeding area. Subsequently, it invades the lower brain stem along the axonal network, then spreads firstly into the mesocortex and eventually into the neocortex, coherently with Braak's staging theory \cite{braak2003}. The concentration of healthy proteins evolves almost symmetrically, giving insight into cell loss.
    
    \subsection{Validation versus Braak's stages}
    \label{subsection:braak_stages}
\begin{figure}[ht]
        \centering
        \includegraphics[width=0.19\textwidth]{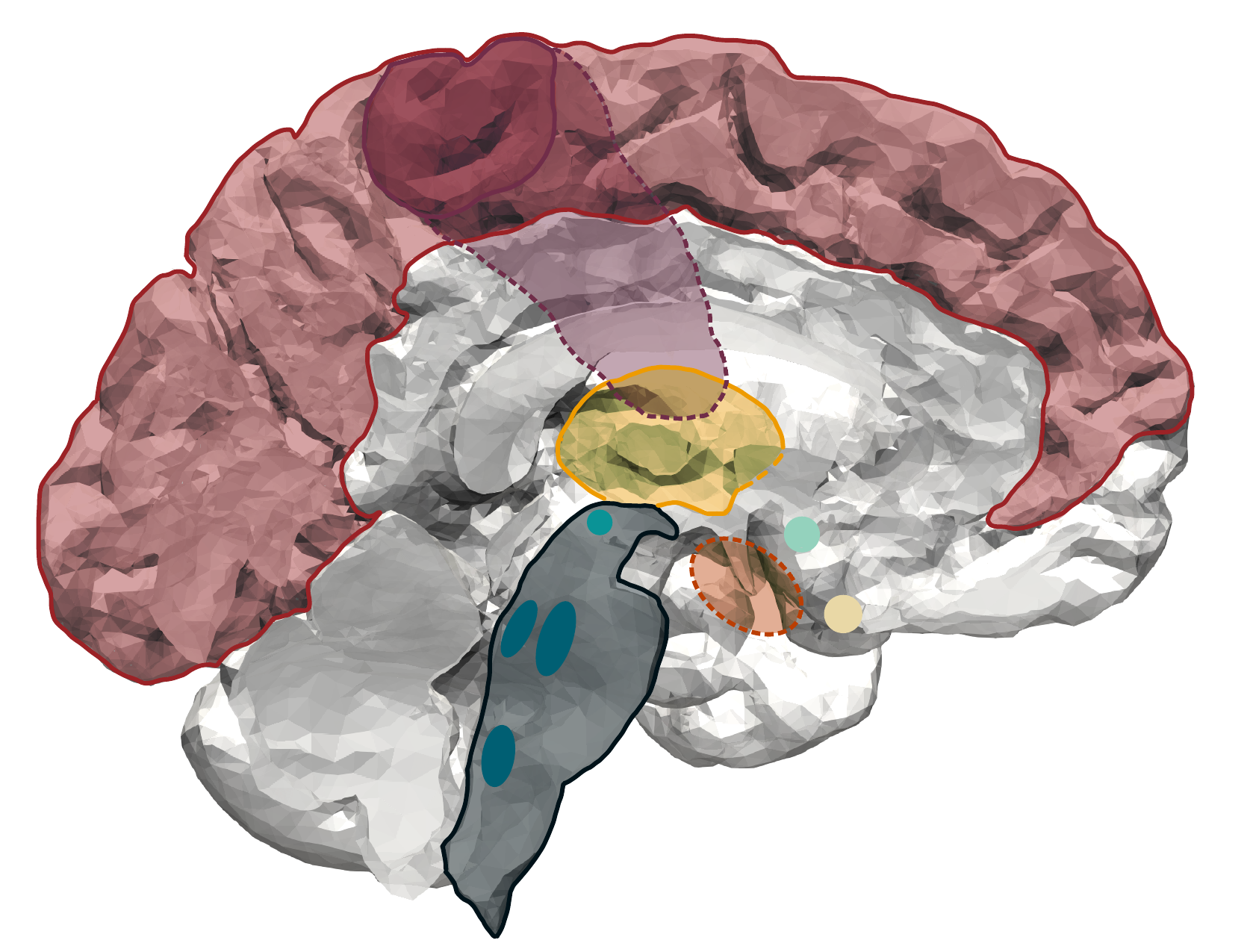}
        \includegraphics[width=0.8\textwidth]{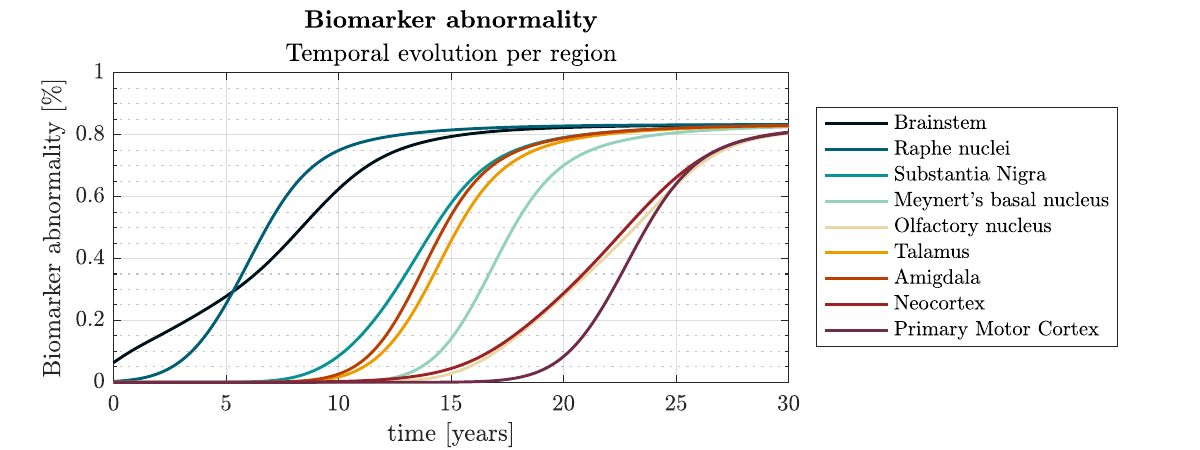}
        \caption{Test case of Section \ref{subsection:braak_stages}. Temporal evolution of biomarker curves for different brain regions of clinical interest \cite{braak2003}.}
        \label{fig:res8}
    \end{figure}

    Refer to Figures \ref{fig:res8}, \ref{fig:res2}, \ref{fig:res1}, \ref{fig:res3}, and to Figure \ref{fig:res4}: the former depicts biomarker abnormality curves $B^\alpha_i (t)$ from the prion-like spreading of $\alpha$-synuclein across different brain regions $R_i \subset \Omega$ of clinical interest \cite{braak2003}, where we define $B^\alpha_i (t)$ as:
    \begin{equation*}
        B^\alpha_i (t) = \frac{\int_{R_i}q(\boldsymbol{x},t)\mathrm{d}\boldsymbol{x}}{\int_{R_i}c(\boldsymbol{x},t)+q(\boldsymbol{x},t)\mathrm{d}\boldsymbol{x}}.
    \end{equation*}
    All of them are smooth sigmoidal curves, in agreement with clinical biomarkers of neurodegenerative disease progression. It is also possible to identify the six different Braak's stages \cite{braak2003}. During stage one, the lower brainstem is affected: firstly the dorsal motor nucleus of the vagus nerve in the medulla oblongata, and then the pons reach a critical concentration value. After, during stage 2, the Raphe nuclei are also affected and the misfolded proteins travel up the brainstem into the pontine tegmentum. During stage 3, the \textit{substantia nigra} is reached: the \textit{parse compacta}, together with the basal nucleus of Meynert starts showing misfolding. Stage 4 coincides with the mesocortex invasion, while the neocortex is still healthy. However, the amygdala is now compromised, as the anterior olfactory nucleus is. These two stages represent the start of sleep and motor disturbances. Stage 5 also involves the neocortex, thus the temporal, parietal, and frontal lobes. The substantia nigra is almost completely compromised. In stage 6, misfolded proteins have fully invaded the neocortex, affecting also, for example, the primary motor cortex, and other sensory areas. Cognitive impairment is now severe. The olfactory cortex appears to be years late if compared to clinical observations. This discrepancy arises because, in our simulation, we consider only the medulla oblongata as the seeding region. Consequently, the infection of the olfactory cortex occurs significantly later than the average observed in patients, where misfolded proteins are detected early in the olfactory cortex due to its proximity to the olfactory lobe.\vspace{0.5em}
    
\begin{figure}[ht]
        \centering
        \includegraphics[width=0.85\textwidth]{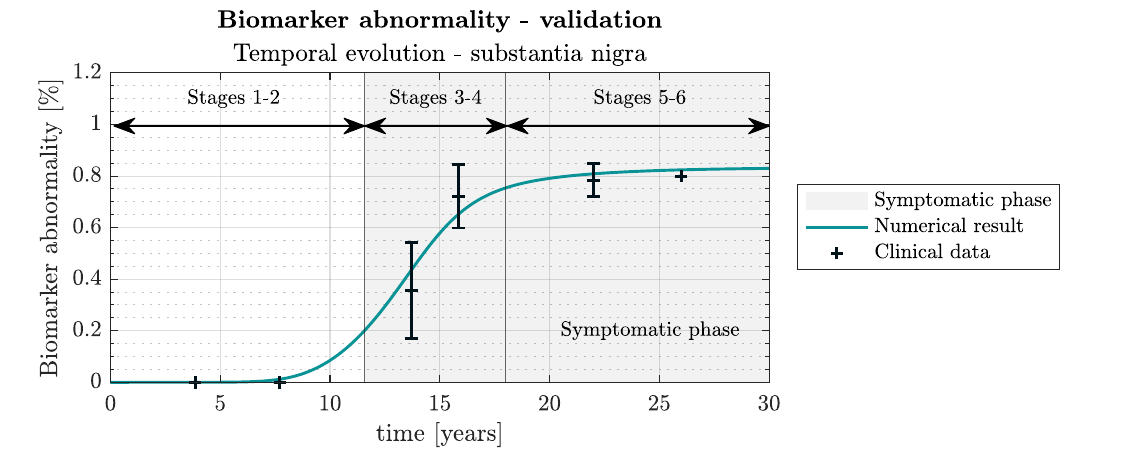}
        \caption{Test case of Section \ref{subsection:braak_stages}. Validation for the biomarker abnormality for the substantia nigra brain region. The green line represents the biomarker abnormality curve computed from the numerical results obtained from the simulation of Section \ref{sec:results}. Clinical data are taken from \cite{braak2003}, Table 2. The severity of the disease is mapped on the range $[0, 0.8]$, and for each stage, we take the average value and show the symmetric interval with amplitude twice the standard deviation.}
        \label{fig:SN}
    \end{figure}

    The substantia nigra impairment is a crucial point in the disease progression. Indeed, at this stage,  symptoms become evident \cite{braak2003}. Consequently, this region is particularly significant in the study of Parkinson's disease evolution. By establishing a critical concentration threshold of $q_\text{crit} = 0.2$, at which point a brain region is deemed compromised, we can determine the temporal duration of each Braak stage in our simulation. Consistently with \cite{braak2003}, we align the onset of stage 3 and stage 5, with substantia nigra and neocortex compromise, respectively. Additionally, by considering the data in Table 2 of  \cite{braak2003}, we can validate the temporal evolution of the biomarker in the substantia nigra, $B^\alpha_\text{SN}(t)$. Figure \ref{fig:SN} illustrates the curve $B^\alpha_\text{SN}(t)$ obtained from the numerical simulation of this Section alongside the averaged values of 110 Parkinson's disease patients at each Braak stage. The result is promising, yet it is heavily constrained by the semi-quantitative nature of the literature-acquired data.
    
    \subsection{Sensitivity analysis}
    \label{subsection:sensitivity_analysis}
\begin{figure}[ht]
        \centering
        \includegraphics[width=0.45\textwidth]{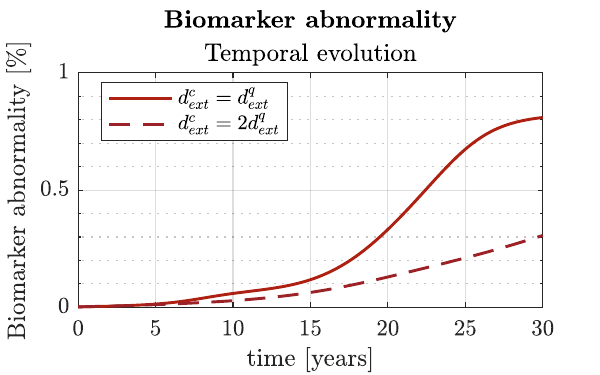}
        \caption{Test case of Section \ref{subsection:sensitivity_analysis}. Sensitivity analysis for two different mechanisms of extracellular diffusion.}
        \label{fig:res7}
    \end{figure}

    This section briefly discusses the results of a sensitivity analysis study that explores two parameter variations not covered in the literature.
    \paragraph{Time step length} First, we address the time step length. Even though a longer $\Delta t$ does not cause spurious oscillations in the numerical solution representation, it induces a much lower spreading over the domain. A proper choice of $\Delta t$ is in the order of $10^{-2} \ \mathrm{years}$.
    \paragraph{Diffusion coefficients}  Second, we study the relative contribution of extracellular diffusion of healthy and misfolded proteins: indeed, the different steric hindrances might induce different gradient-driven velocities. Increasing the extracellular diffusion of the healthy component induces a faster spreading of the healthy proteins among the domain. Thus, it takes a longer time to locally accumulate misfolded proteins and cause the disease. As an instance, the choice $d_\text{ext}^c = 2 d_\text{ext}^q = 16 \mathrm{mm}^2\mathrm{years}^{-1}$ induces a temporal delay for the biomarker curve of over 10 years. Figure \ref{fig:res7} depicts the temporal evolution of the global biomarker abnormality curve versus the period of $30$ years for the two choices $d_\text{ext}^c = d_\text{ext}^q$ and $d_\text{ext}^c = 2d_\text{ext}^q$. In the latter case, the saturation time is reached much later than expected by medical literature. Therefore, deviating from clinical observations over several years, it does not appear to be a good modeling choice.
   {
   \newgeometry{top=0cm,bottom=1.5cm,left=2.5cm,right=2.5cm,marginparwidth=0cm}
    \thispagestyle{empty}
    \begin{figure}[H]
        \centering
        \begin{subfigure}{0.12\textwidth}
            \includegraphics[width=0.7\textwidth]{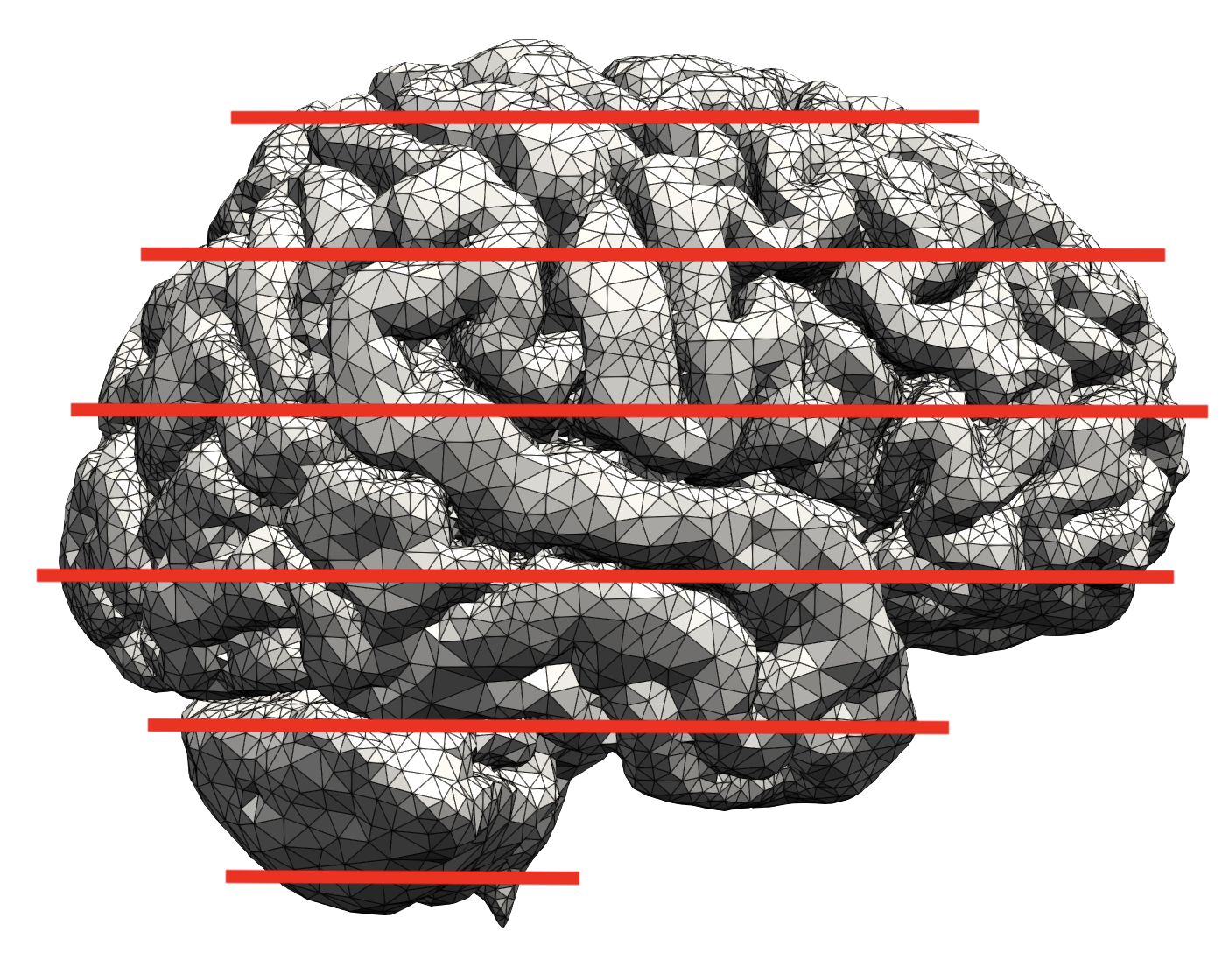}
        \end{subfigure}
        \begin{subfigure}{0.33\textwidth}
            \includegraphics[width=0.32\textwidth]{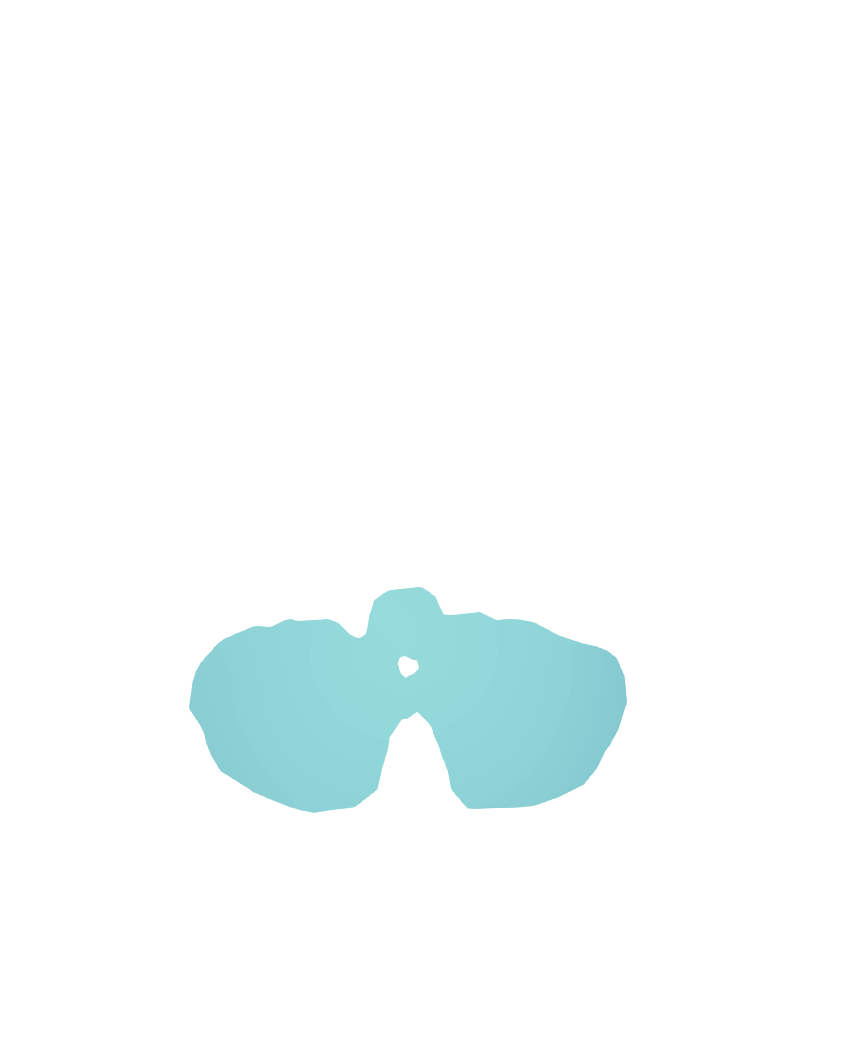}
            \includegraphics[width=0.32\textwidth]{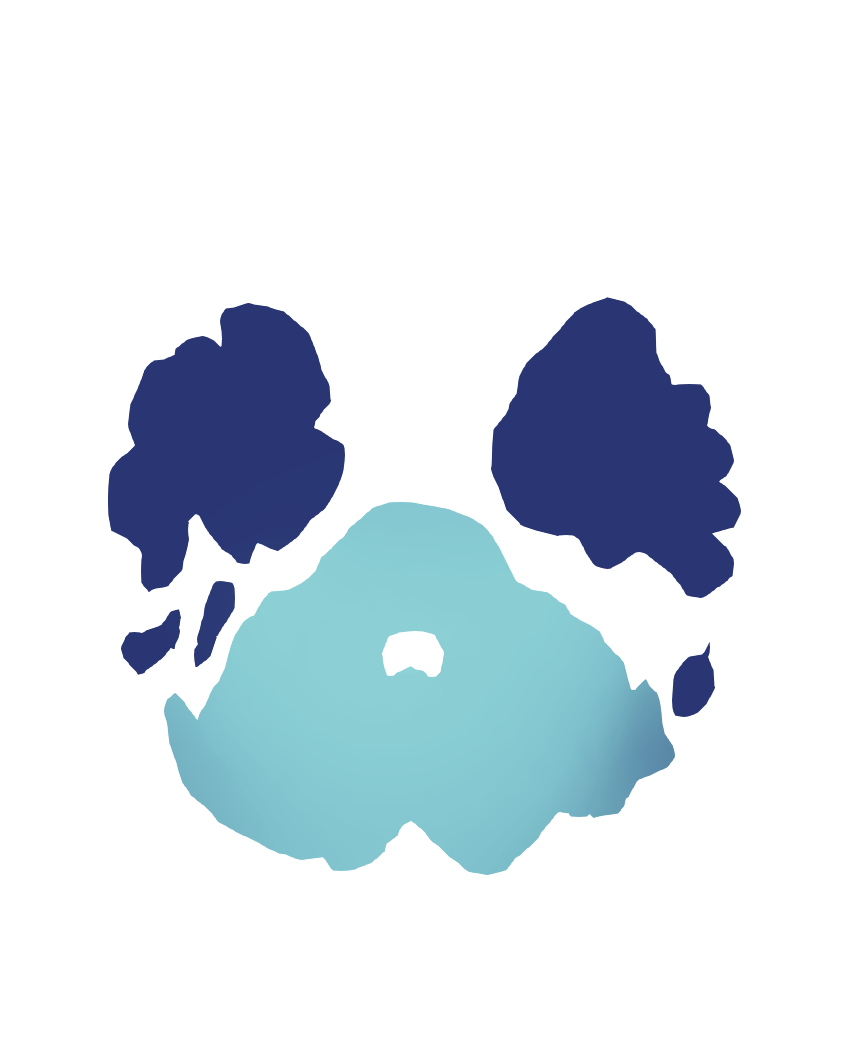}
            \includegraphics[width=0.32\textwidth]{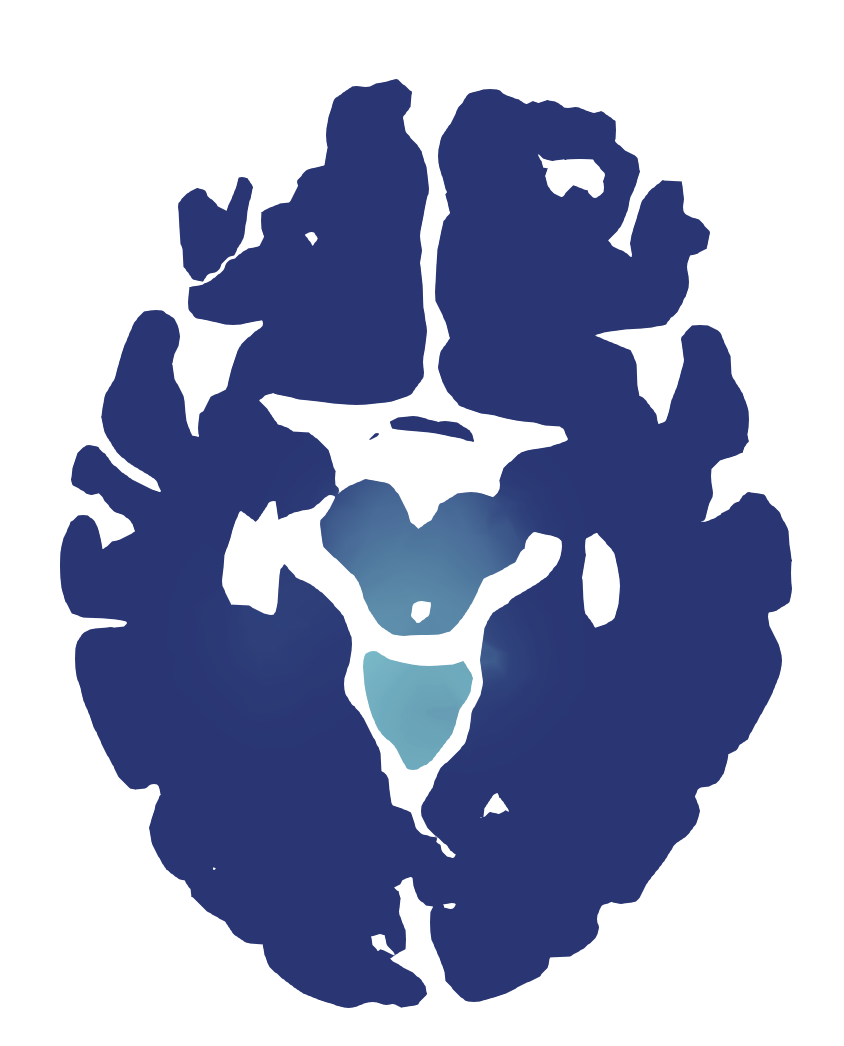}\\
            \includegraphics[width=0.32\textwidth]{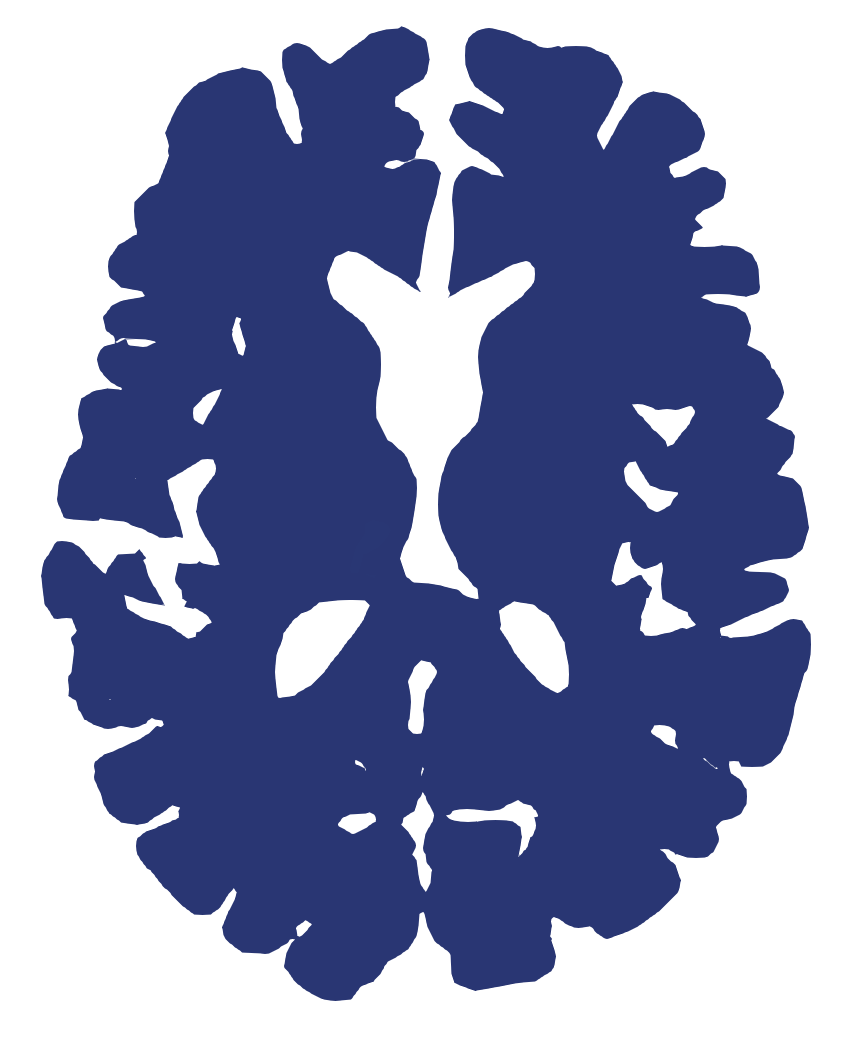}
            \includegraphics[width=0.32\textwidth]{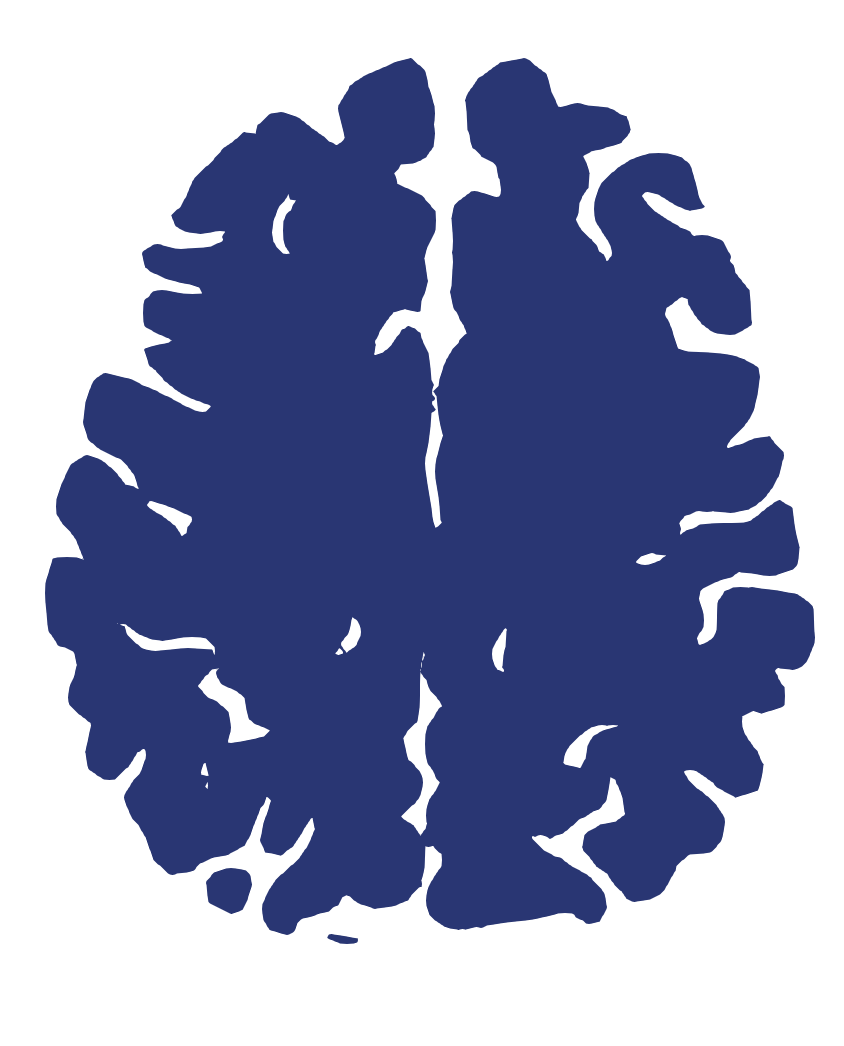}
            \includegraphics[width=0.32\textwidth]{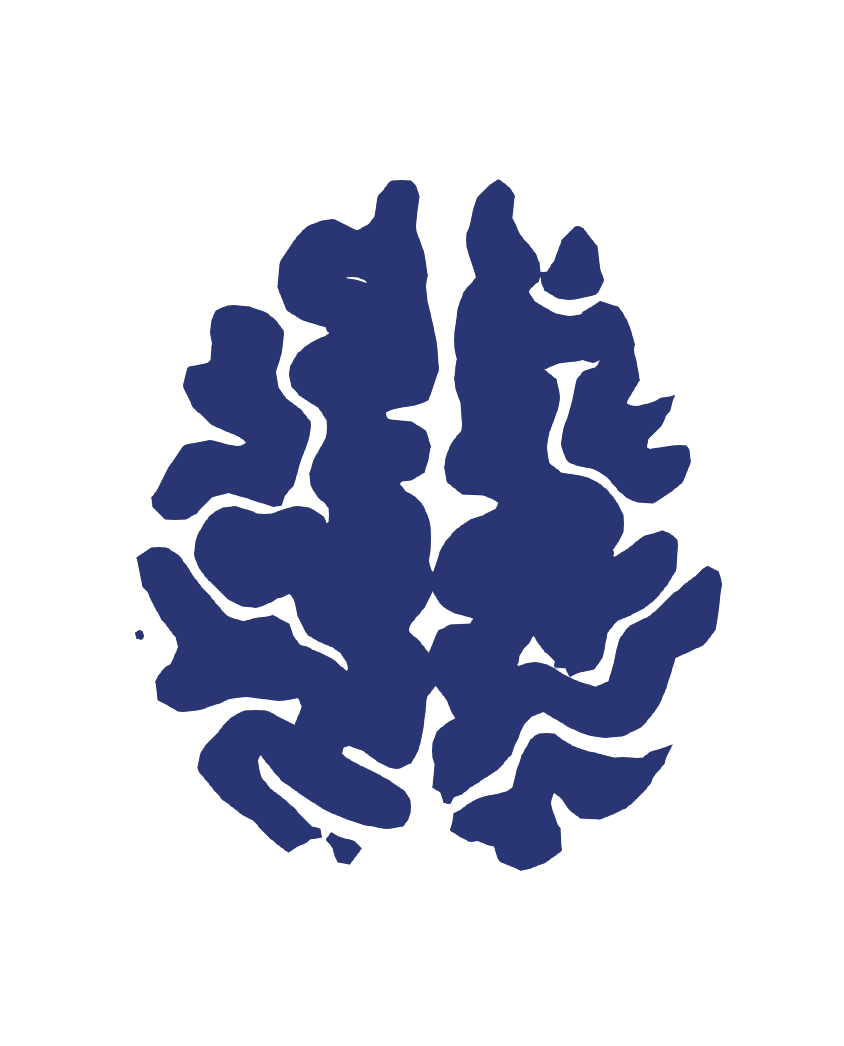}
            \caption{$t=10$ years}
            \label{fig:res2a}
        \end{subfigure}
        \hspace{2em}
        \begin{subfigure}{0.33\textwidth}
            \includegraphics[width=0.32\textwidth]{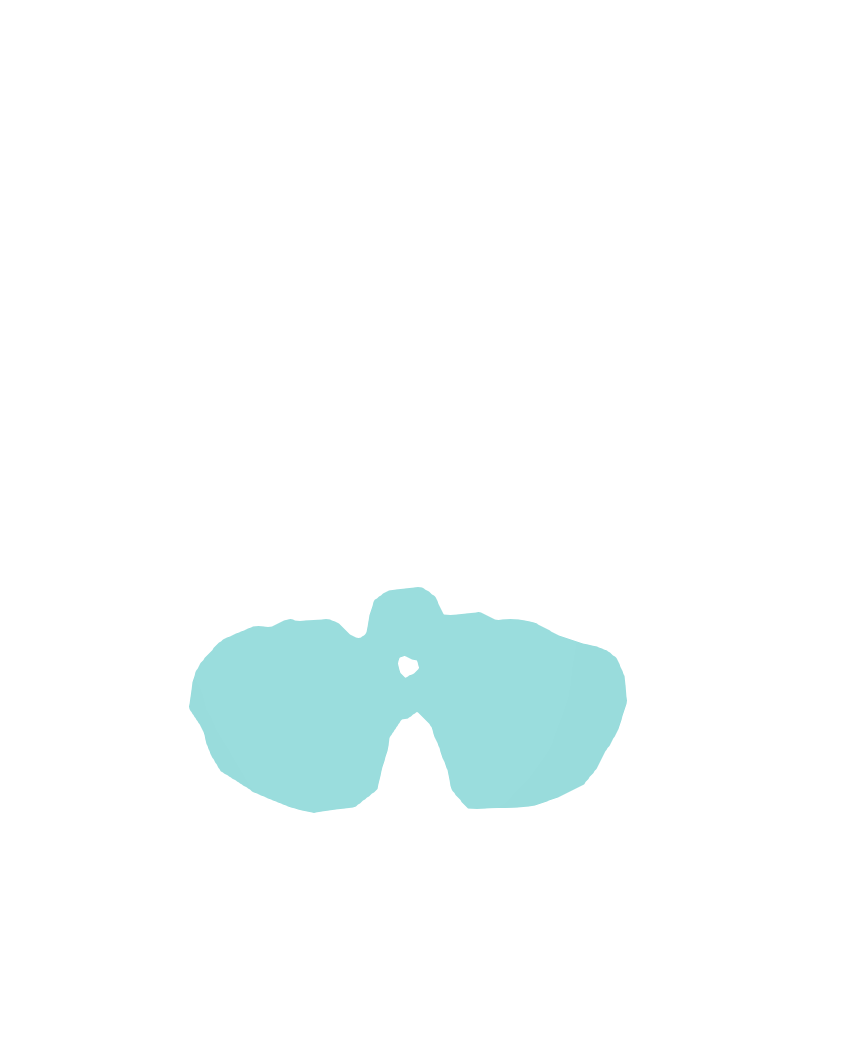}
            \includegraphics[width=0.32\textwidth]{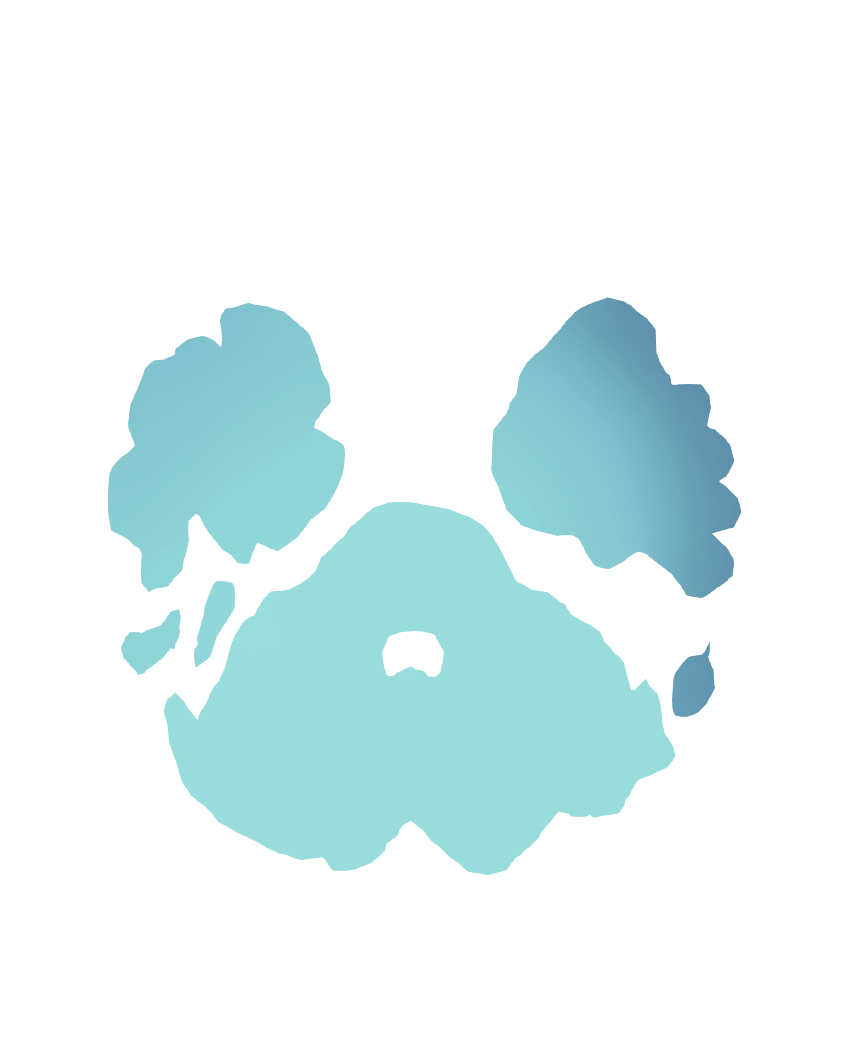}
            \includegraphics[width=0.32\textwidth]{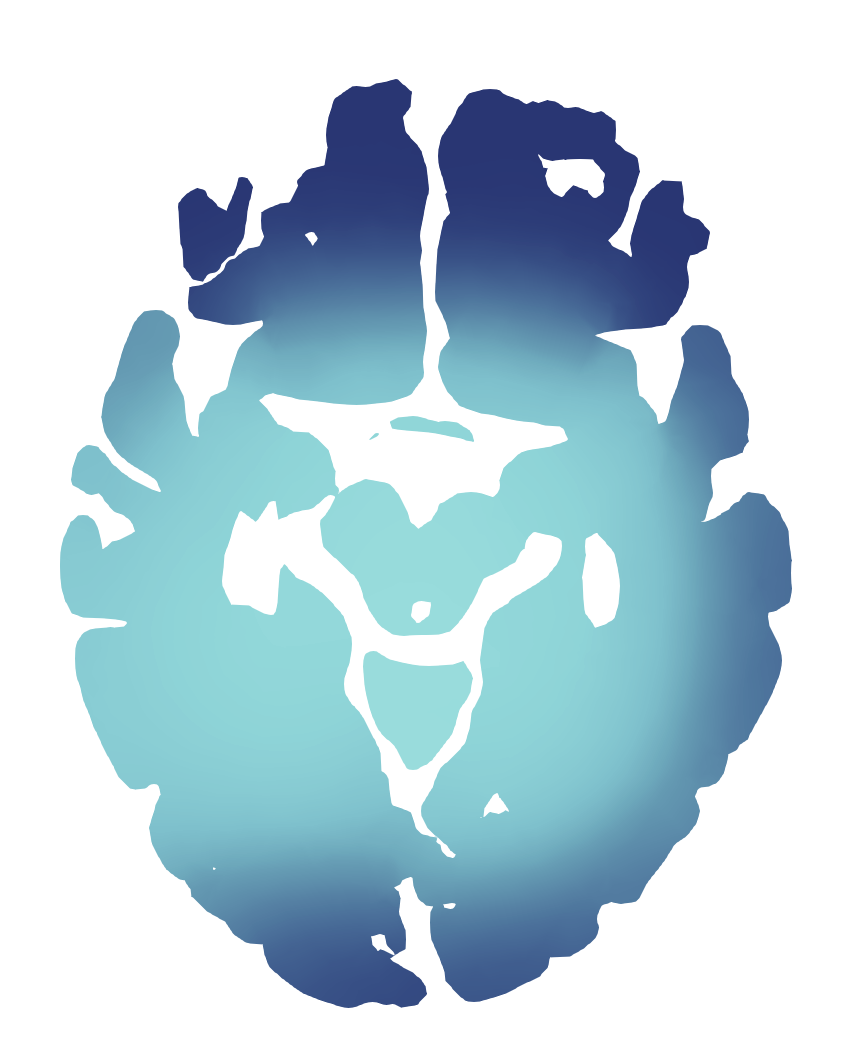}\\
            \includegraphics[width=0.32\textwidth]{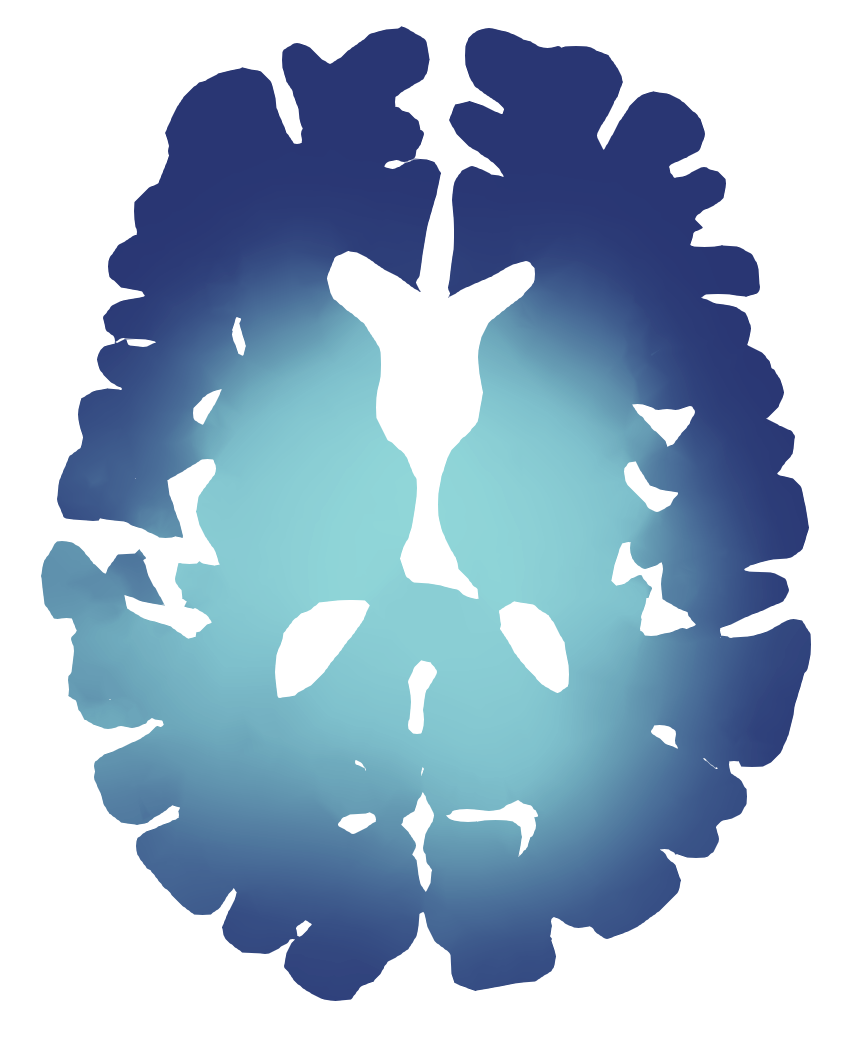}
            \includegraphics[width=0.32\textwidth]{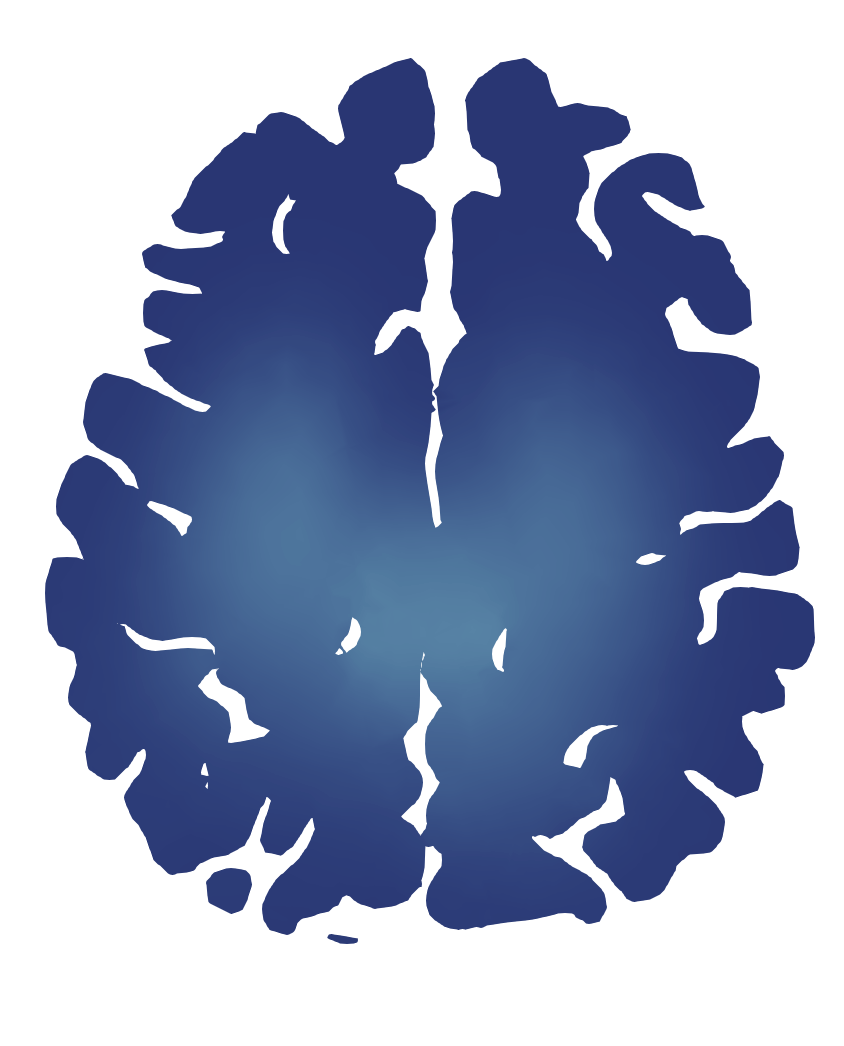}
            \includegraphics[width=0.32\textwidth]{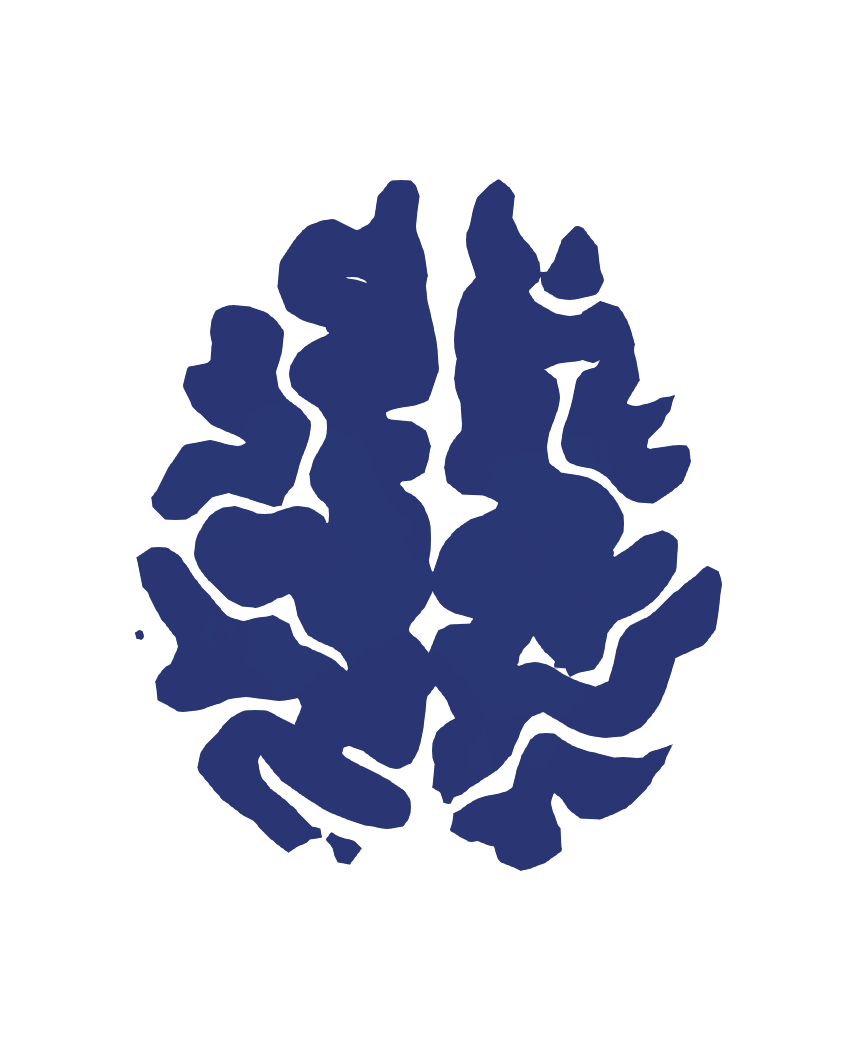}
            \caption{$t=20$ years}
            \label{fig:res2b}
        \end{subfigure}
        \begin{subfigure}{0.12\textwidth}
            \includegraphics[height=4cm]{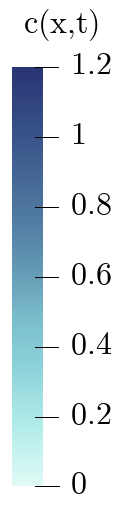}
        \end{subfigure}
        \caption{Test case of Section \ref{sec:results}. Evaluation of healthy protein loss at 10 and 20 years, observed across six transversal sections within the three-dimensional domain shown on the left.}
        \label{fig:res2}
    \end{figure}
    \begin{figure}[H]
        \centering
        \begin{subfigure}{0.12\textwidth}
            \includegraphics[width=0.7\textwidth]{final_simulation/mesh.png}
        \end{subfigure}
        \begin{subfigure}{0.33\textwidth}
            \includegraphics[width=0.32\textwidth]{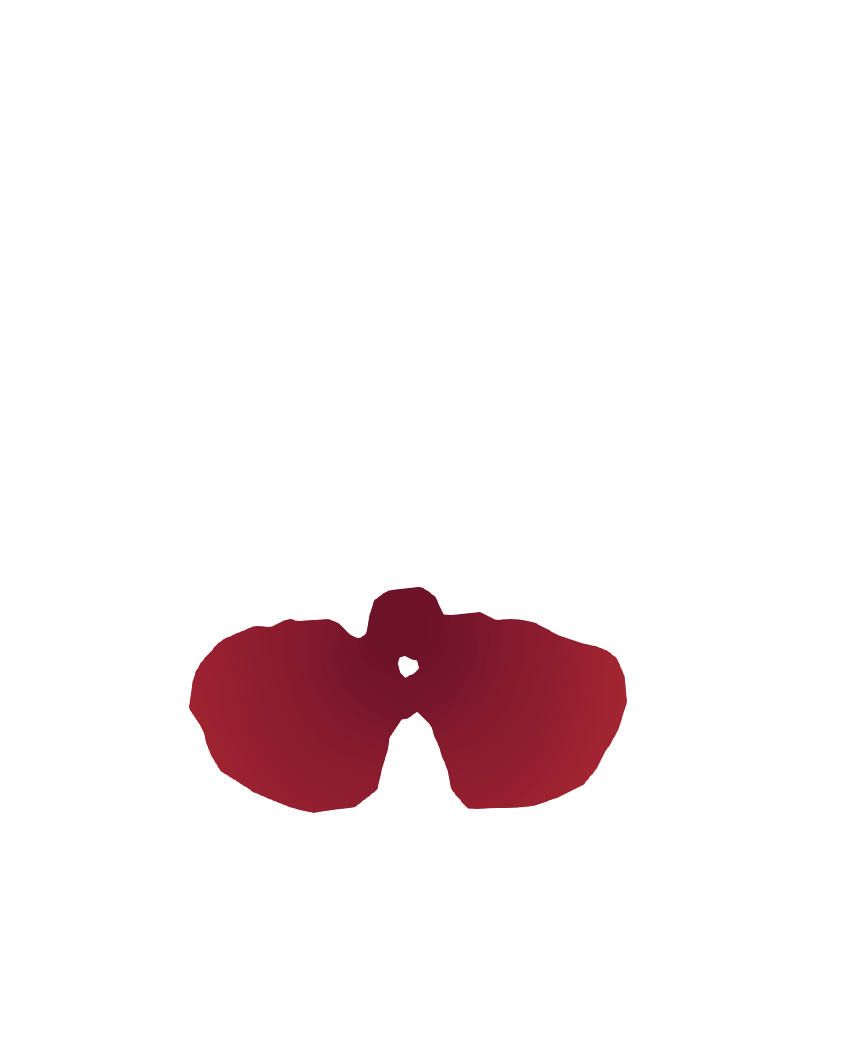}
            \includegraphics[width=0.32\textwidth]{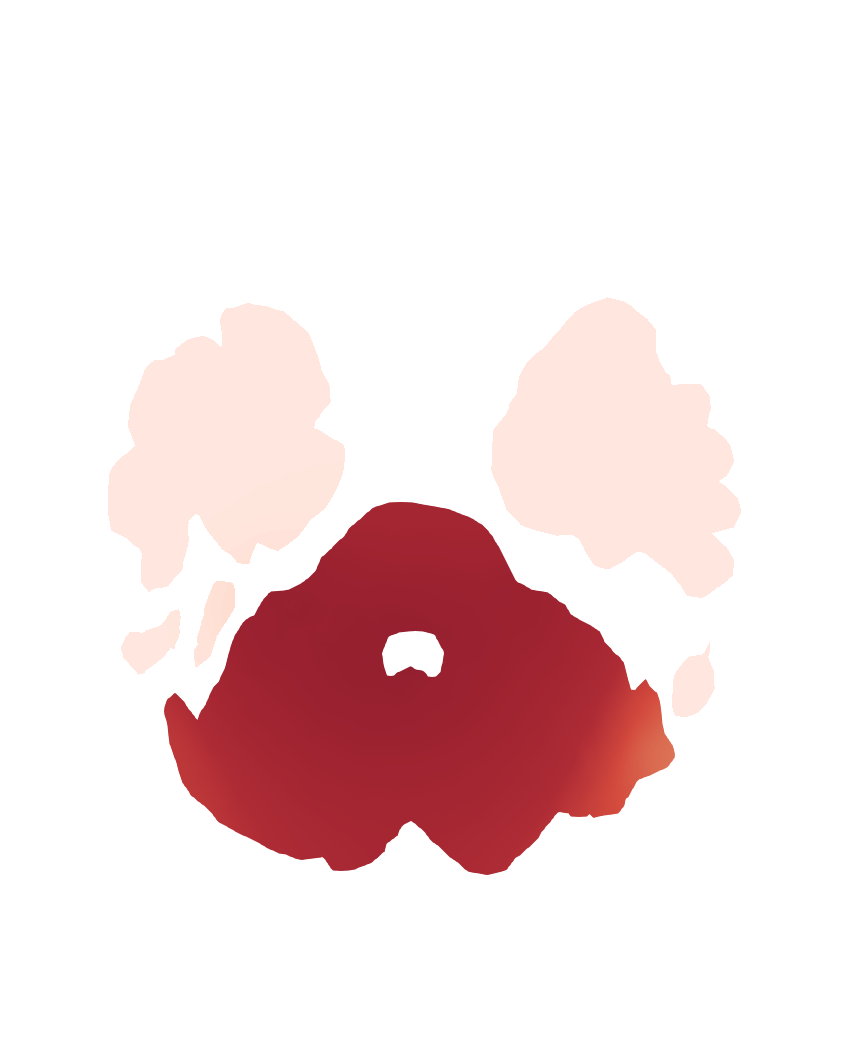}
            \includegraphics[width=0.32\textwidth]{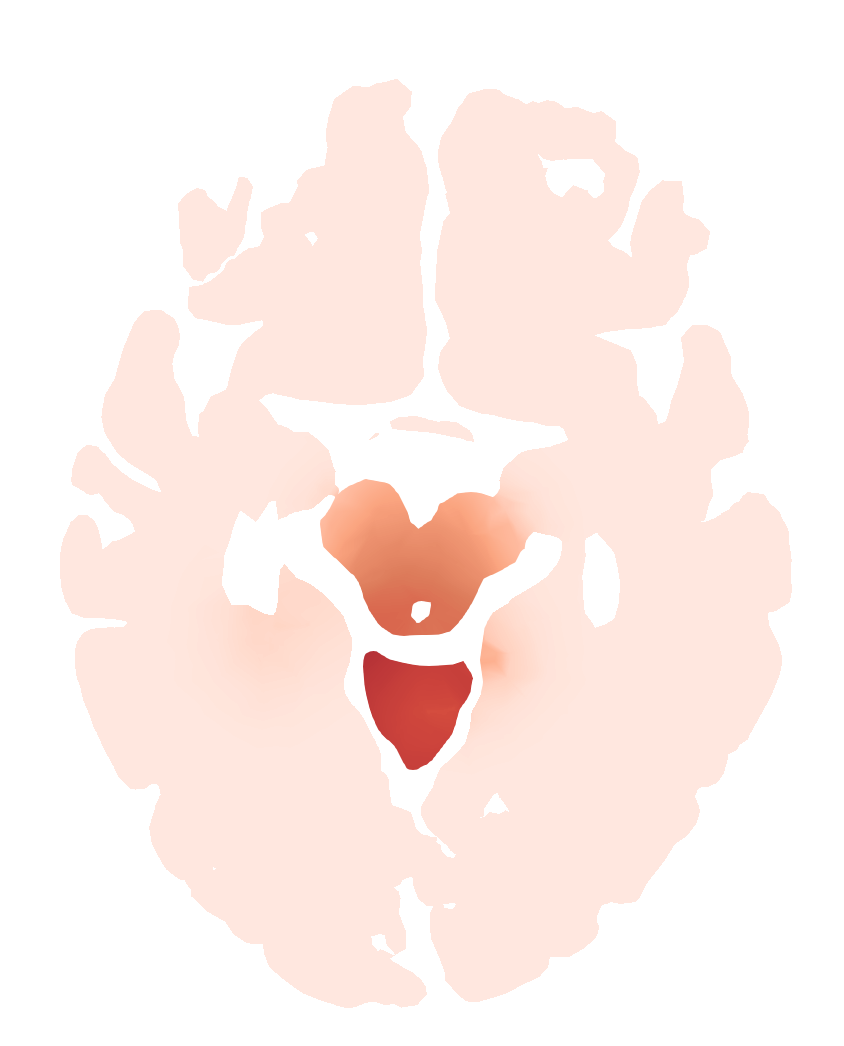}\\
            \includegraphics[width=0.32\textwidth]{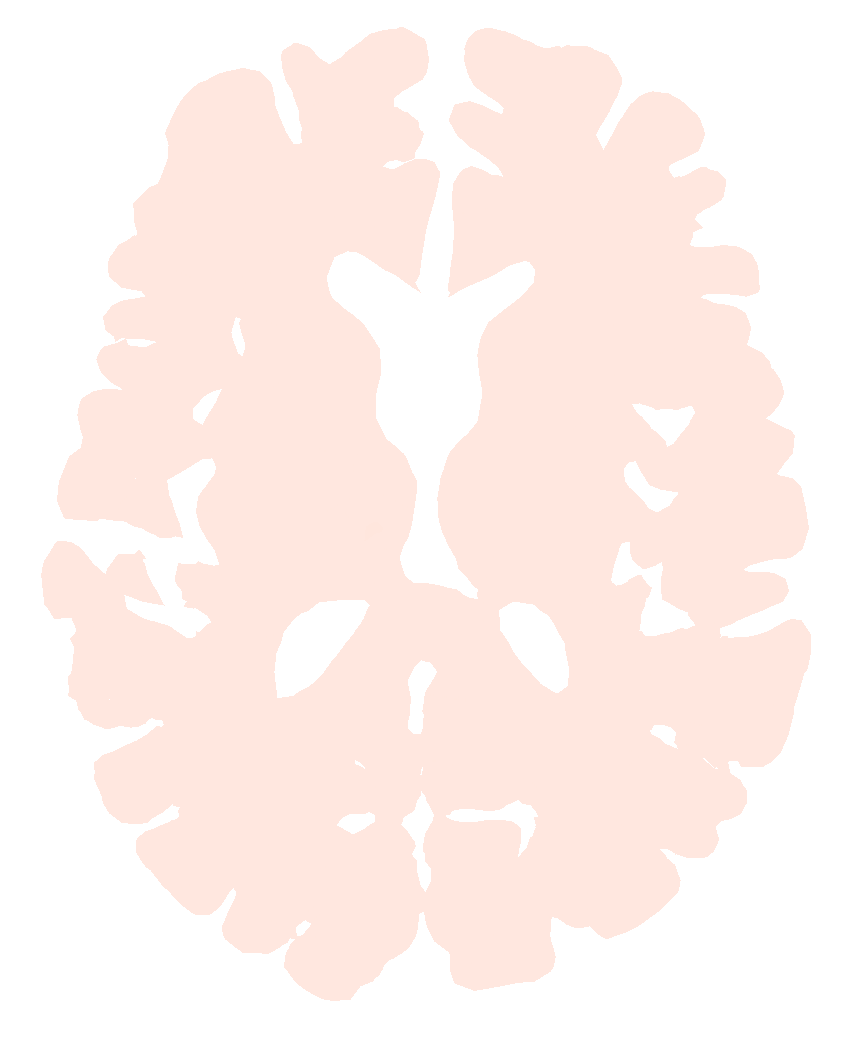}
            \includegraphics[width=0.32\textwidth]{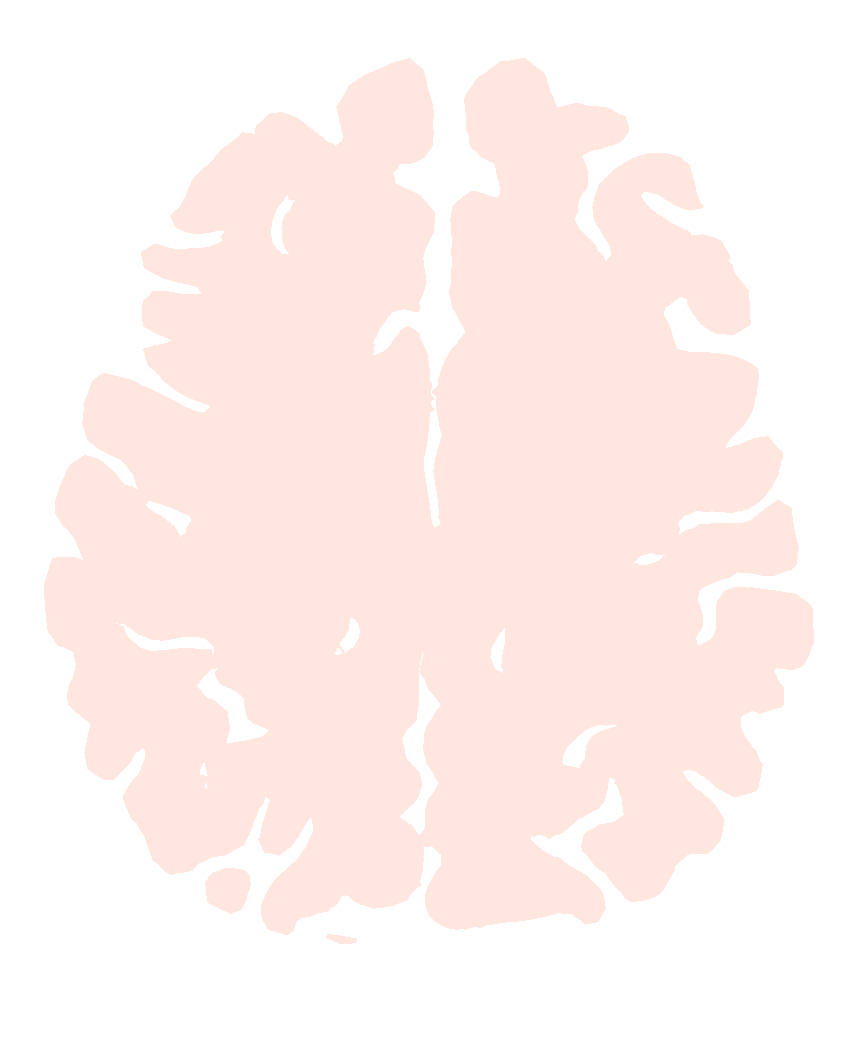}
            \includegraphics[width=0.32\textwidth]{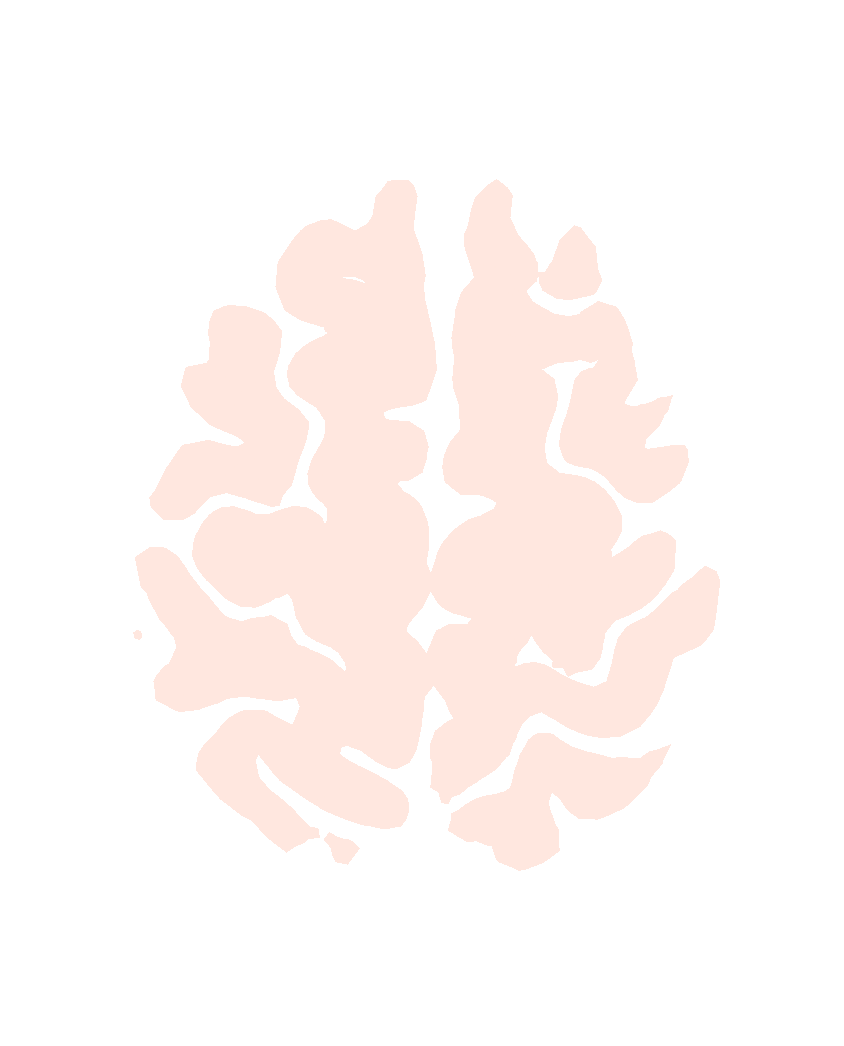}
            \caption{$t=10$ years}
            \label{fig:res1a}
        \end{subfigure}
        \hspace{2em}
        \begin{subfigure}{0.33\textwidth}
            \includegraphics[width=0.32\textwidth]{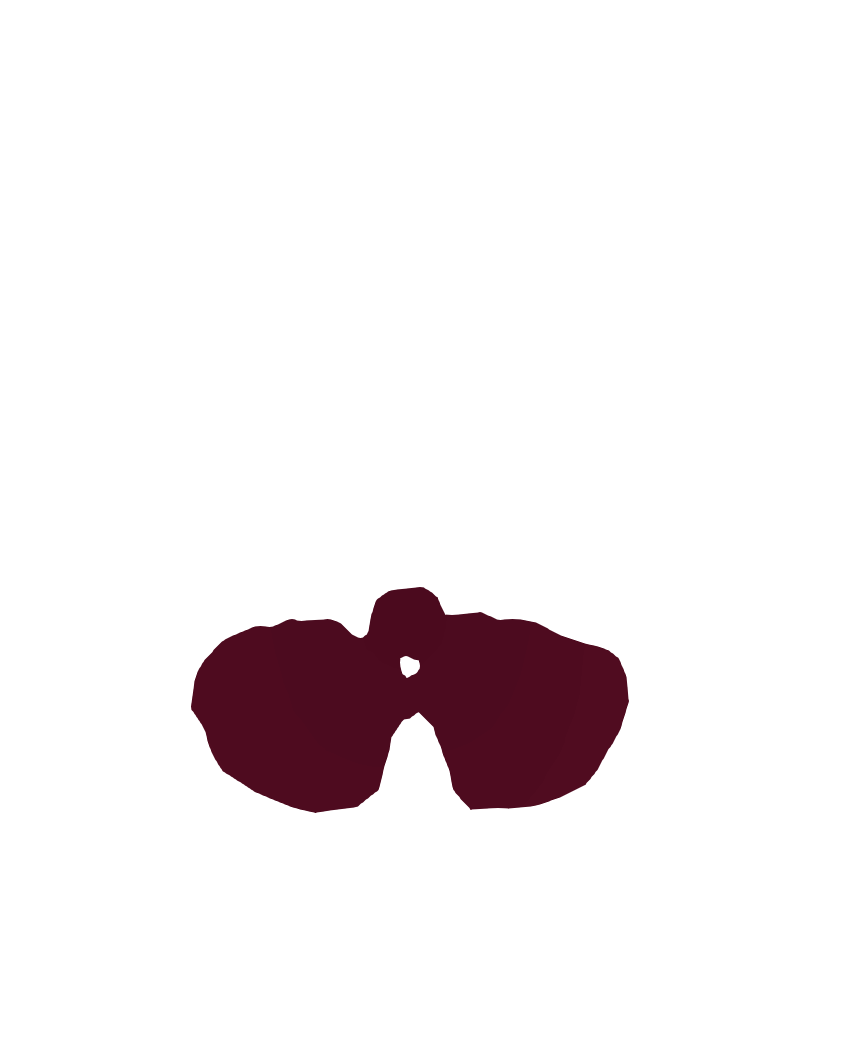}
            \includegraphics[width=0.32\textwidth]{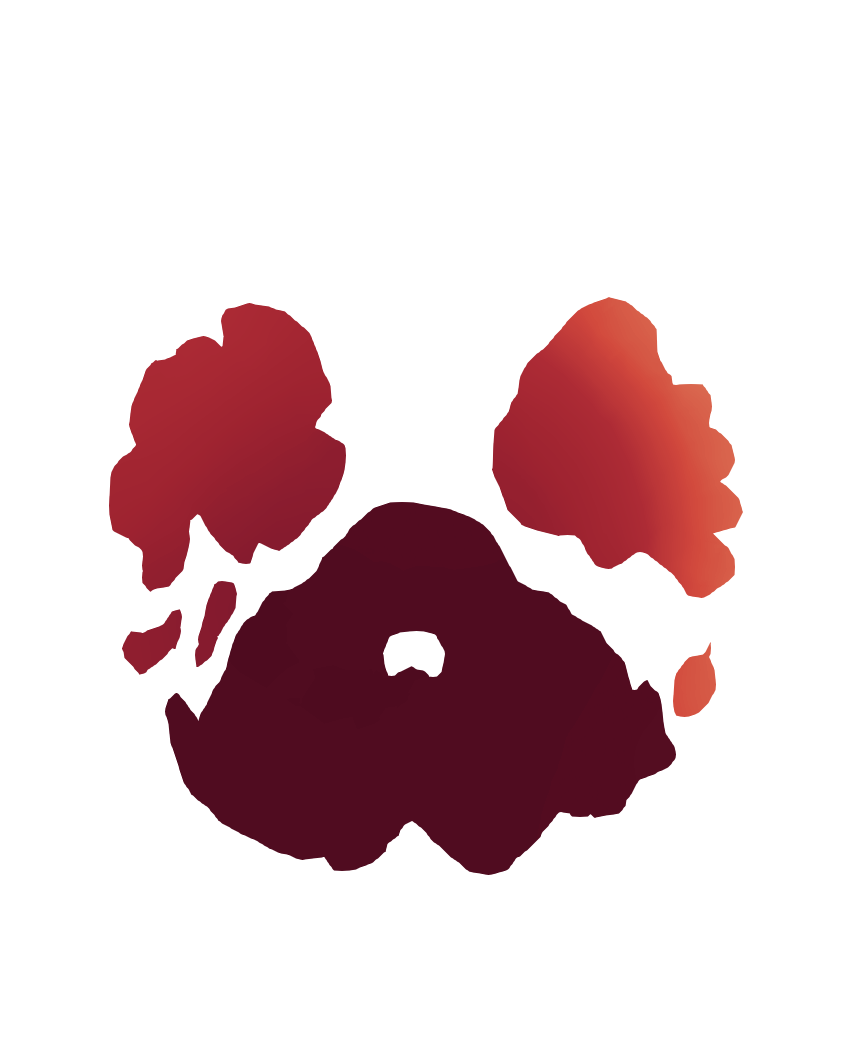}
            \includegraphics[width=0.32\textwidth]{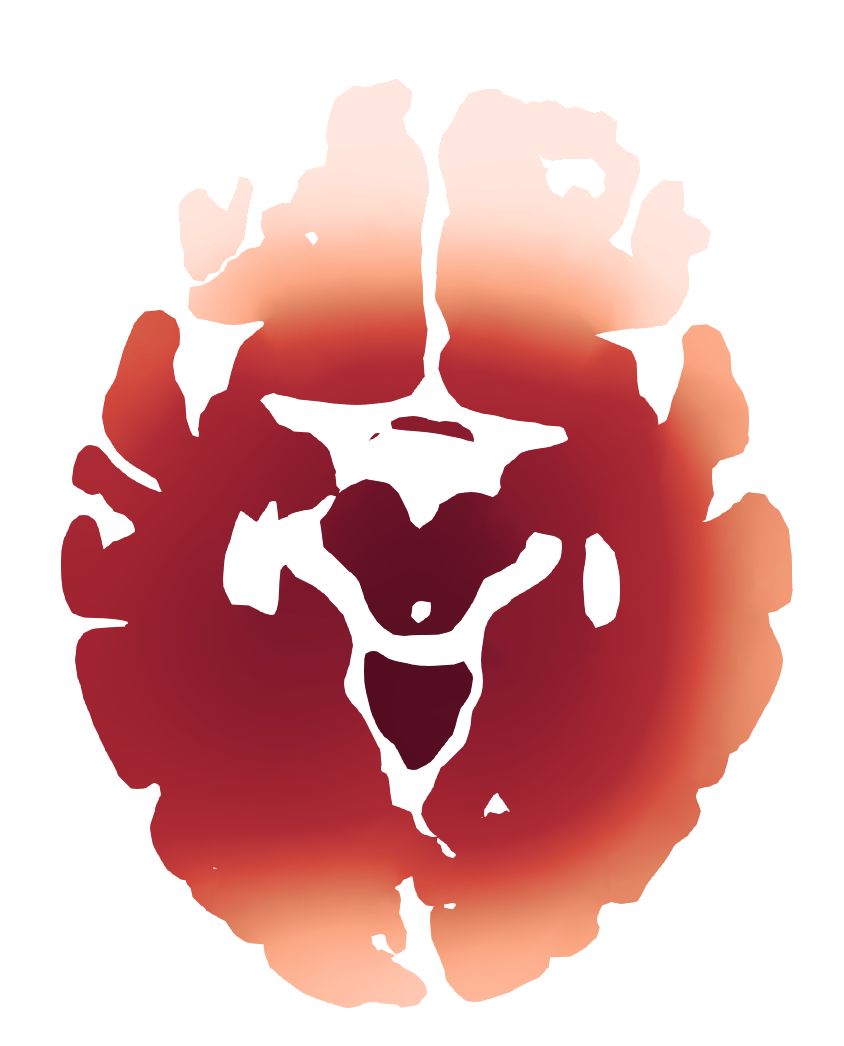}\\
            \includegraphics[width=0.32\textwidth]{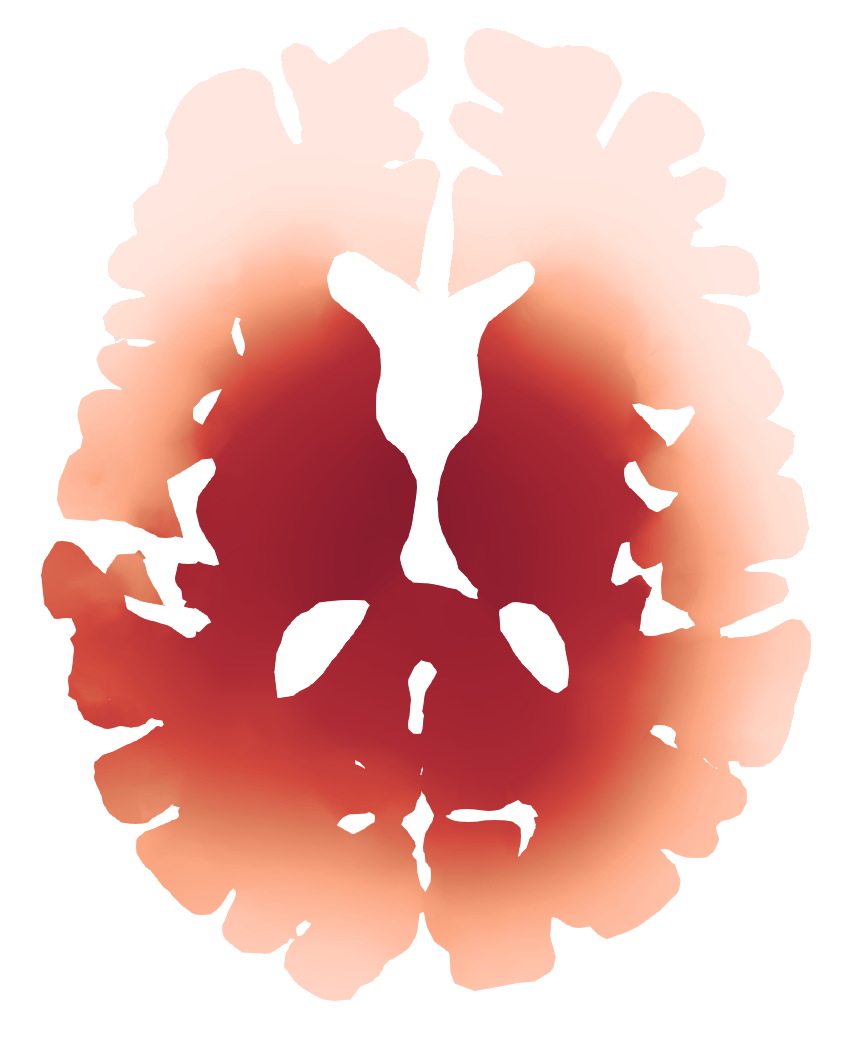}
            \includegraphics[width=0.32\textwidth]{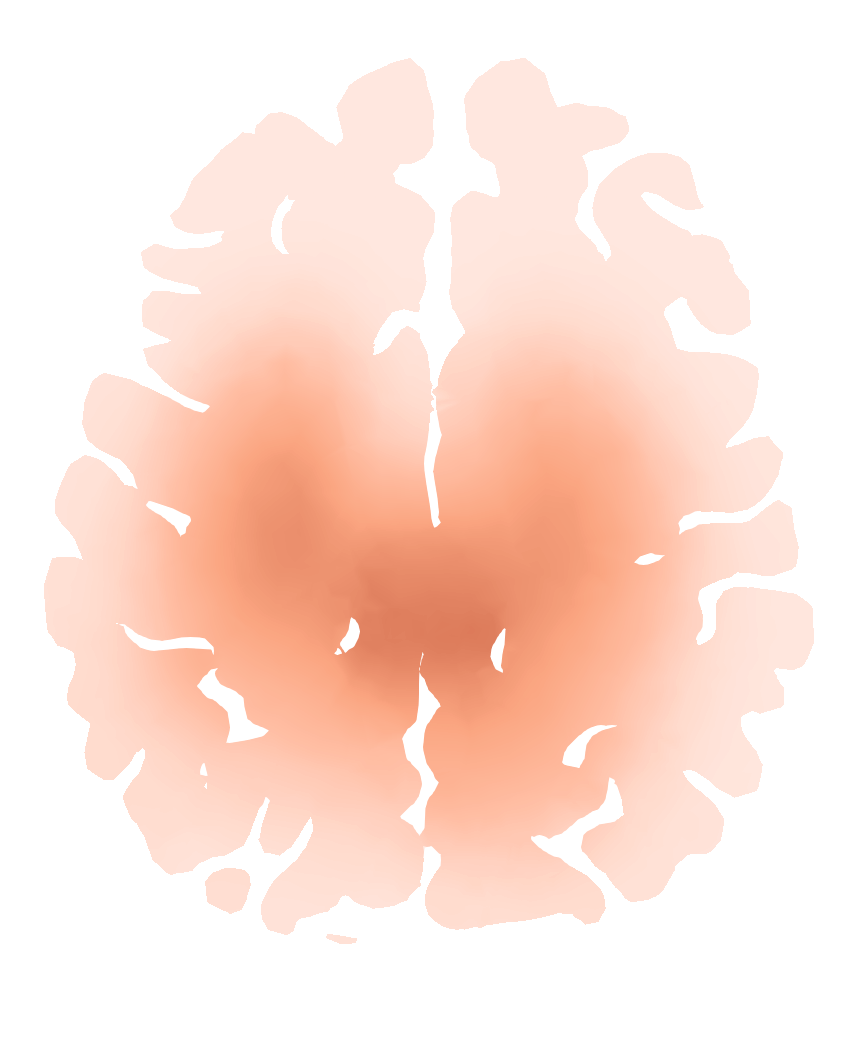}
            \includegraphics[width=0.32\textwidth]{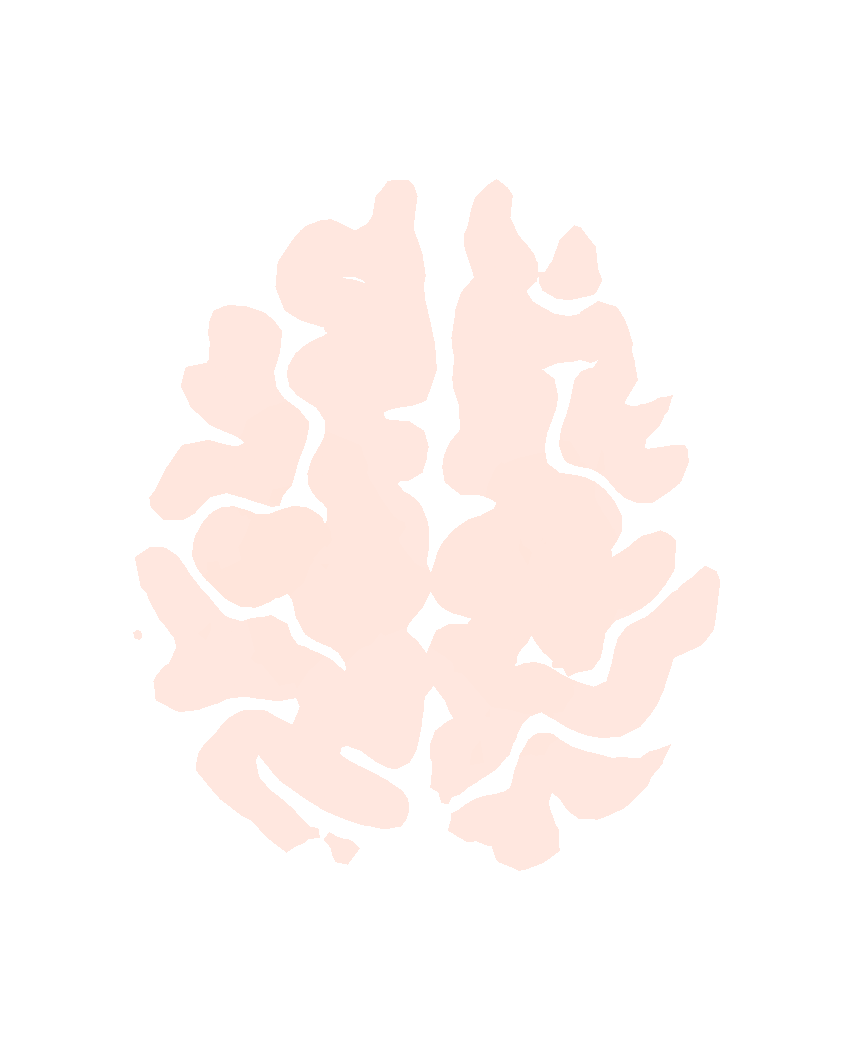}
            \caption{$t=20$ years}
            \label{fig:res1b}
        \end{subfigure}
        \begin{subfigure}{0.12\textwidth}
            \includegraphics[height=4cm]{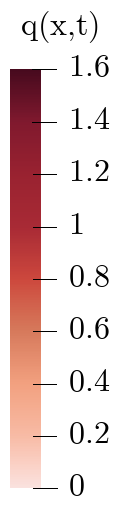}
        \end{subfigure}
        \caption{Test case of Section \ref{sec:results}. Evaluation of misfolded protein diffusion at 10 and 20 years, observed across six transversal sections within the three-dimensional domain shown on the left.}
        \label{fig:res1}
    \end{figure}
    \begin{figure}[H]
        \centering
        \begin{subfigure}{0.12\textwidth}
            \includegraphics[width=0.7\textwidth]{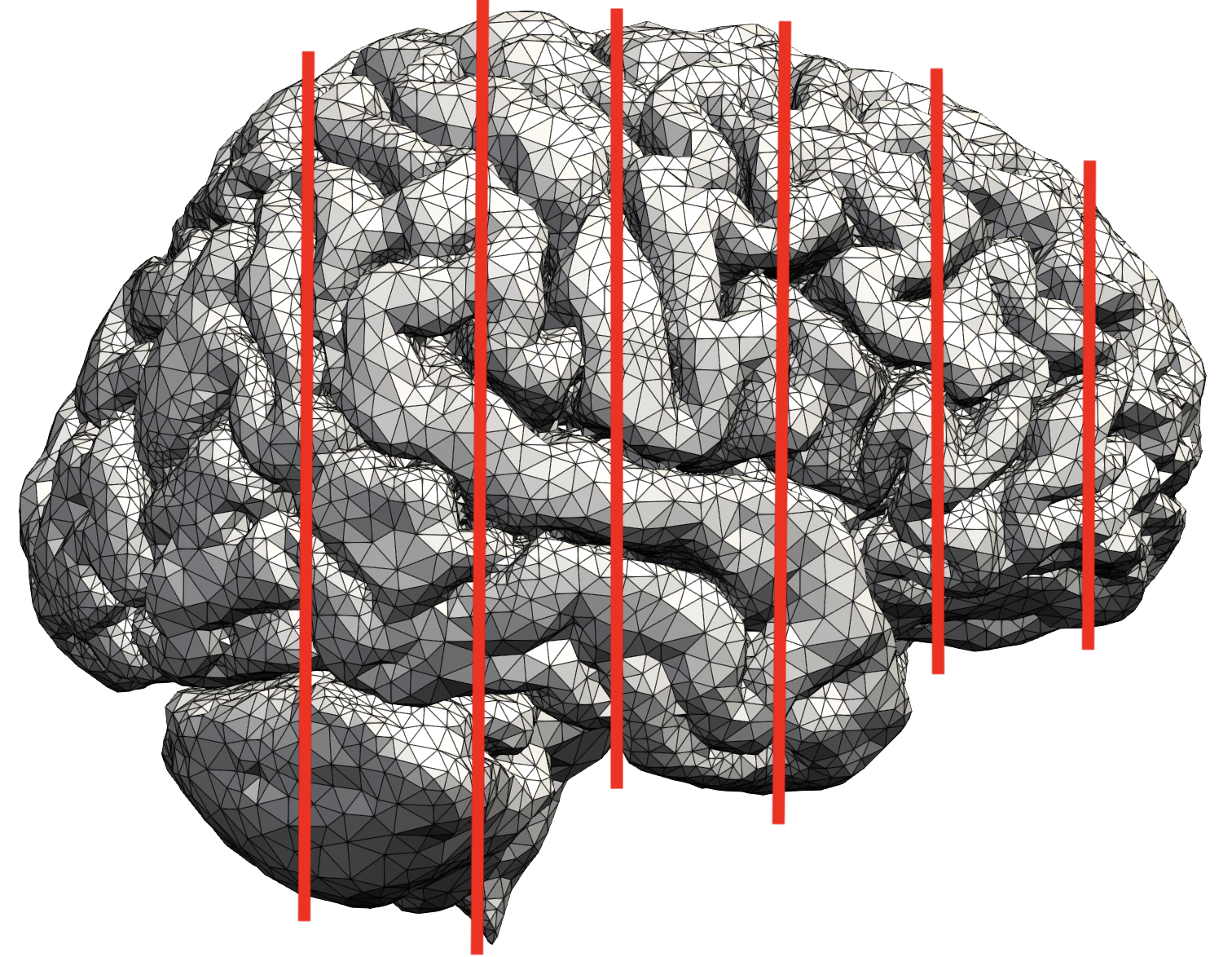}
        \end{subfigure}
        \begin{subfigure}{0.33\textwidth}
            \includegraphics[width=0.32\textwidth]{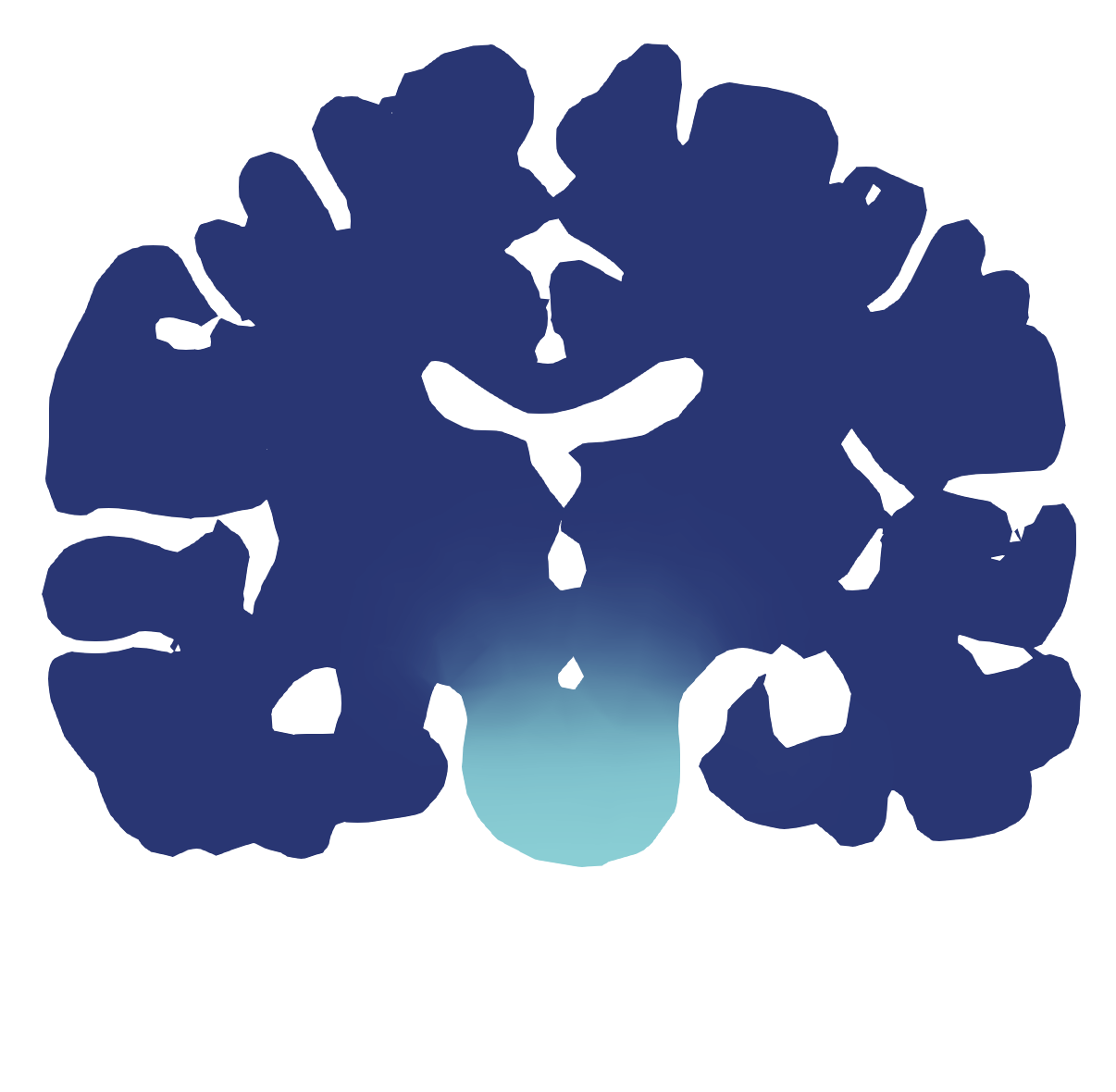}
            \includegraphics[width=0.32\textwidth]{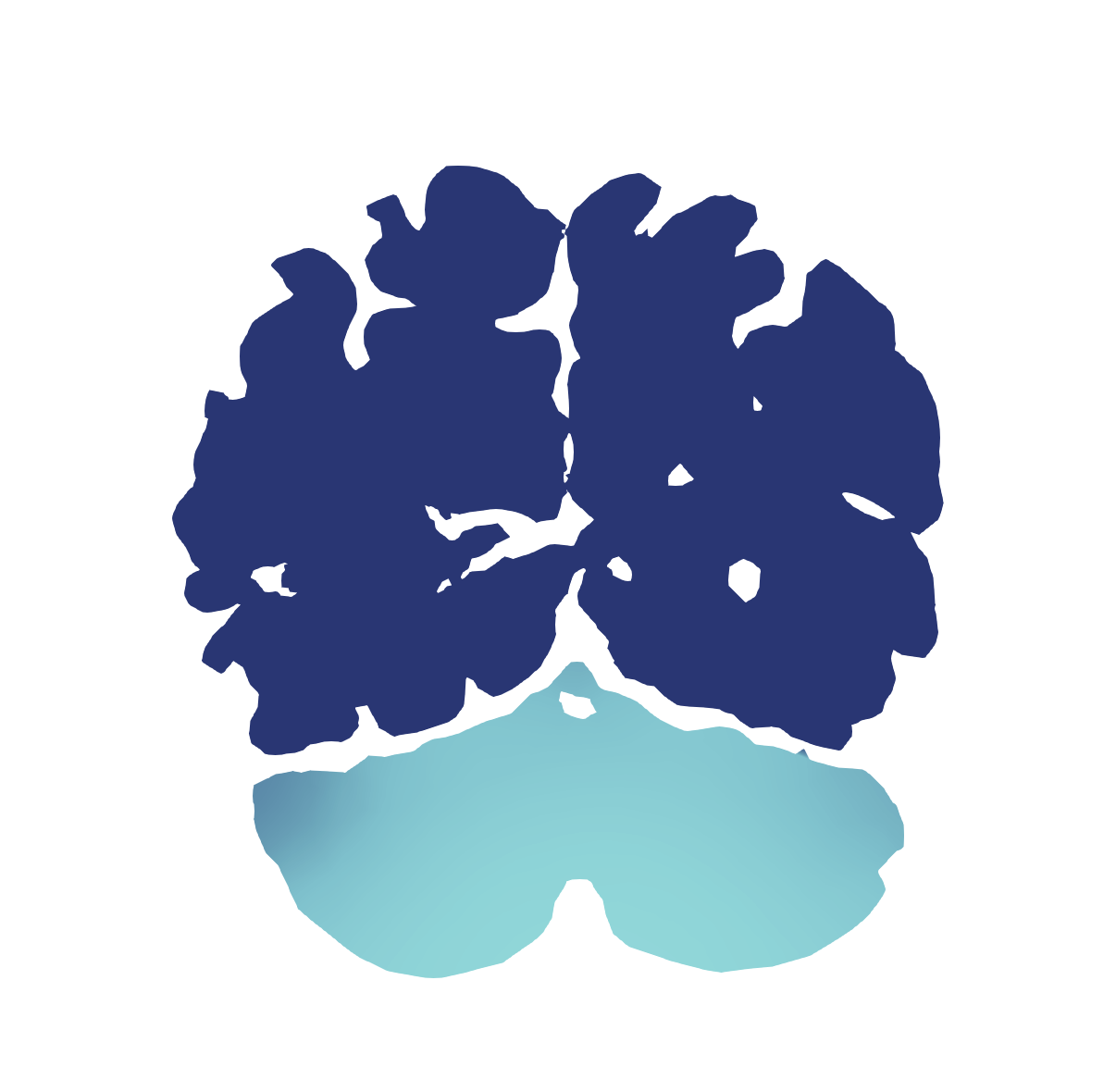}
            \includegraphics[width=0.32\textwidth]{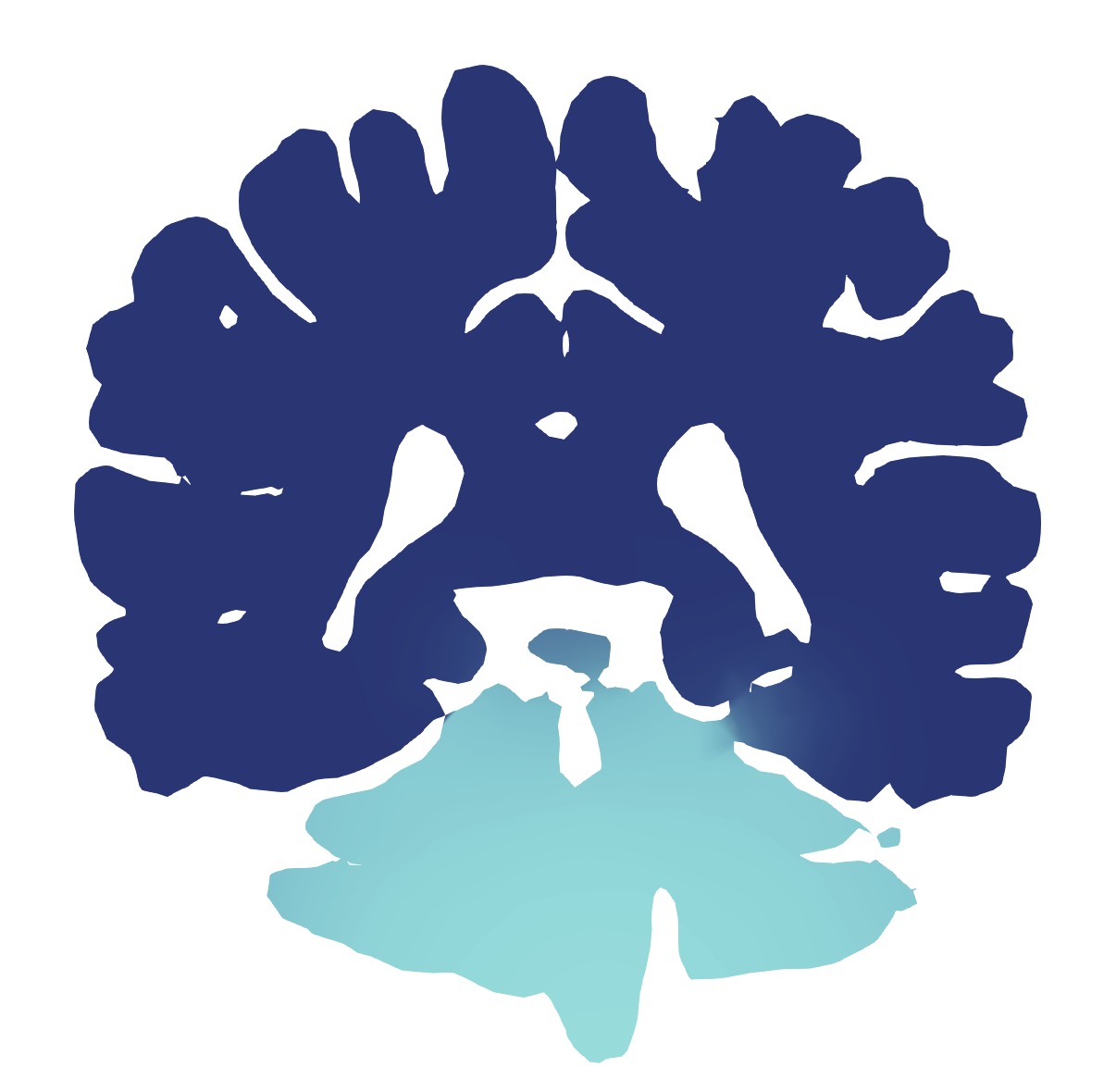}\\
            \includegraphics[width=0.32\textwidth]{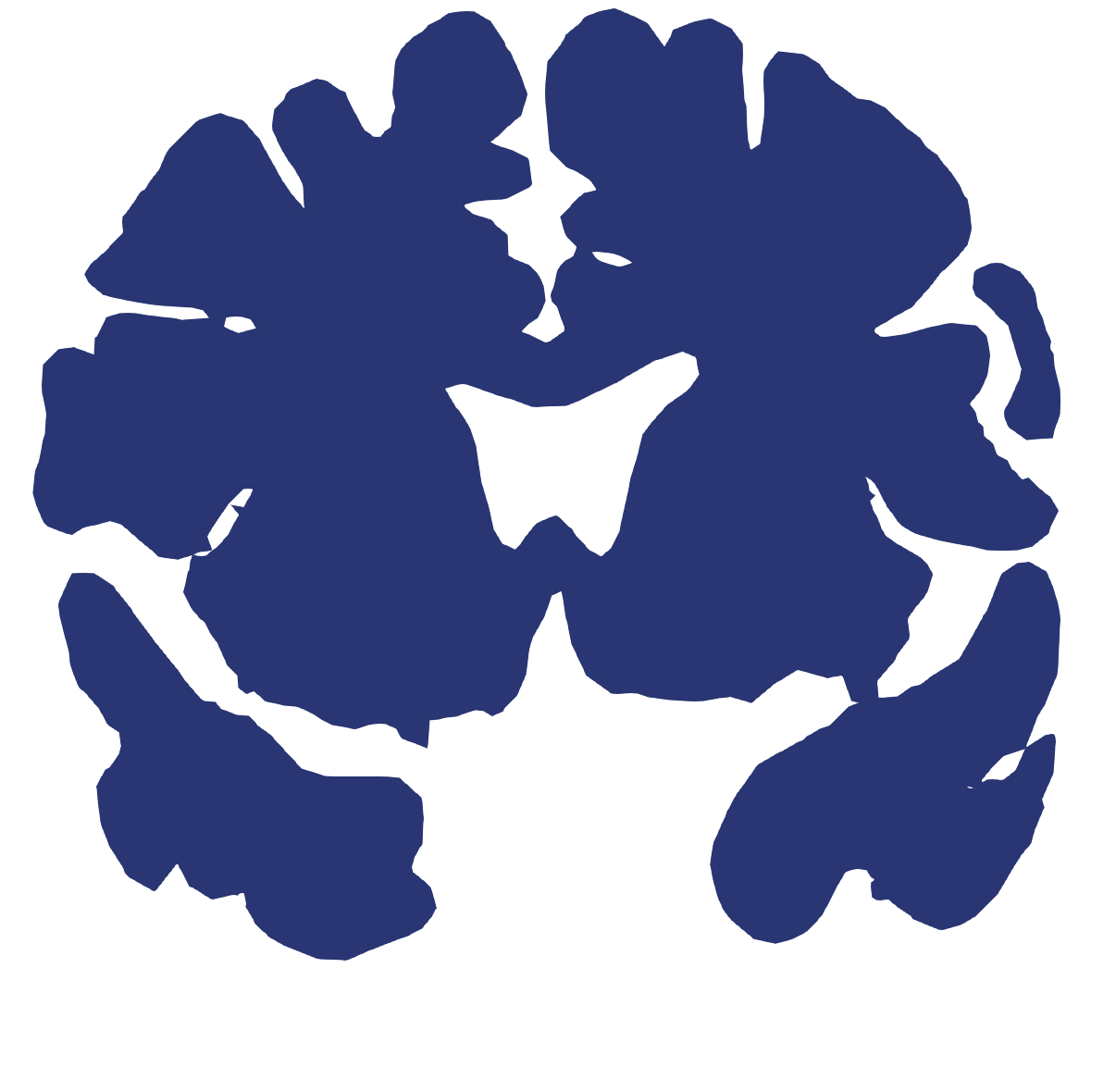}
            \includegraphics[width=0.32\textwidth]{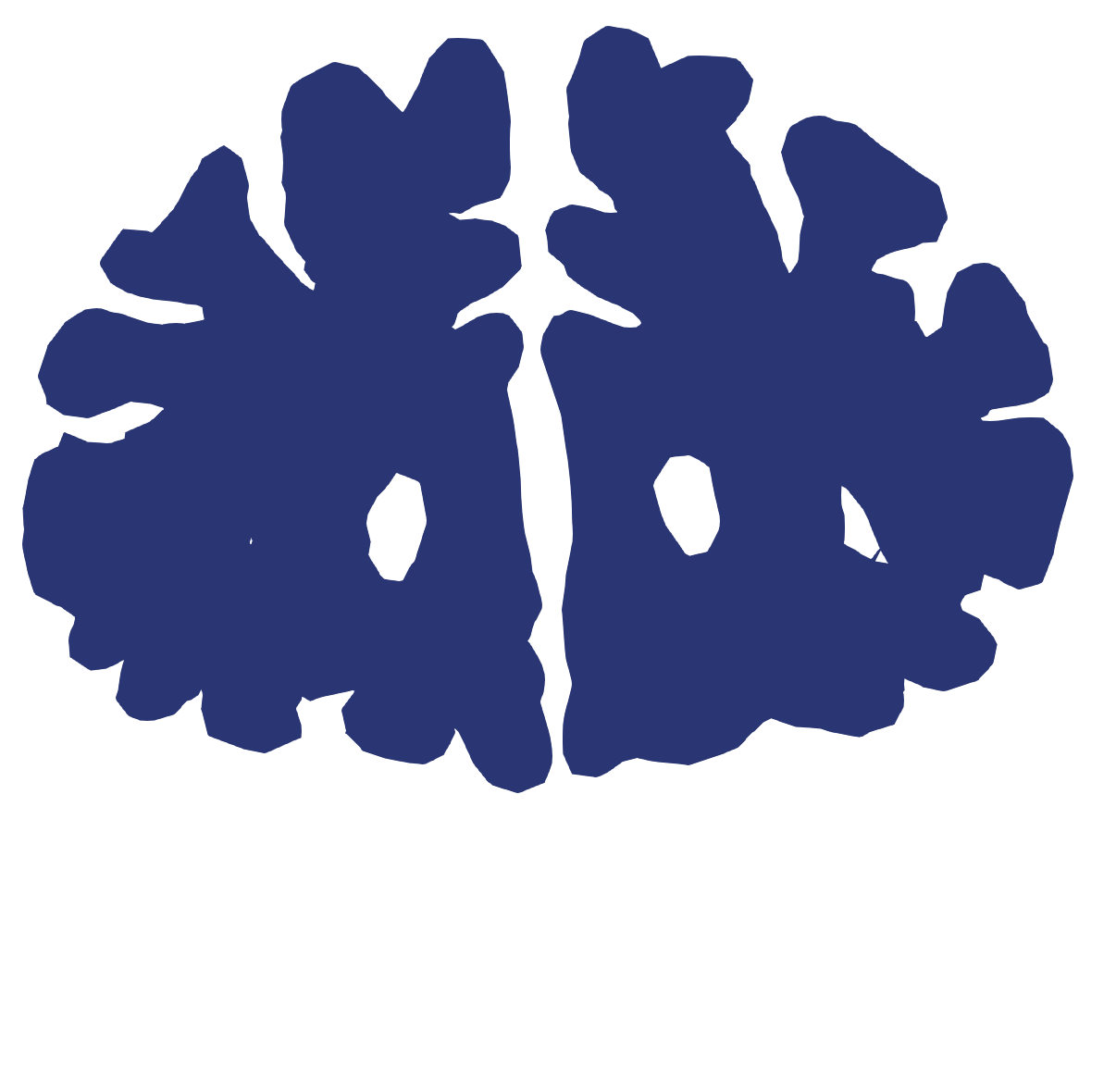}
            \includegraphics[width=0.32\textwidth]{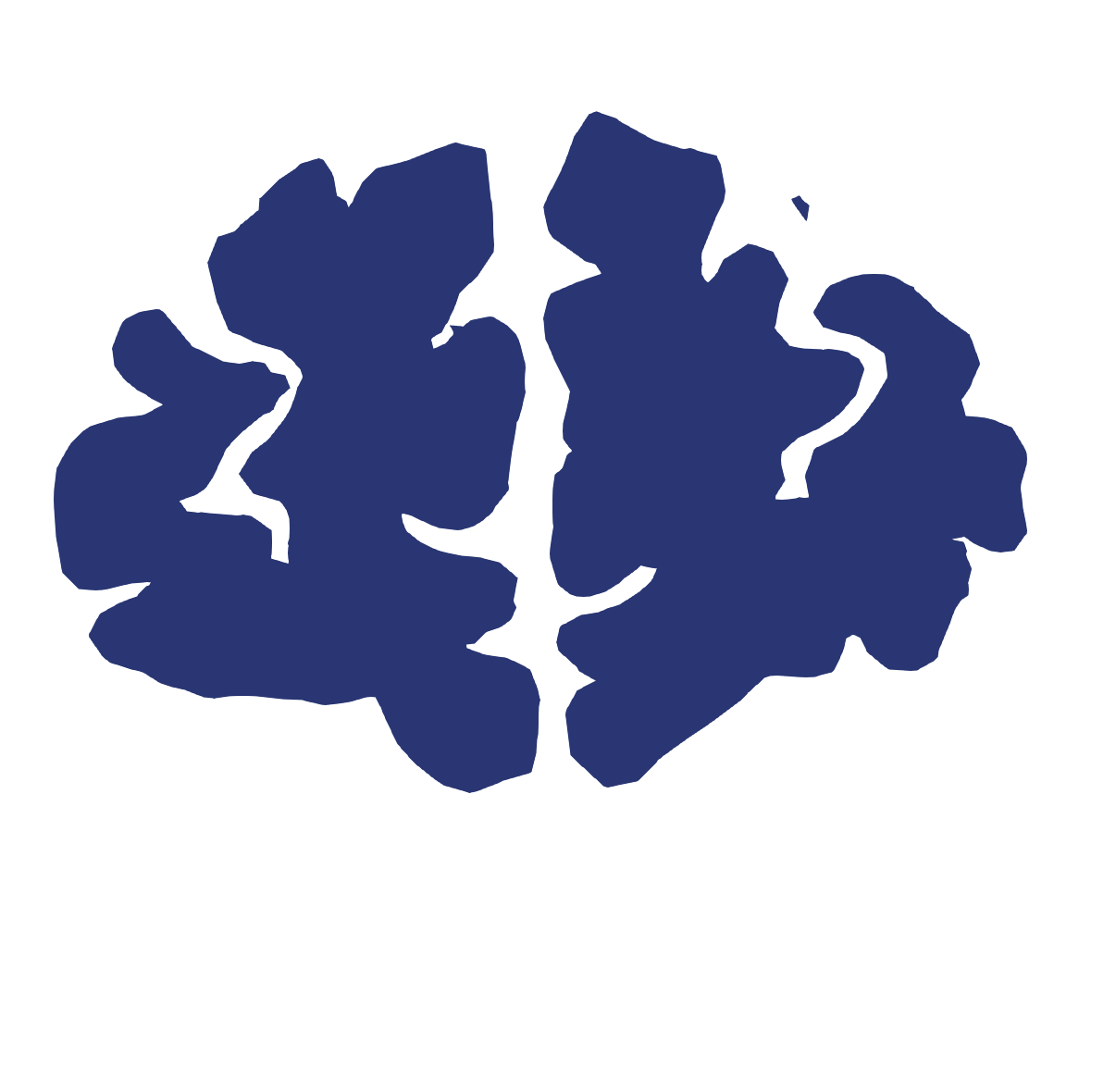}
            \caption{$t=10$ years}
            \label{fig:res3a}
        \end{subfigure}
        \hspace{2em}
        \begin{subfigure}{0.33\textwidth}
            \includegraphics[width=0.32\textwidth]{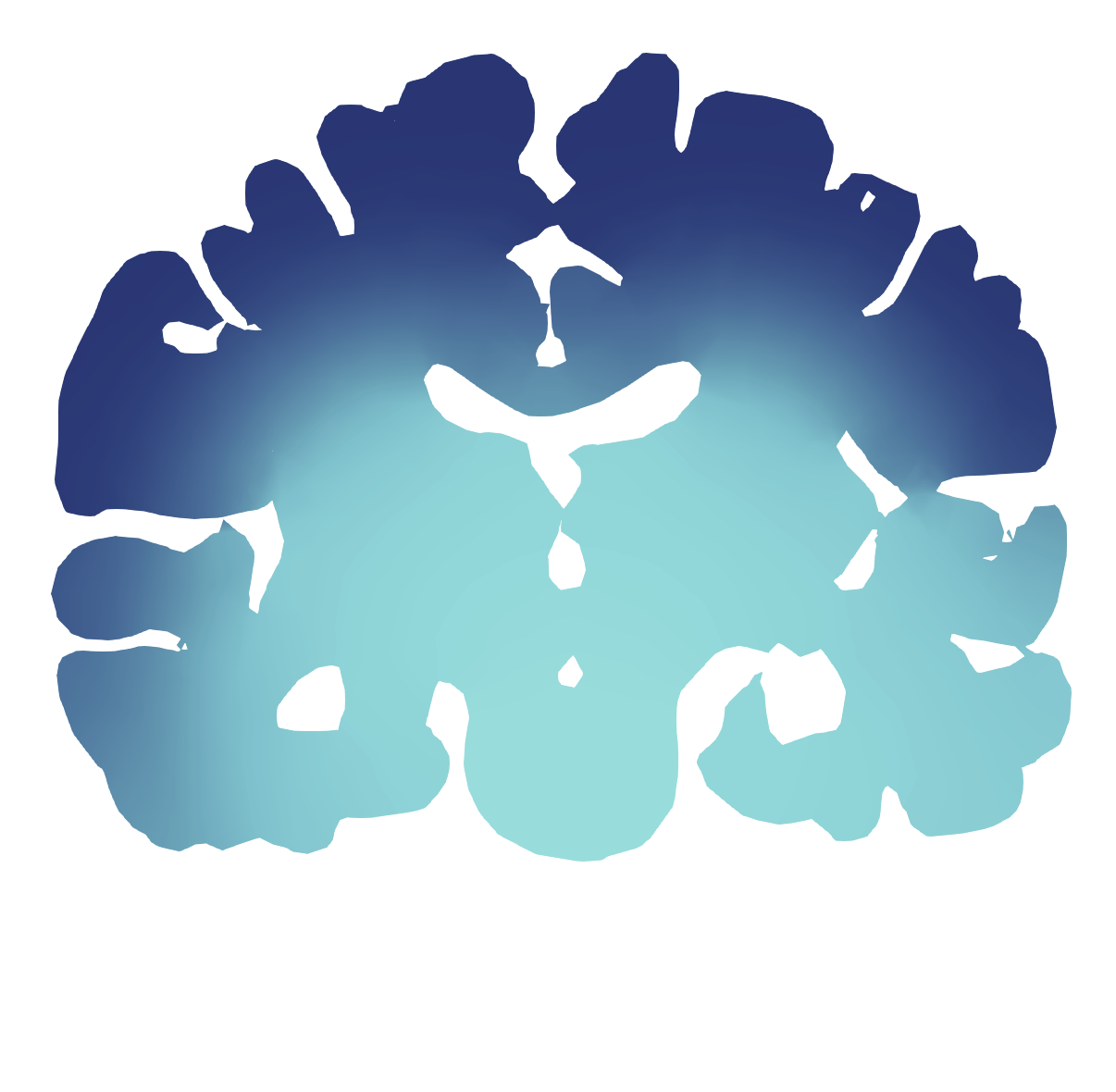}
            \includegraphics[width=0.32\textwidth]{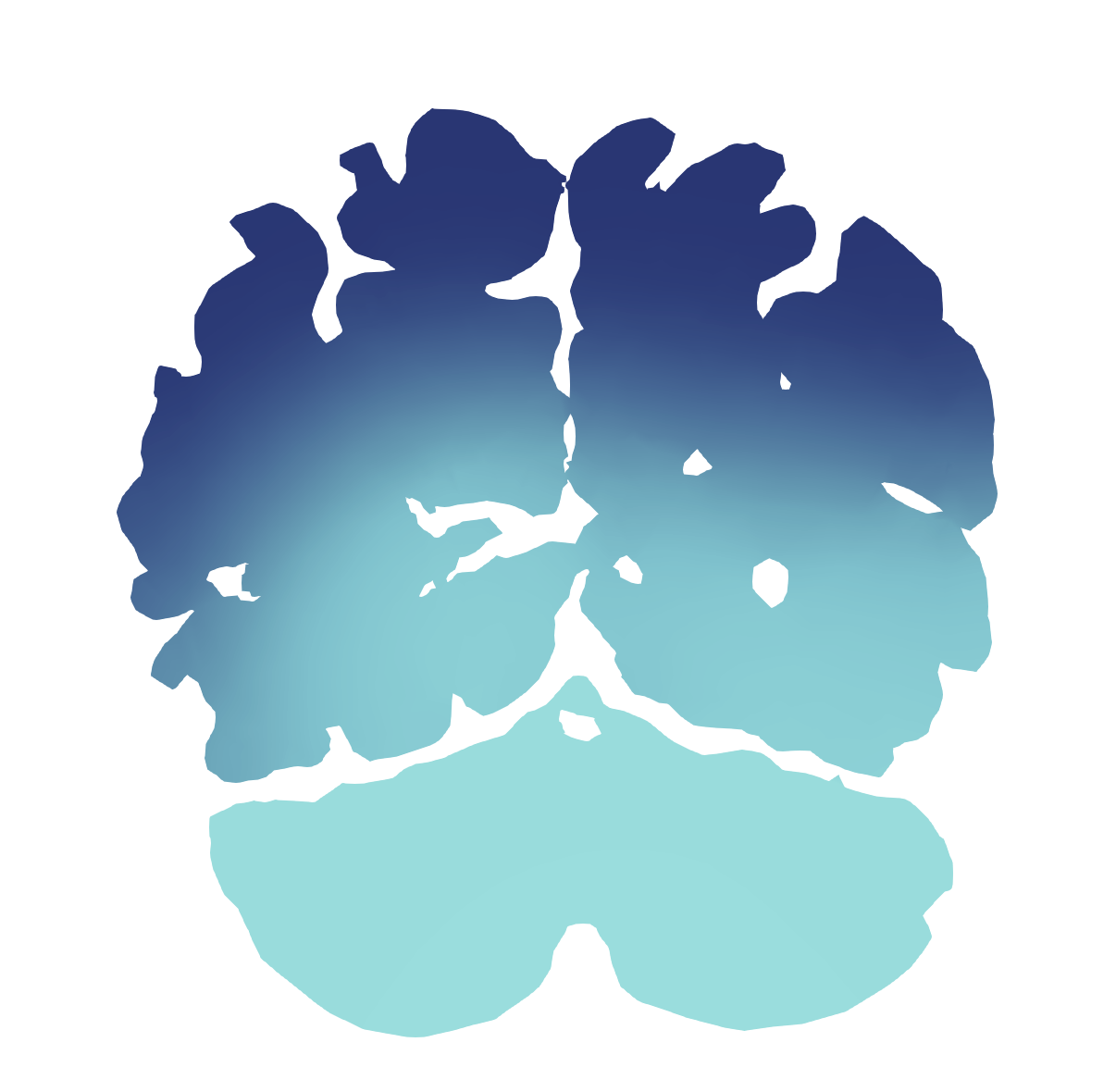}
            \includegraphics[width=0.32\textwidth]{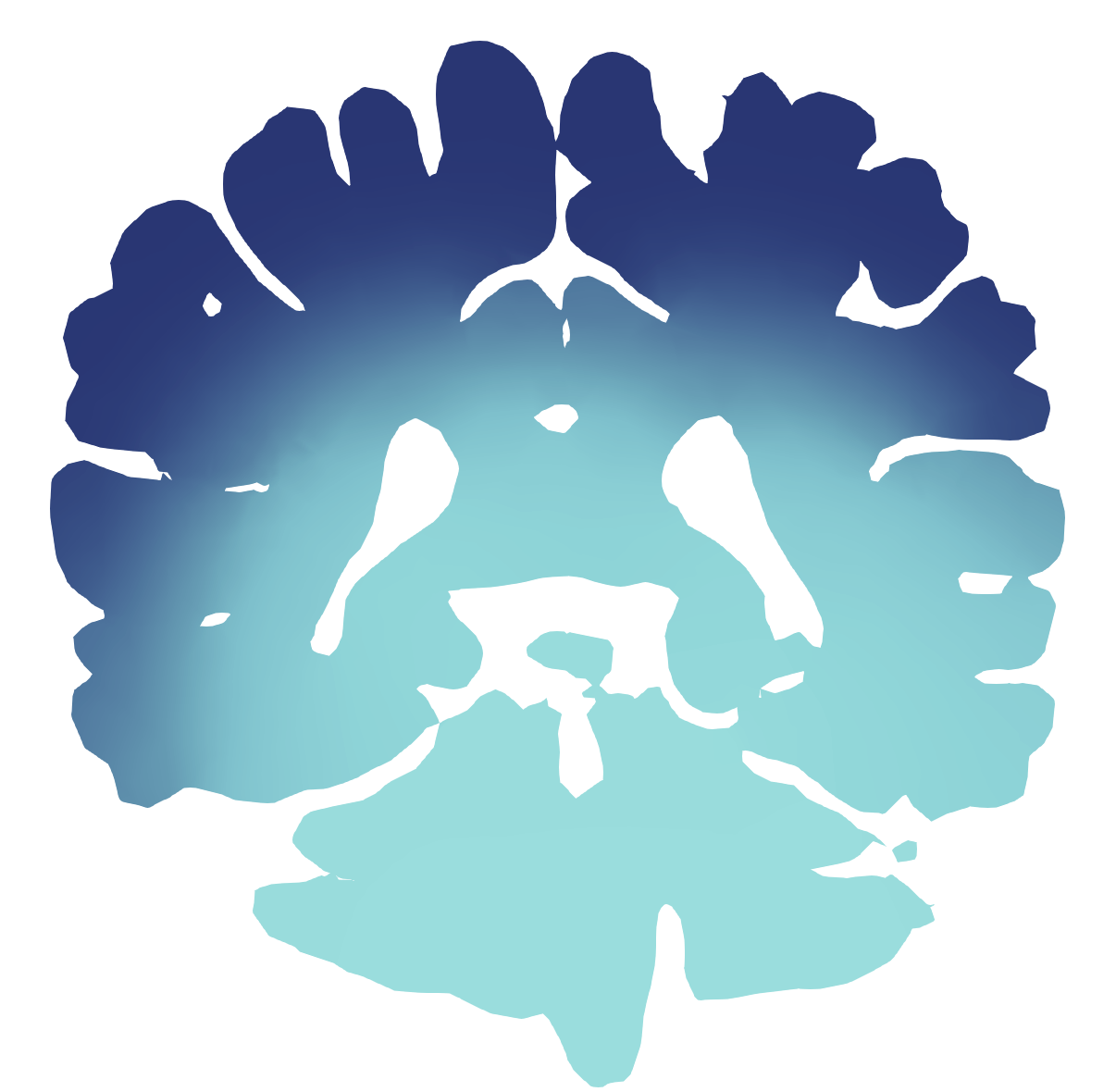}\\
            \includegraphics[width=0.32\textwidth]{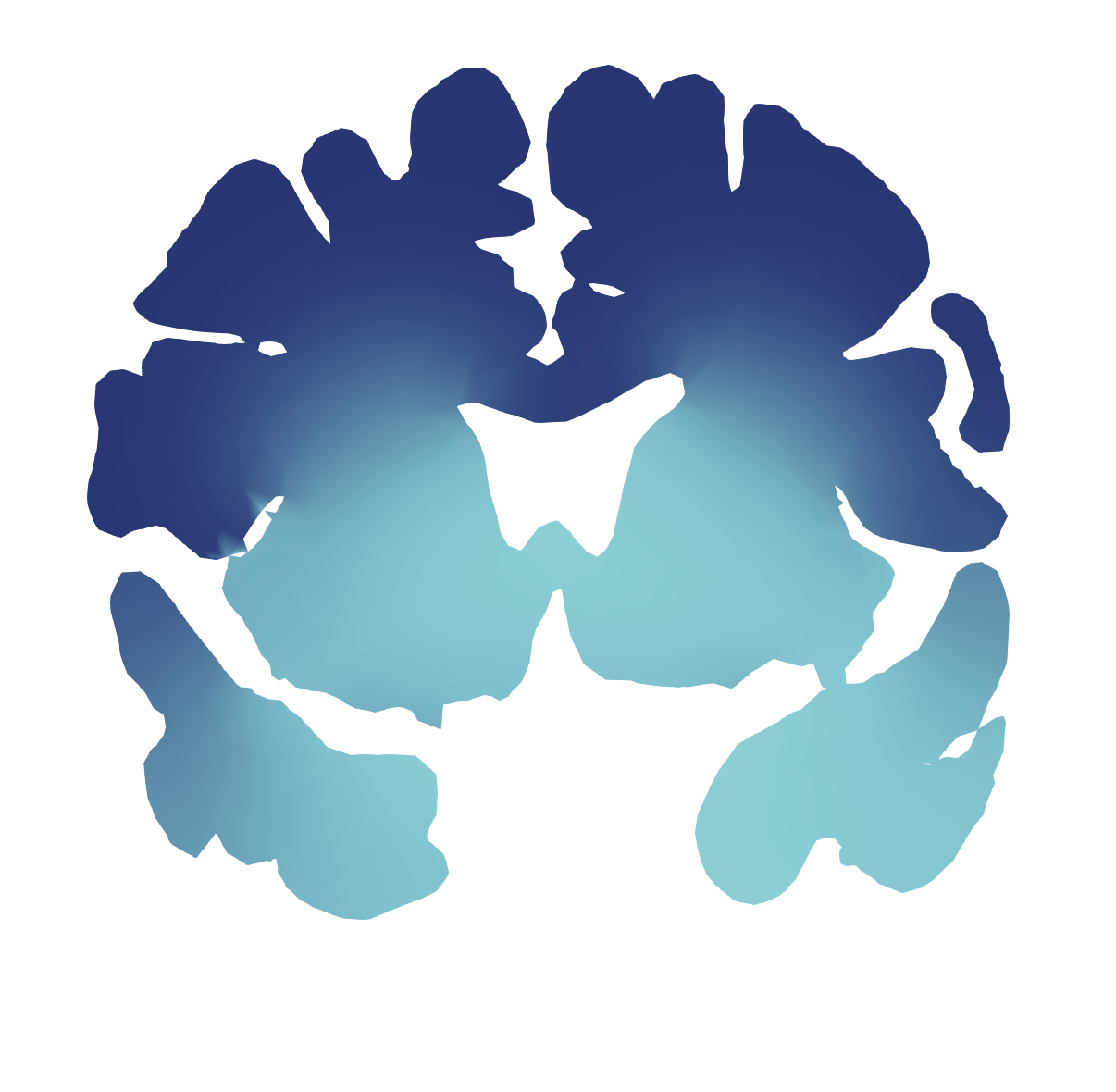}
            \includegraphics[width=0.32\textwidth]{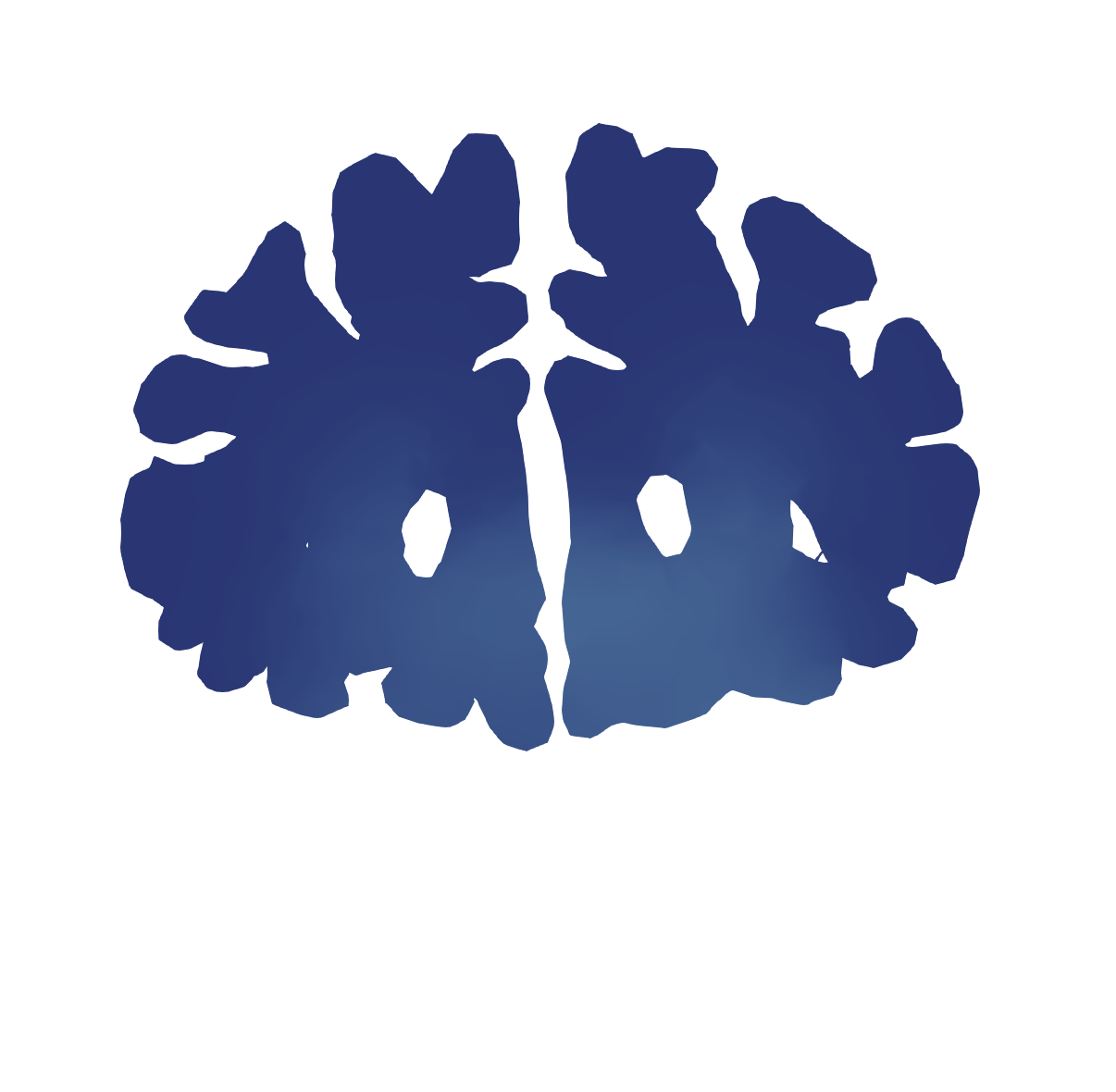}
            \includegraphics[width=0.32\textwidth]{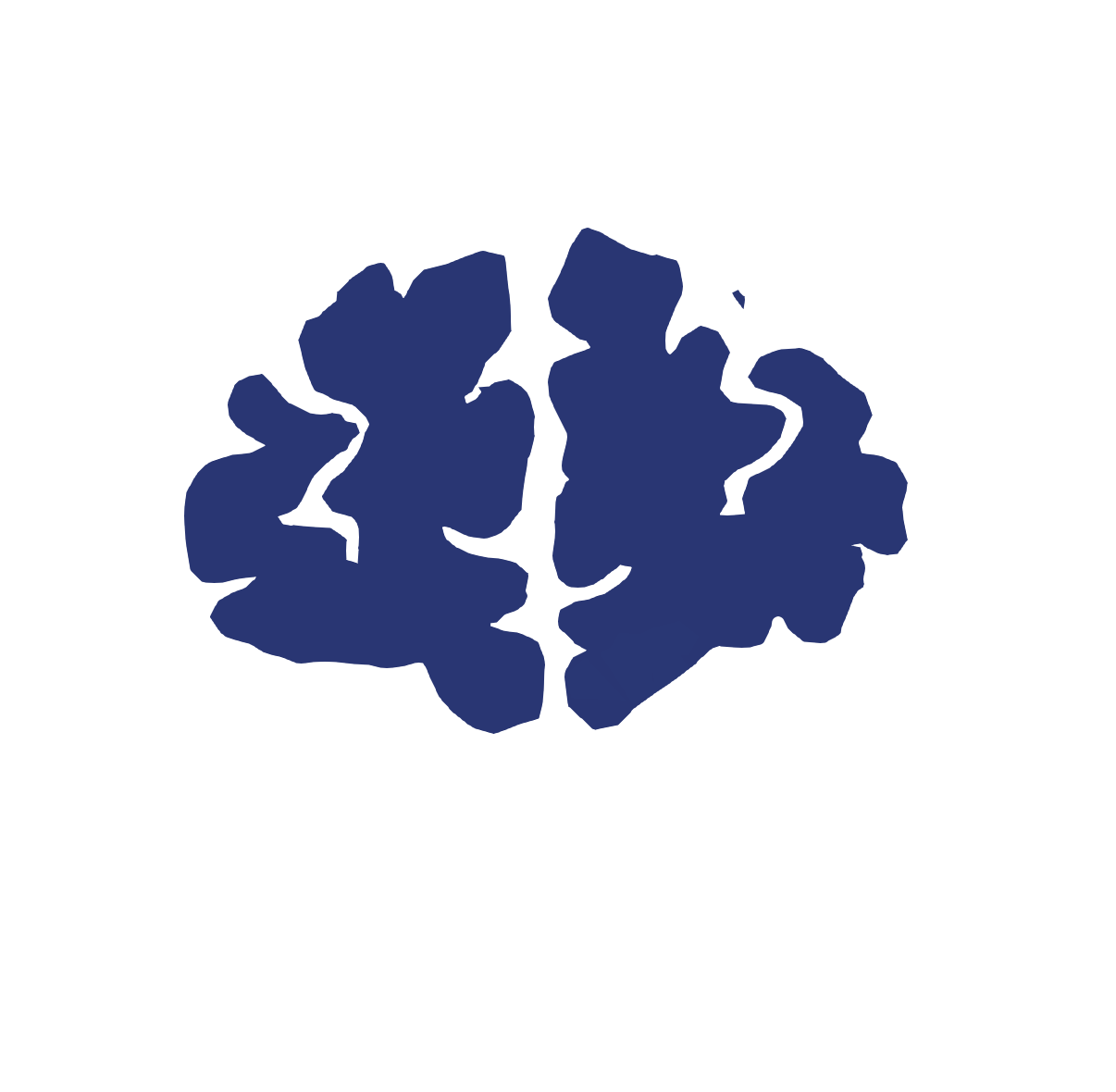}
            \caption{$t=20$ years}
            \label{fig:res3b}
        \end{subfigure}
        \begin{subfigure}{0.12\textwidth}
            \includegraphics[height=4cm]{final_simulation/scale_c.png}
        \end{subfigure}
        \caption{Test case of Section \ref{sec:results}. Evaluation of healthy protein loss at 10 and 20 years, observed across six coronal sections within the three-dimensional domain shown on the left.}
        \label{fig:res3}
    \end{figure}
    \begin{figure}[H]
        \centering
        \begin{subfigure}{0.12\textwidth}
            \includegraphics[width=0.7\textwidth]{final_simulation/sections_f.png}
        \end{subfigure}
        \begin{subfigure}{0.33\textwidth}
            \includegraphics[width=0.32\textwidth]{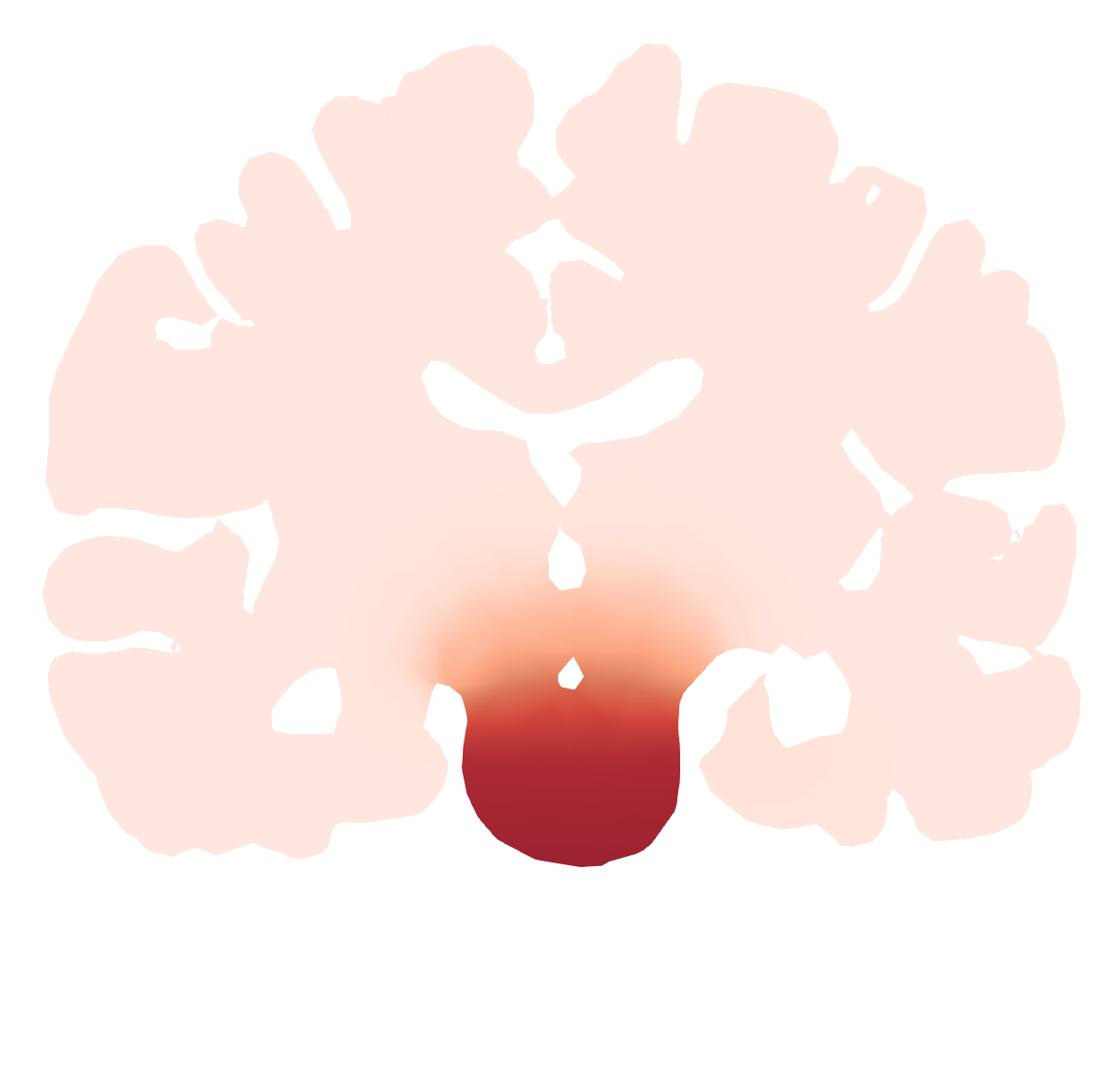}
            \includegraphics[width=0.32\textwidth]{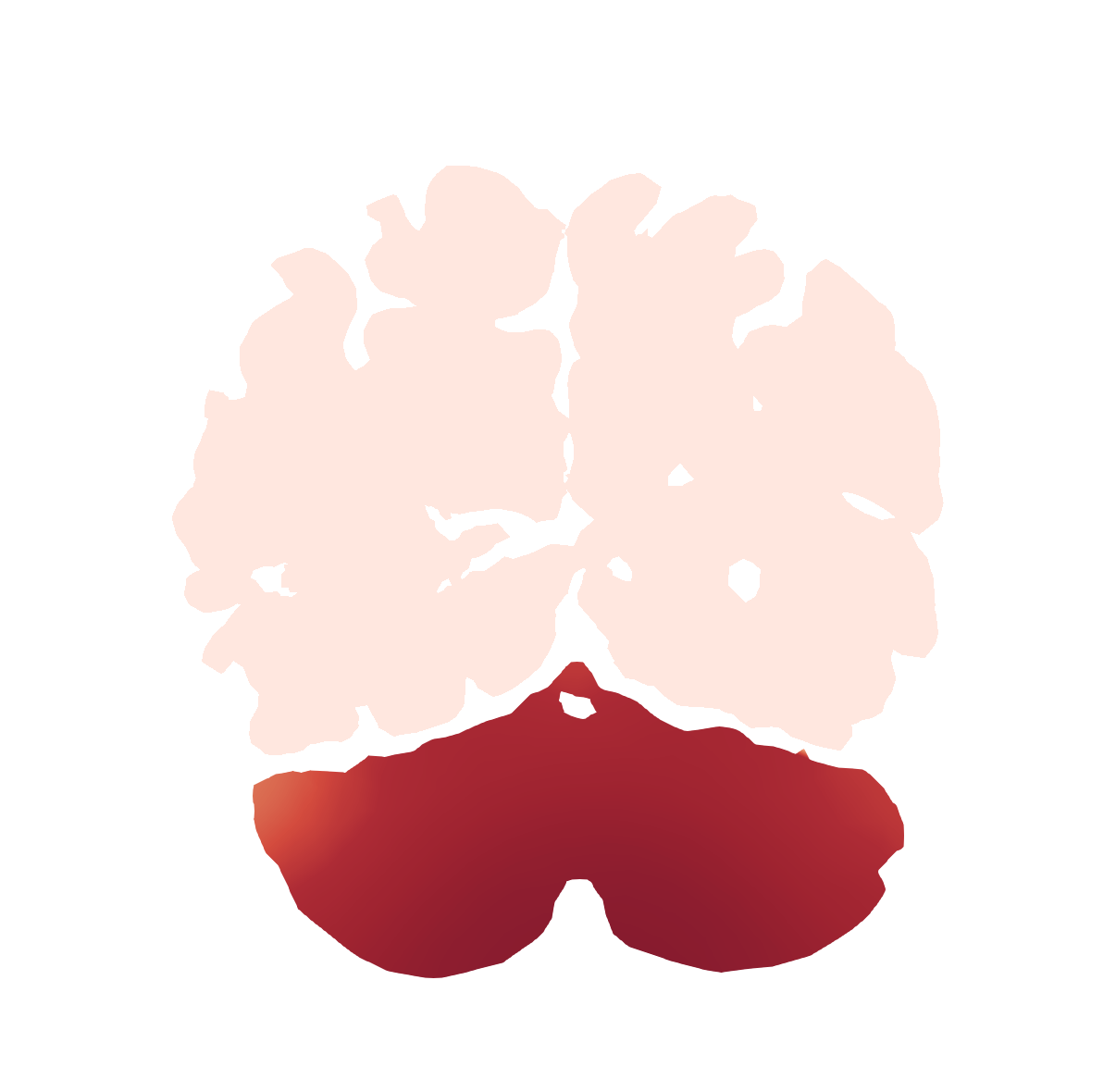}
            \includegraphics[width=0.32\textwidth]{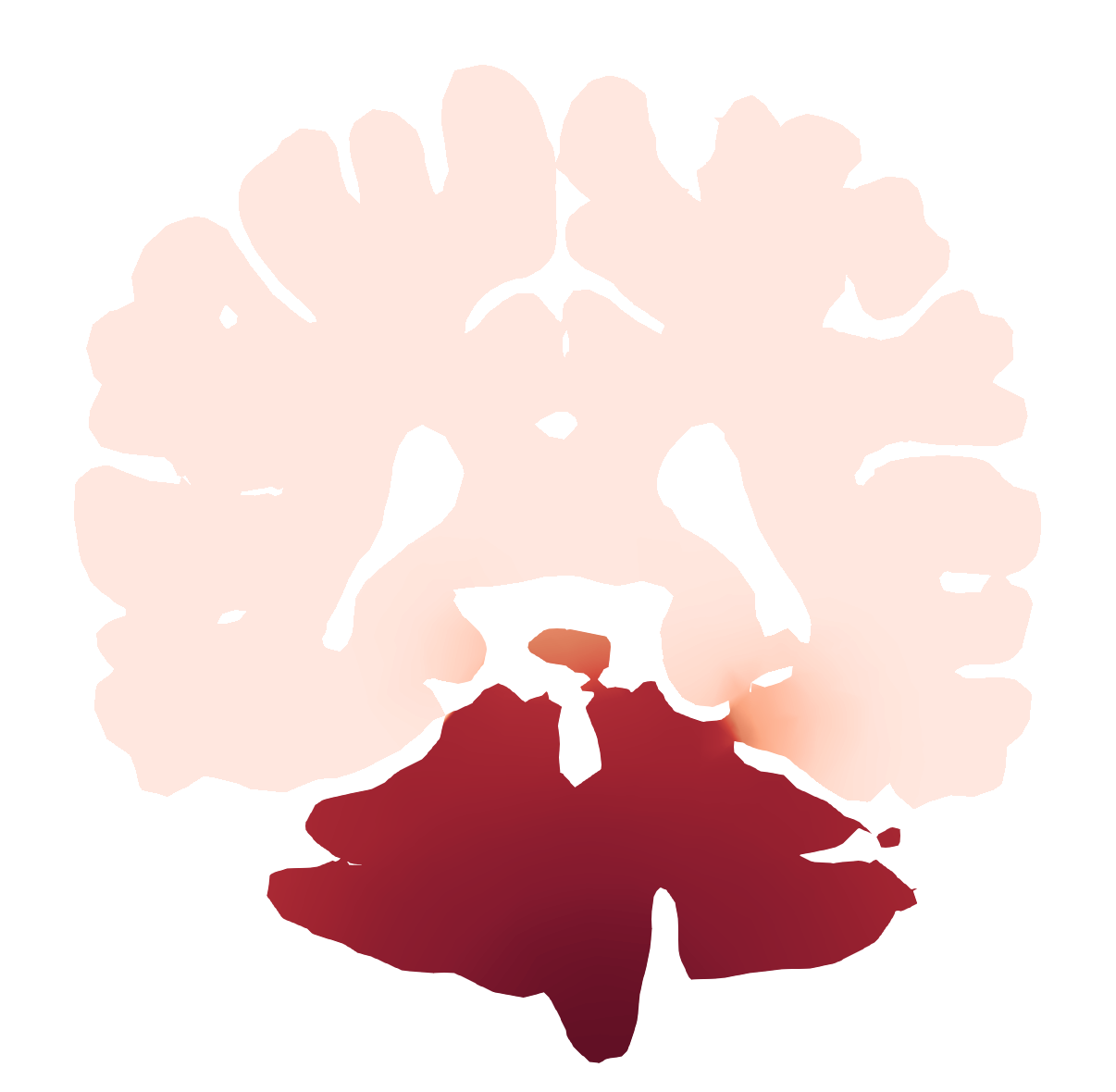}\\
            \includegraphics[width=0.32\textwidth]{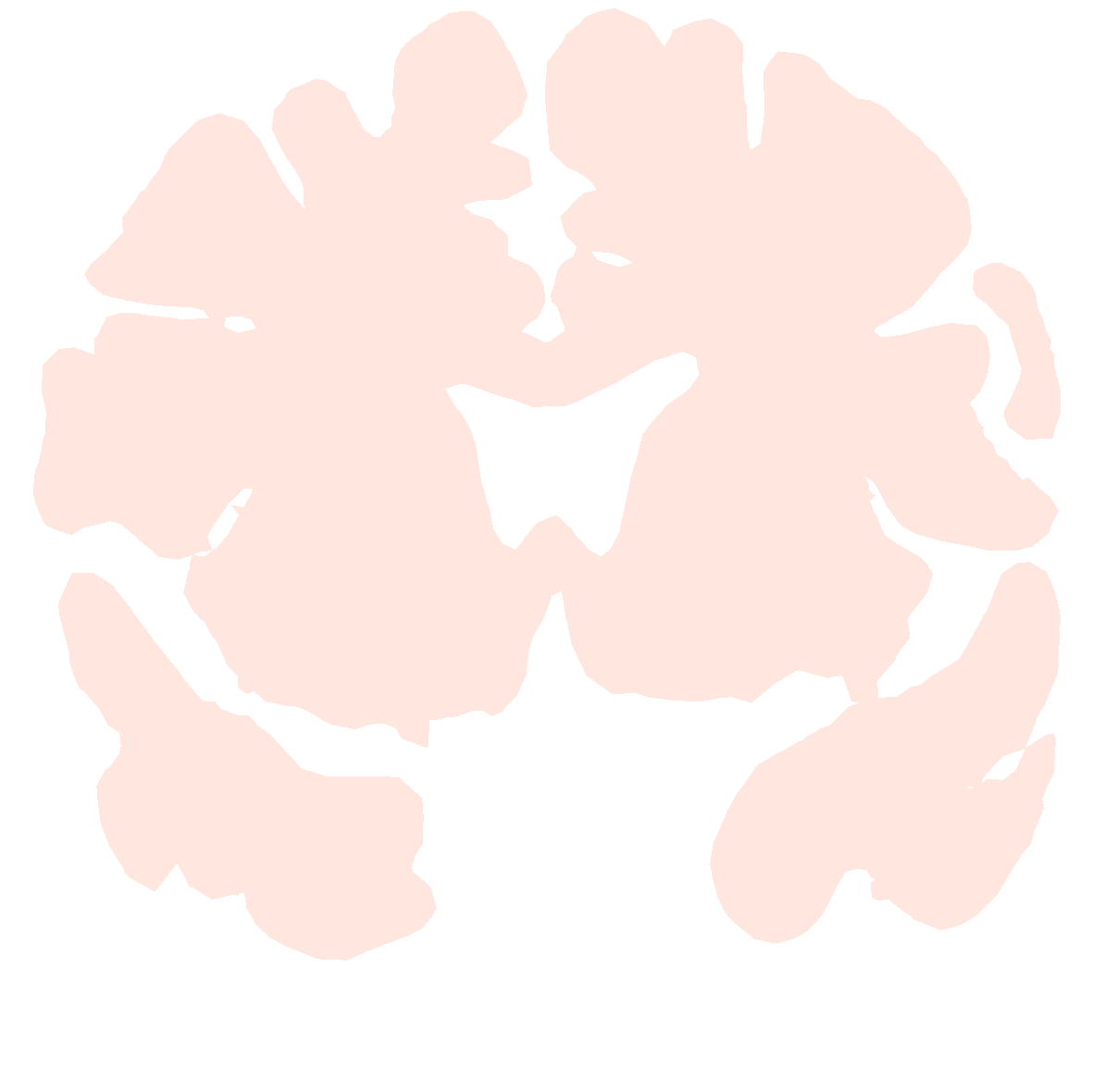}
            \includegraphics[width=0.32\textwidth]{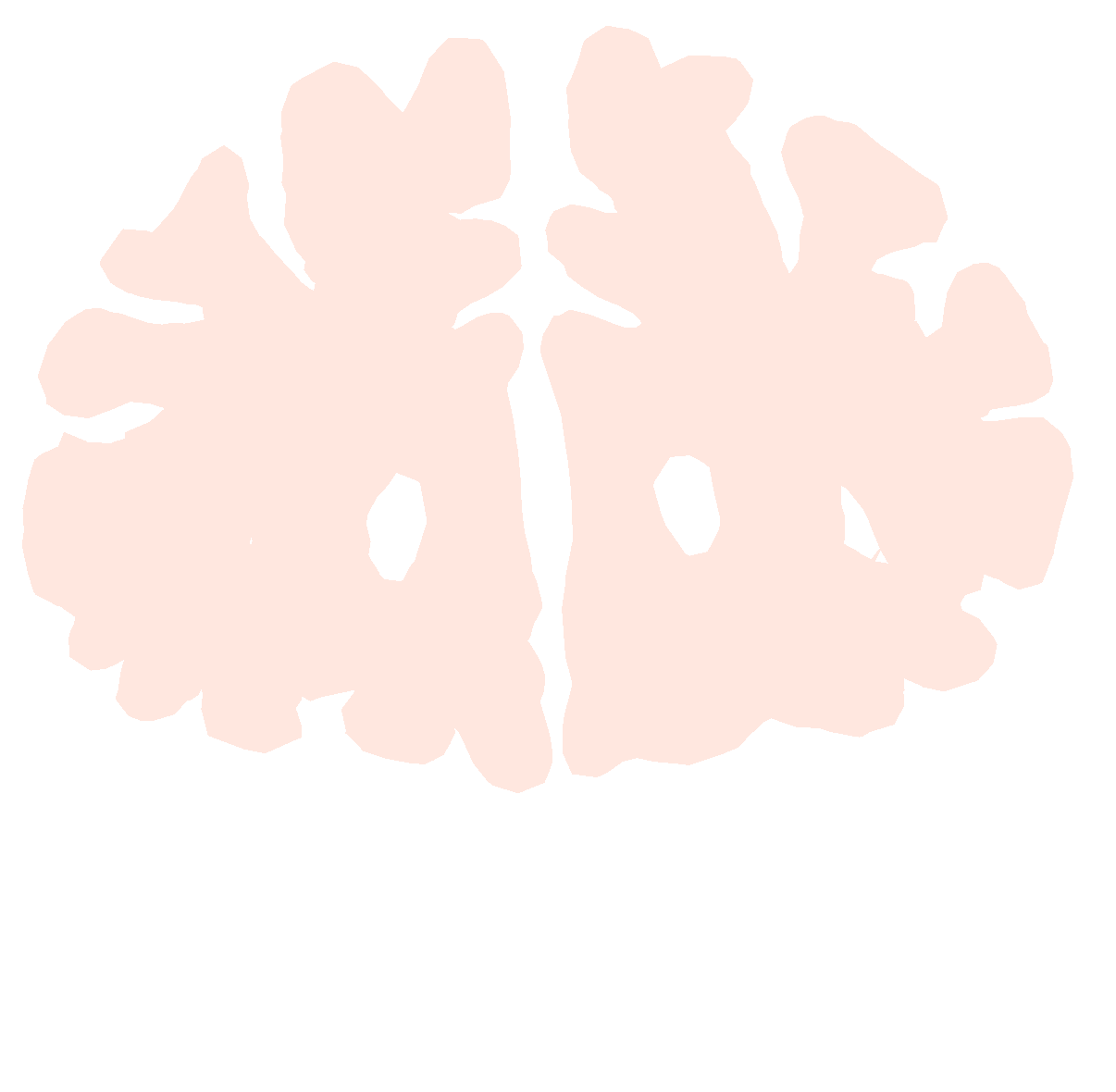}
            \includegraphics[width=0.32\textwidth]{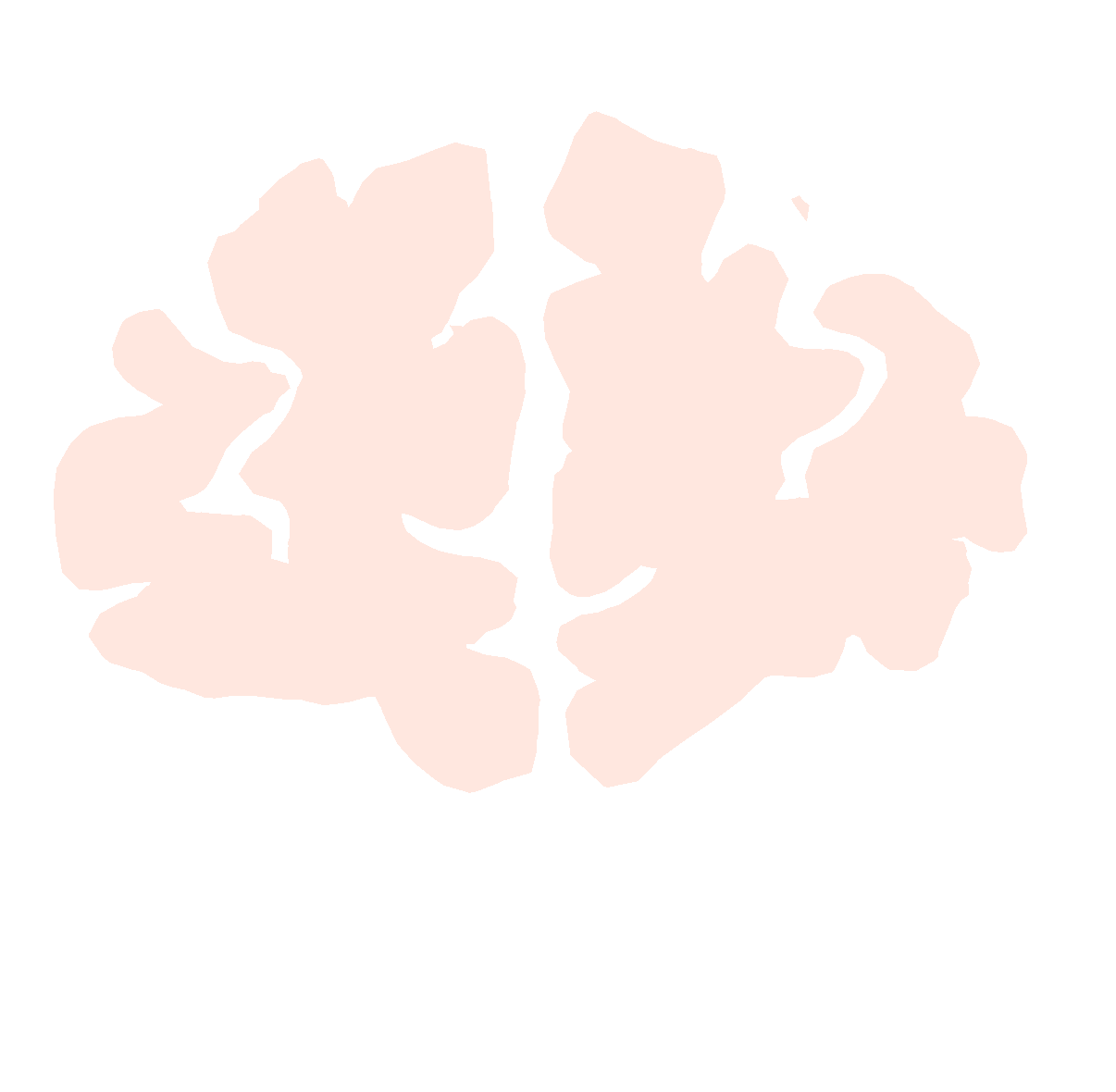}
            \caption{$t=10$ years}
            \label{fig:res4a}
        \end{subfigure}
        \hspace{2em}
        \begin{subfigure}{0.33\textwidth}
            \includegraphics[width=0.32\textwidth]{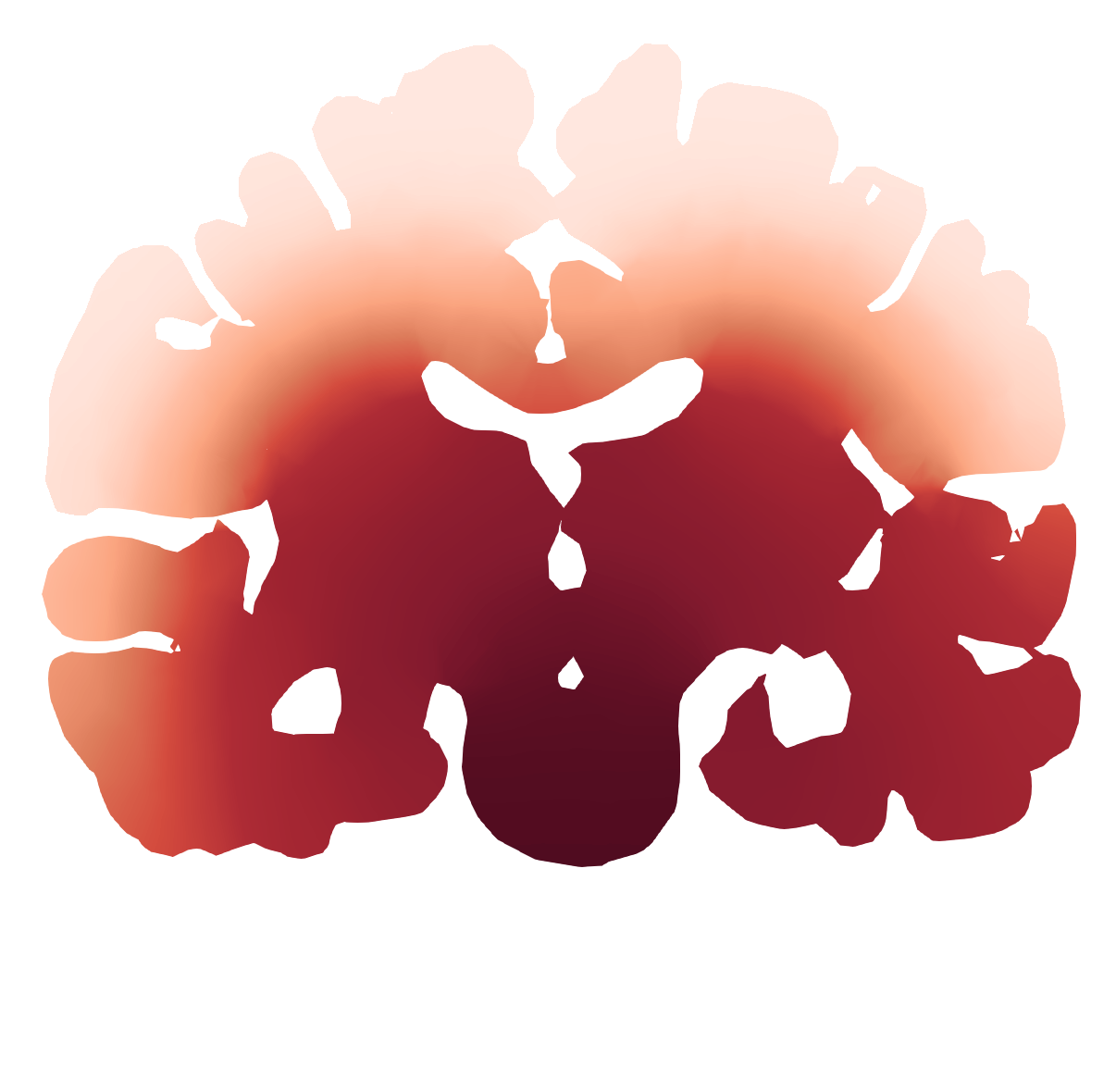}
            \includegraphics[width=0.32\textwidth]{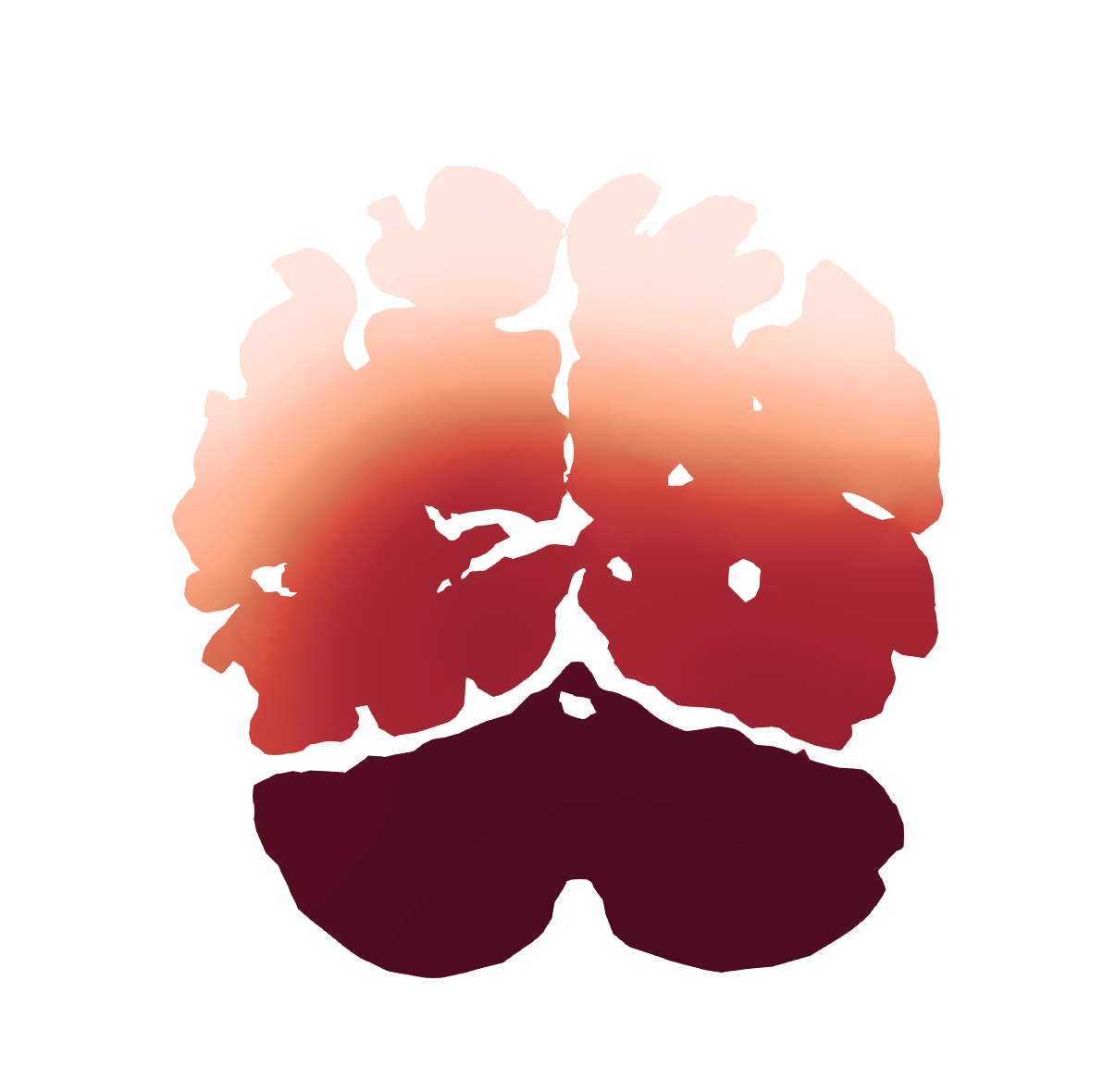}
            \includegraphics[width=0.32\textwidth]{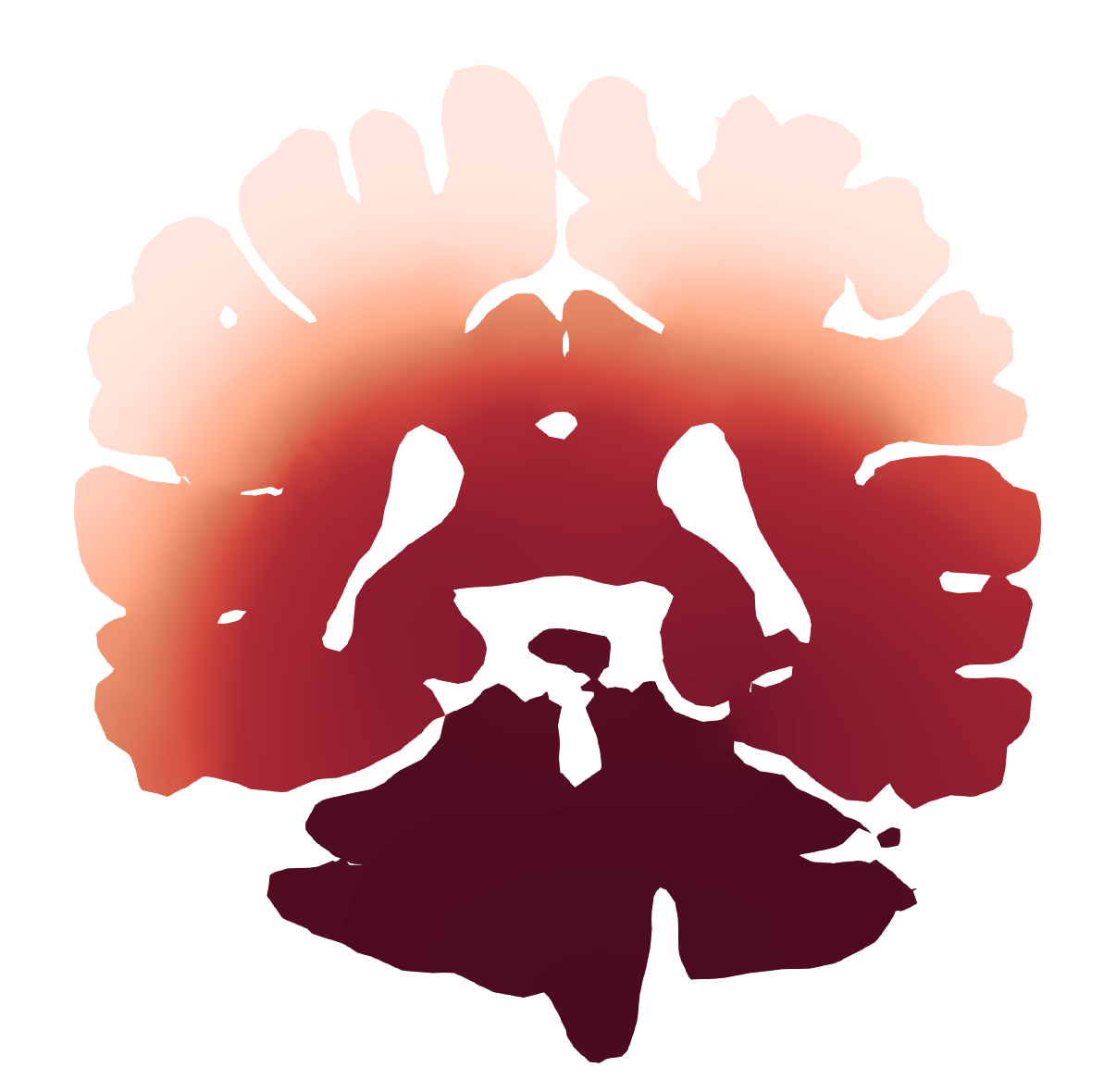}\\
            \includegraphics[width=0.32\textwidth]{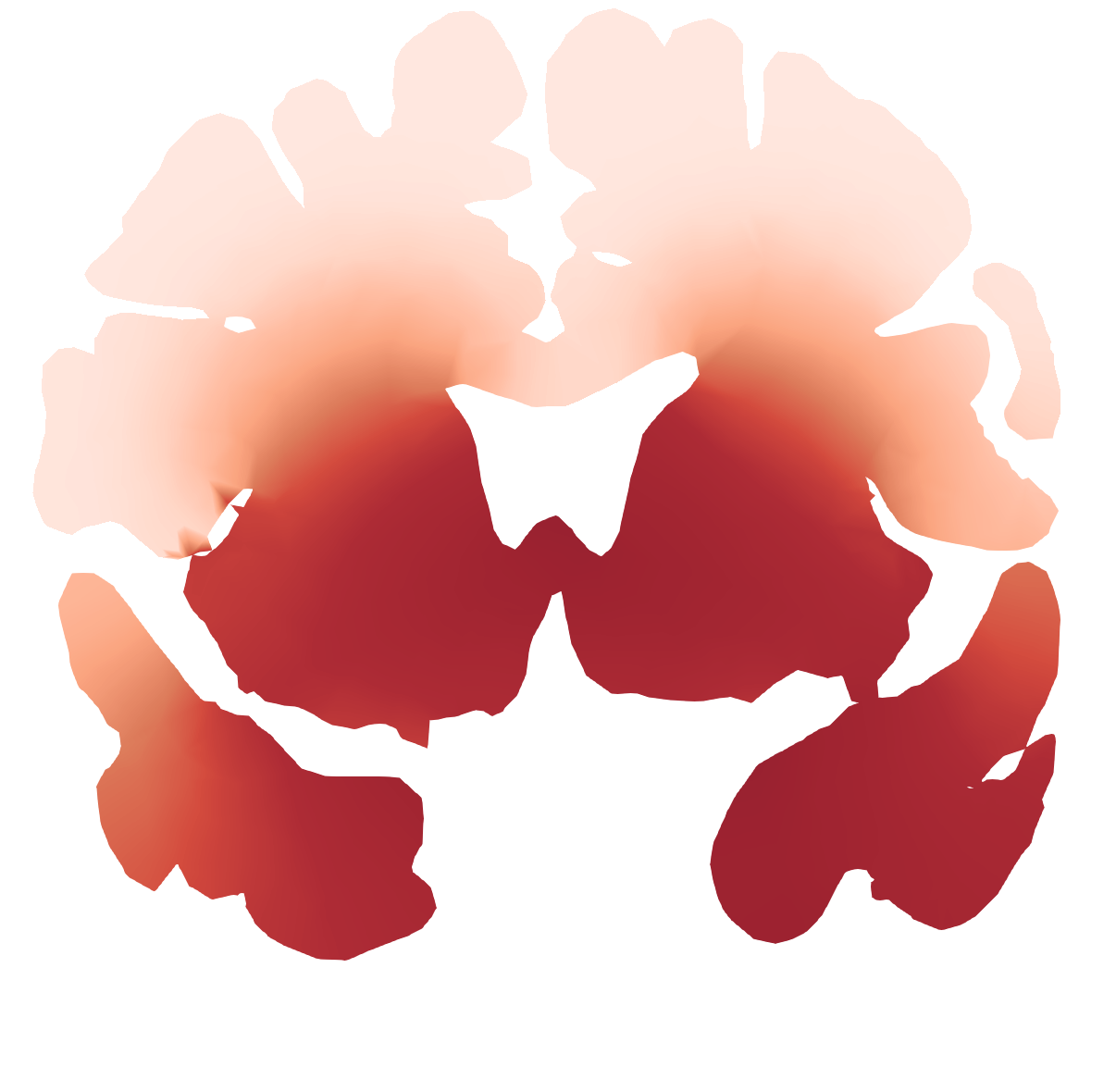}
            \includegraphics[width=0.32\textwidth]{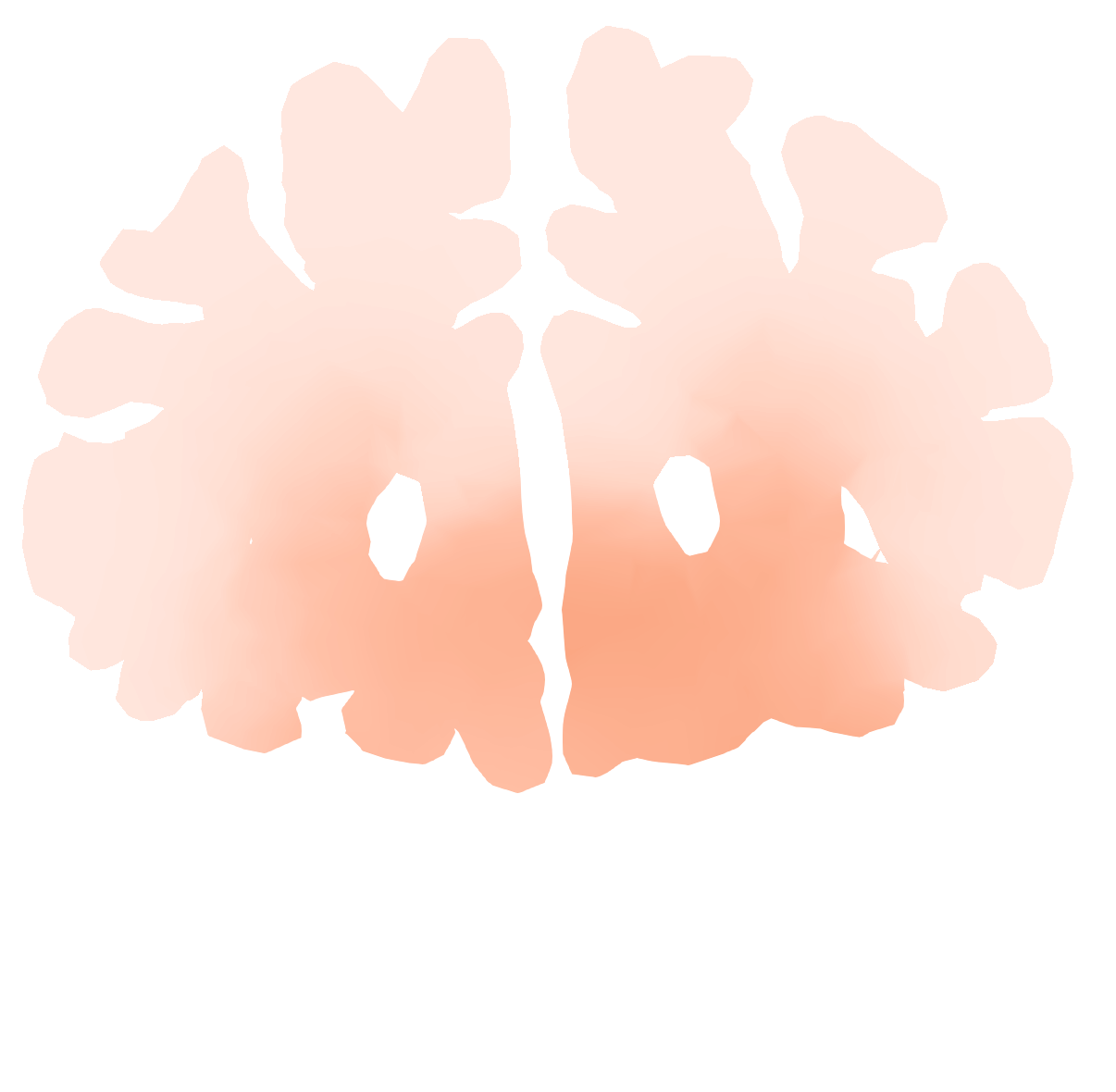}
            \includegraphics[width=0.32\textwidth]{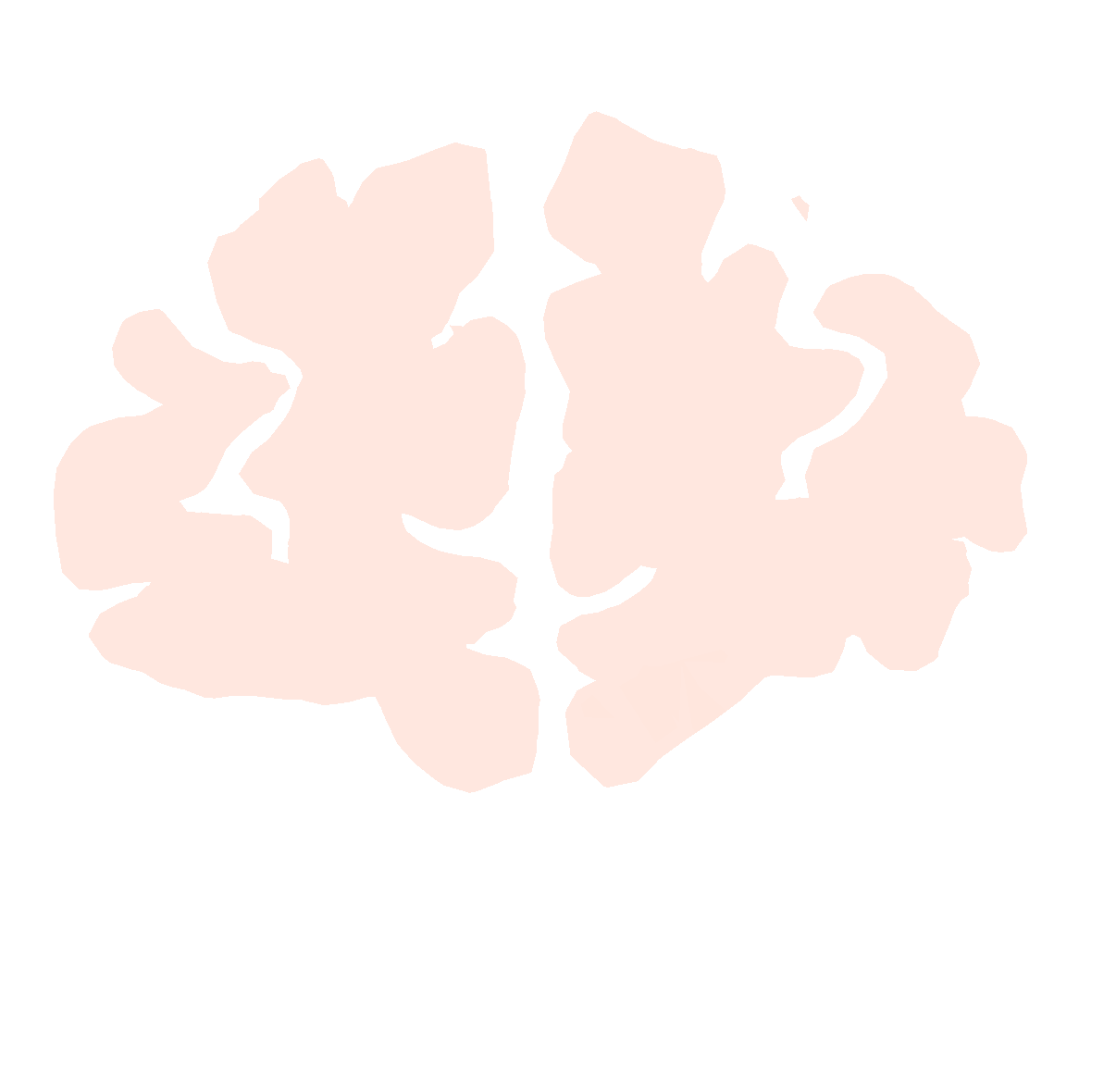}
            \caption{$t=20$ years}
            \label{fig:res4b}
        \end{subfigure}
        \begin{subfigure}{0.12\textwidth}
            \includegraphics[height=4cm]{final_simulation/scale_q.png}
        \end{subfigure}
        \caption{Test case of Section \ref{sec:results}. Evaluation of misfolded protein diffusion at 10 and 20 years, observed across six coronal sections within the three-dimensional domain shown on the left.}
        \label{fig:res4}
    \end{figure}
    }
    \restoregeometry
    
    \section{Conclusions}
    \label{sec:conclusions}
    In this study, we applied the heterodimer model to the study of the diffusion of the $\alpha$-synuclein protein throughout the brain in Parkinson's disease. We discussed some theoretical results concerning the heterodimer model and we derived its polytopal discontinuous Galerkin formulation. Convergence tests in two and three dimensions confirm the theoretical results for implicit and semi-implicit treatments of the nonlinear term in trivial domains. Furthermore, we conducted a convergence analysis with respect to the mesh size $h$ and the order of approximation $p$. Lastly, we performed a realistic simulation of $\alpha$-synuclein protein propagation across a three-dimensional brain mesh obtained from structural Magnetic Resonance Imaging. Our results have been validated by comparison with clinical data from the literature, showing good agreement. Additionally, biomarker curves correspond well to available clinical data. The choice of using the polytopal discontinuous Galerkin method has been justified by both the high rate of convergence achieved in numerical tests and the effective description of propagating fronts in the three-dimensional domain. Indeed, our results align well with the postmortem histological analyses and Braak's staging theory. The simulation also underscored the importance of detailed spatial simulations and the limitations of two-dimensional ones on brain sections. Moreover, it represents a step forward compared to the Fisher-Kolmogorov model, offering insights into the concentration of healthy proteins within the brain during the pathology.
    \par
    The computational framework enables the study of potential biomarkers for neurodegeneration not in vivo. This approach could lead to benefits in terms of early intervention in therapy. The mechanistic description of the phenomenon enables examination of the individual contributions' relative significance, facilitating the customization of therapy approaches. Furthermore, the complete spatial description also allows studying a priori the timing and the location of interventions. In the future, it is reasonable to envision patient-specific therapy based on quantitative models for prion-like diseases.
    
    \section*{Acknowledgements}
    The brain MRI images were provided by OASIS-3: Longitudinal Multimodal Neuroimaging: Principal Investigators: T. Benzinger, D. Marcus, J. Morris; NIH P30 AG066444, P50 AG00561, P30 NS09857781, P01 AG026276, P01 AG003991, R01 AG043434, UL1 TR000448, R01 EB009352. AV-45 doses were provided by Avid Radiopharmaceuticals, a wholly-owned subsidiary of Eli Lilly.
    
    \section*{Funding}
    This work received funding from the European Union (ERC SyG, NEMESIS, project number 101115663). Views and opinions expressed are however those of the authors only and do not necessarily reflect those of the European Union or the European Research Council Executive Agency. Neither the European Union nor the granting authority can be held responsible for them. P.F.A. has been partially funded by the research project PRIN 2020 (n. 20204LN5N5) funded by the Italian Ministry of University and Research (MUR).  P.F.A. and M.C. are members of the INdAM Research group GNCS.
    
    \section*{Declaration of competing interests}
    The authors declare that they have no known competing financial interests or personal relationships that could have appeared to influence the work reported in this article.
    
    \printbibliography
    
\end{document}